\pdfoutput=1
\RequirePackage{ifpdf}
\ifpdf 
\documentclass[pdftex]{sigma}
\else
\documentclass{sigma}
\fi

\usepackage{tikz,subcaption}
\usepackage[labelsep=period,labelfont=bf,font=small]{caption}

\usetikzlibrary{decorations.pathreplacing}
\usepackage{tikz-cd}

\numberwithin{equation}{section}

\newtheorem{Theorem}{Theorem}[section]
\newtheorem*{Theorem*}{Theorem}
\newtheorem{Corollary}[Theorem]{Corollary}
\newtheorem{Lemma}[Theorem]{Lemma}
\newtheorem{Proposition}[Theorem]{Proposition}

 { \theoremstyle{definition}
\newtheorem{Definition}[Theorem]{Definition}

\newtheorem{Remark}[Theorem]{Remark}

\newtheorem{Question}[Theorem]{Question}
}

\newenvironment{psmallmatrix}
 {\left(\begin{smallmatrix}}
 {\end{smallmatrix}\right)}

\definecolor{darkgreen}{rgb}{0,.5,0}

\colorlet{line1}{orange}
\colorlet{fill1}{line1!60}

\definecolor{line2}{rgb}{0,0,1}
\colorlet{fill2}{line2!60}

\colorlet{line3}{red}
\colorlet{fill3}{line3!60}

\colorlet{line4}{darkgreen!80}
\colorlet{fill4}{line4!60}

\newcommand{\triangleline}{dashed}

\def\C{{\mathbb C}}

\def\R{{\mathbb R}}

\def\Z{{\mathbb Z}}

\newcommand{\De}{\Delta}
\newcommand{\om}{\omega}
\newcommand{\area}{\operatorname{area}}
\newcommand{\PT}{\mathcal{P}_{\rm T}}
\newcommand{\PTtilde}{\widetilde{\PT}}
\newcommand{\PST}{\mathcal{P}_{\rm ST}}
\newcommand{\PSTtilde}{\widetilde{\mathcal{P}_{\rm ST}}}
\newcommand{\SL}{{\rm SL}_2(\Z)}
\newcommand{\SLn}{{\rm SL}_n(\Z)}

\newcommand{\stpoly}{[\De,\vec{c},\vec{\epsilon}\,]}
\newcommand{\strep}{(\De,\vec{c},\vec{\epsilon}\,)}

\newcommand{\vect}[2]{\begin{psmallmatrix} #1 \\ #2 \end{psmallmatrix}}

\newcommand{\vectbig}[2]{\begin{pmatrix} #1 \\ #2 \end{pmatrix}}

\begin{document}
\allowdisplaybreaks

\newcommand{\arXivNumber}{2210.06415}

\renewcommand{\PaperNumber}{081}

\FirstPageHeading

\ShortArticleName{Packing Densities of Delzant and Semitoric Polygons}

\ArticleName{Packing Densities of Delzant and Semitoric Polygons}

\Author{Yu DU~$^{\rm a}$, Gabriel KOSMACHER~$^{\rm a}$, Yichen LIU~$^{\rm a}$, Jeff MASSMAN~$^{\rm a}$, Joseph PALMER~$^{\rm ab}$,\newline Timothy THIEME~$^{\rm a}$, Jerry WU~$^{\rm a}$ and Zheyu ZHANG~$^{\rm a}$}

\AuthorNameForHeading{Y.~Du et al.}

\Address{$^{\rm a)}$~Department of Mathematics, University of Illinois at Urbana-Champaign,\\
\hphantom{$^{\rm a)}$}~Urbana, Illinois, USA}
\EmailD{\href{mailto:yudu3@illinois.edu}{yudu3@illinois.edu},
\href{mailto:gkk2@illinois.edu}{gkk2@illinois.edu},
\href{mailto:yichen23@illinois.edu}{yichen23@illinois.edu},
\newline
\hspace*{17.5mm}\href{mailto:massman4@illinois.edu}{massman4@illinois.edu}, \href{mailto:jpalmer5@illinois.edu}{jpalmer5@illinois.edu},
\href{mailto:tthiem3@illinois.edu}{tthiem3@illinois.edu},
\newline
\hspace*{17.5mm}\href{mailto:yihanwu4@illinois.edu}{yihanwu4@illinois.edu}, \href{mailto:zheyu3@illinois.edu}{zheyu3@illinois.edu}}

\Address{$^{\rm b)}$~Department of Mathematics, University of Antwerp,Antwerp, Belgium}

\ArticleDates{Received November 22, 2022, in final form October 20, 2023; Published online October 29, 2023}

\Abstract{Exploiting the relationship between $4$-dimensional toric and semitoric integrable systems with Delzant and semitoric polygons, respectively, we develop techniques to compute certain equivariant packing densities and equivariant capacities of these systems by working exclusively with the polygons. This expands on results of Pelayo and Pelayo--Schmidt. We compute the densities of several important examples and we also use our techniques to solve the equivariant semitoric perfect packing problem, i.e., we list all semitoric polygons for which the associated semitoric system admits an equivariant packing which fills all but a~set of measure zero of the manifold. This paper also serves as a concise and accessible introduction to Delzant and semitoric polygons in dimension four.}


\Classification{37J35; 53D20; 37J06}

\section{Introduction}

The space of possible states of a classical mechanical
system naturally comes with the structure of a so-called
\emph{symplectic manifold},
and for this reason the study of symplectic manifolds has been an active area of research for many years.
A dynamical system on
such a space is called an \emph{integrable system} if, roughly speaking, it has the maximal
possible number of quantities preserved by the dynamics (such as total energy or angular momentum).
In some highly symmetric situations, integrable systems can be represented by lower dimensional objects, such as polygons.
In this paper we will exploit this representation to
work entirely with certain $2$-dimensional polygons in order to
explore questions motivated by $4$-dimensional symplectic geometry.
Since we avoid having to work directly with any symplectic manifold, the bulk of this paper\footnote{With the notable exception of Section~\ref{sec:symplectic-background}, which explains the background and motivation from symplectic geometry.} is accessible to anyone with a background in linear algebra and multivariable calculus.
Though the required background is minimal, the proofs are non-trivial and we are able to develop an algorithm to explicitly calculate certain invariants of integrable systems, which are typically difficult to compute.

Since their discovery in the 1980s, \emph{symplectic capacities}
have been an important class of invariants in symplectic
geometry.
For a nice overview of this area see~\cite{CiHoLaSc2007} and the references therein. In~\cite{FiPaPe2016}, a~type of symplectic
invariants called \emph{equivariant capacities} are introduced for symplectic manifolds equipped with a symplectic group action. Integrable systems naturally come with such an action, and in~\cite{FiPaPe2016} it is shown that in some cases certain equivarant capacities of integrable systems can be computed directly from associated polygonal invariants of the systems, motivated by similar invariants studied earlier by Pelayo~\cite{Pe2006,Pe2007} and Pelayo--Schmidt~\cite{PeSc2008}.

\begin{figure}
\centering
 \begin{subfigure}{.45\linewidth}
 \centering
 \begin{tikzpicture}[scale = .9]
 \draw (0,0)--(0,3)--(3,3)--(6,0)--cycle;
 \draw [line1,line width = 1,fill = fill1] (0,1) -- (0,0) -- (1,0);
 \draw [\triangleline, line1,line width = 1] (0,1) -- (1,0);
 \draw [line2,line width = 1,fill = fill2] (1,3) -- (0,3) -- (0,2);
 \draw [\triangleline, line2,line width = 1] (0,2) -- (1,3);
 \draw [line3,line width = 1,fill = fill3] (5,0) -- (6,0) -- (5,1);
 \draw [\triangleline, line3,line width = 1] (5,1) -- (5,0);
 \draw [line4,line width = 1,fill = fill4] (1,3) -- (3,3) -- (5,1);
 \draw [\triangleline, line4,line width = 1] (5,1) -- (1,3);
 \end{tikzpicture}
 \caption{A toric packing}
 \label{fig:delz-pack-example}
 \end{subfigure}
 \begin{subfigure}{.45\linewidth}
 \centering
 \begin{tikzpicture}[scale = 1.7]
 \draw (0,0)--(2,1)--(4,0)--cycle;
 \draw[line1, line width = 1, fill=fill1] (1,0.5)--(0,0)--(0.5,0);
 \draw [\triangleline, line1,line width = 1] (0.5,0) -- (1,0.5);
 \fill[fill=fill2] (2.70,0)--(4,0)--(2,1)--(1.4,0.70)--(2,0.70);
 \draw[line2, line width = 1] (2.70,0)--(4,0)--(2,1)--(1.4,0.70);
 \draw [\triangleline, line2,line width = 1] (1.4,0.70) -- (2,0.70)--(2.70,0);
 \node at (2,0.3) {$\times$};
 \draw [dashed](2,1)--(2,0.3);
 \end{tikzpicture}
 \caption{A semitoric packing}
 \label{fig:st-pack-example}
 \end{subfigure}
 \caption{A toric packing
 of a Delzant polygon by four triangles (left) and a semitoric packing of a~semitoric polygon representative by two ``triangles'' (right). Neither packing is maximal. Notice that one of the packed triangles in the semitoric packing is deformed by the presence of the cut (the dotted vertical line), as we describe in Section~\ref{sec:semitoric}.}
 \label{fig:delz-st-pack-example}
\end{figure}
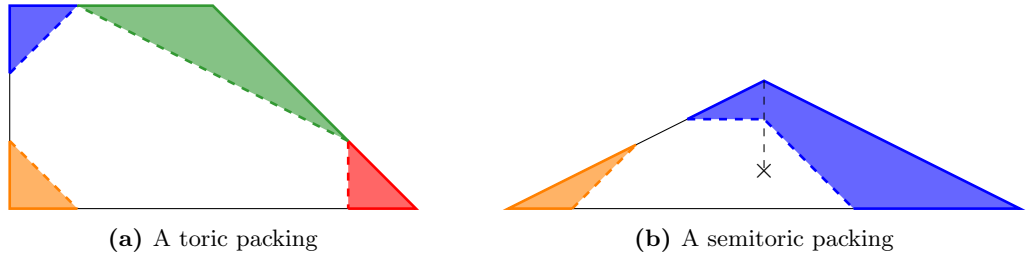

All of the invariants that we will be concerned with in this paper are \emph{packing-type invariants}. These invariants ask how much of the volume of a given manifold can be filled with disjoint images of embedded balls satisfying certain conditions related to the integrable system structure, and in this situation in fact this is equivalent to asking how much of a given Delzant or semitoric polygon can be filled with disjoint triangles satisfying certain conditions.
Such
a~packing of a~Delzant polygon is shown in Figure~\ref{fig:delz-pack-example} and such a packing of a~semitoric polygon is shown in Figure~\ref{fig:st-pack-example} (note that one of the triangles packed in the semitoric polygon is deformed by the ``cut'' present in the polygon).
This relationship between the packing of the manifold and the packing of the polygon was proven for toric integrable systems by Pelayo~\cite{Pe2006} and for semitoric integrable systems by Figalli--Pelayo--Palmer~\cite{FiPaPe2016}, but until now these equivariant capacities were not explicitly computed in many cases.
In fact, we will actually spend most of the paper computing \emph{packing densities}, which are the maximal proportion of a polygon (or manifold) which can be filled by admissibly packed triangles (or equivariantly embedded balls).
The capacities in question can be easily determined from the packing densities, and vice versa, via equations~\eqref{eqn:toric-cap-to-density} and~\eqref{eqn:semitoric-cap-to-density} given in Section~\ref{sec:capacities}.

The works of Pelayo~\cite{Pe2006, Pe2007} and Pelayo--Schmidt~\cite{PeSc2008} were to our knowledge the first papers in this direction, exploring torus-equivariant packing-type invariants of toric integrable systems. Inspired by these works, Figalli--Palmer--Pelayo~\cite{FiPaPe2016} defined the general notion of equivariant capacities and presented an analogous construction of such a capacity for semitoric systems. The present paper computes the toric packing invariants from~\cite{Pe2006, Pe2007, PeSc2008} and the semitoric packing invariants from~\cite{FiPaPe2016} in some new examples.
In order to perform these computations, we recover a result of Pelayo--Schmidt~\cite{PeSc2008} about how to find maximal toric packings and generalize it to the semitoric case.

Of course, the more difficult and subtle version of this question is the non-equivariant case: i.e.,~asking questions related to how large a ball may be symplectically embedded into a given symplectic manifold (as in the work of Gromov~\cite{Gromov1985}). One motivation for studying the equivariant versions of these questions, as we do in this paper, is that equivariant results naturally provide bounds for the non-equivariant questions. One way to think of it is that combinatorial techniques, like those in this paper, can provide constructions of explicit examples of symplectically embedded objects (or families of objects), which therefore gives a lower bound on the largest such objects which can be symplectically embedded (forgetting about the equivariant structure).
In some cases it can be shown that the bound obtained from equivariant techniques is sharp, such as in~\cite{HS17}.
Moreover, obtaining such bounds is similar to the point of view taken in several recent works~\cite{CGHMP,CV,four-periodic,macgill,MMW-staircasepatterns} in which a sequence of bounds is obtained by constructing a sequence of different systems on a given symplectic manifold. We discuss these papers in more detail in Remarks~\ref{rmk:ATF}, and~\ref{rmk:S1-actions-sharp}, and Section~\ref{sec:ATFs}.

The purpose of this paper is two-fold. Firstly, we compute these packing densities in as many cases as possible, and describe strategies to compute them in further examples.
Secondly, this paper also serves as a concise introduction to Delzant and semitoric polygons, carefully describing these objects and giving many examples. It is our goal to provide a resource for researchers from outside the semitoric community to quickly and efficiently learn about these polygons, especially semitoric polygons.

Specifically, in this paper we
\begin{enumerate}\itemsep=0pt
 \item[(1)] provide a thorough introduction to Delzant and semitoric polygons (see Section~\ref{sec:polygons-intro});
 \item[(2)] build a framework to study and compute the packing densities of Delzant and semitoric polygons (see Proposition~\ref{prop:toric} and Theorem~\ref{thm:semitoric});
 \item[(3)] explicitly compute the packing density in many important cases (see Theorems~\ref{thm:toric-examples} and~\ref{thm:semitoric-examples});
 \item[(4)] produce an exhaustive list of all semitoric polygons which admit a ``perfect packing'', which fills the polygon up to a set of measure zero (see Theorem~\ref{thm:perfect}).
\end{enumerate}
In Section~\ref{sec:symplectic-background}, we explain
the motivation for this research from the point of view of symplectic geometry, and, in particular, in Section~\ref{sec:capacities}, we interpret our main results into that context, but outside Section~\ref{sec:symplectic-background} no prior knowledge of symplectic geometry (or indeed any differential geometry at all) will be necessary.
Furthermore, a reader with no background in symplectic geometry can easily skip Section~\ref{sec:symplectic-background} and still understand the rest of the paper.
In particular, in Section~\ref{sec:polygons-intro}, we introduce Delzant polygons and semitoric polygons independently of their relationship to symplectic geometry.

\subsection{Motivation}
Toric integrable systems
were classified in terms of a polytope in the 1980s by the results of Atiyah~\cite{Atiyah}, Guillemin--Sternberg~\cite{GuSt1982}, and Delzant~\cite{De1988},
and about 15 years ago this classification was extended to semitoric integrable systems
by Pelayo--V\~{u} Ng\d{o}c~\cite{PeVN2009,PeVN2011} by including a more complicated polygon invariant and also several other invariants which do not appear in the toric case.

The details of these classifications, and indeed the definitions of toric and semitoric integrable systems, are not strictly necessary to study certain invariants of these systems.
This is because, due to the results of Pelayo~\cite{Pe2006, Pe2007}, Pelayo--Schmidt~\cite{PeSc2008}, and Figalli--Palmer--Pelayo~\cite{FiPaPe2016}, the invariants we are interested in can be computed \emph{directly from the associated polygons}.
This is the perspective that we will take in this paper, completely circumventing the symplectic geometry, and allowing us to use careful application of relatively elementary techniques to obtain very non-trivial results.

We will be working with two classes of polygons: Delzant polygons and semitoric polygons. Delzant polygons are a specific class of convex polygons, but semitoric polygons are more complicated. They include the extra data of ``marked points'' and ``cuts'' which we will explain in Section~\ref{sec:semitoric-poly},
and a ``semitoric polygon'' technically refers to an equivalence class of polygons related by the application of certain operations, such as changing the direction of the cut.
Semitoric polygons were introduced by V\~{u} Ng\d{o}c~\cite{VN2007} and are a special case of the almost toric base diagrams introduced by Symington~\cite{Sy2003}.

Semitoric polygons have been studied in many contexts: they have been computed for specific systems~\cite{ADH-momenta,HohMeu,HoPa2017,LFPal, LFP}, compared with invariants of $S^1$-spaces~\cite{HP-extend, HSS2015}, extended to more general types of systems~\cite{PPT-nonsimple, PeRaVN2015},
used to define a metric on the space of semitoric systems~\cite{PaSTMetric2015},
used to find minimal models of semitoric systems~\cite{KPP2015,KPP2018}, and used to study the inverse question regarding quantum integrable systems~\cite{LFPeVN2016,LFVN2021}, among many other applications.
In the present paper we define them as objects of interest in their own right, independent of the underlying symplectic geometry. We do this to provide a resource for readers interested in studying these objects that may not possess the background in symplectic geometry which is typically necessary.

Motivated by trying to understand packing toric integrable systems by equivariantly embedded balls, Pelayo~\cite{Pe2006, Pe2007} determined a certain set of rules for how to pack a Delzant polygon $\Delta$ by triangles.
Roughly, the triangles must be disjoint and each be nested inside a corner of $\Delta$, satisfying certain conditions on the relative lengths of their edges, see in Figure~\ref{fig:delz-st-pack-example}.
Any set of triangles which follow
these rules is called
an \emph{admissible packing of $\Delta$}.
In Figalli--Palmer--Pelayo~\cite{FiPaPe2016}, these rules were extended to obtain rules for admissible packings of semitoric polygons, taking into account the marked points and cuts.
We will explain these rules later in the paper, but now we can state our main question:

\begin{Question}
 Given a Delzant or semitoric polygon, what is the maximal fraction of the area of the polygon which can be filled by an admissible packing of triangles?
\end{Question}

By examining the problem carefully, in each case we are able to transform the question into a question about maximizing the magnitude of a given vector in $\mathbb{R}^d$ inside a convex compact polytope defined by certain inequalities, where $d$ is the number of vertices of the given Delzant or semitoric polygon.
The case of Delzant polygons was already understood by Pelayo--Schmidt~\cite{PeSc2008}, but the case of semitoric polygons is completely new.
We show that the maximum must occur on one of the vertices of this $d$-dimensional polytope, and we describe an algorithm to use a~computer to exhaustively check all such vertices.
Note that while we are using a computer to compute the maximal density, it is \emph{not} an approximation, since it is simply comparing the finitely-many candidates for the maximal packing.
Moreover, in several important cases we use our techniques to find the maximal packing by hand explicitly, and thus obtain the
maximal packing density for many important Delzant and semitoric polygons.

\begin{Remark}\label{rmk:ATF}
Semitoric systems are a special case of \emph{almost-toric fibrations}~\cite{Sy2003}.
In Section~\ref{sec:ATFs}, we discuss the relationship between the techniques used in the present paper and certain techniques recently employed by various authors to use almost toric fibrations to obtain lower-bounds for ellipsoid packing capacities. In particular, such lower bounds can be used to study so-called ``infinite staircase'' behavior, as in McDuff--Schlenk~\cite{MS12}.
\end{Remark}

\begin{Remark}\label{rmk:S1-actions-sharp}
In \cite{HS17}, the authors used equivariant techniques to obtain an upper bound, which in some cases is sharp, for the Gromov width of a Fano symplectic manifold with a semi-free $S^1$-action.
The Gromov width measures the maximal radius of a single ball which can be symplectically embedded, and in the present paper we are concerned with equivariantly symplectically embedding multiple balls.
Since both toric and semitoric systems admit a global $S^1$-action, it may be possible to adapt the technique used in \cite{HS17} in order to give a criterion on when the lower bound obtained by our equivariant estimation of packing by multiple balls is sharp.
\end{Remark}

\subsection{Results}

Given any line segment in $\R^2$ which is either vertical or has rational slope, a number $\ell\geq 0$ called its \emph{$\SL$-length} can be defined, see Definition~\ref{def:sl2z-length}.
The triangles that appear in toric and semitoric packings are \emph{$\SL$-equilateral triangles}, which means the $\SL$-length of all three of their sides are equal. We define the \emph{$\SL$-size} of such a triangle to be the common $\SL$-length of its sides. For instance, for any $\lambda>0$ we say that the triangle with vertices at $(0,0)$, $(\lambda,0)$, $(0,\lambda)$ has $\SL$-size $\lambda$.

\subsubsection{Toric packing}
Let $\De$ be a Delzant polygon with $d$ vertices.
We label the vertices in clockwise order, starting with the vertex with lexicographically minimal coordinates,\footnote{That is, we take the vertex which is minimal under the ordering given by $(x,y)<(x',y')$ if and only if either~$x<x'$ or both $x=x'$ and $y< y'$. This is the ``bottom left'' vertex of the polygon.} as $p_1,\ldots,p_d$.
Let $\ell_i$ denote the $\SL$-length of the edge connecting $p_i$ to $p_{i+1}$, taking $p_{d+1} = p_1$ by convention.
Let $\PT(\De)$ denote the set of all admissible packings of $\De$.
For any $P\in\PT(\De)$, we let $\textrm{area}(P)$ denote the sum of the areas of the triangles in the packing.
We define a map
 \begin{align} \Lambda \colon\ \PT(\De) & \to \R^d, \label{eqn:Lambda} \\ \bigcup_{i=1}^d B_i & \mapsto (\lambda_1,\ldots, \lambda_d),\nonumber
\end{align}
 where $\lambda_i$ is the $\SL$-size of the triangle $B_i$, which is packed at the vertex $p_i$ (taking $\lambda_i=0$ if there is no triangle packed at that vertex).
Define
\begin{equation}\label{eqn:tildePT}
 \PTtilde(\De) = \big\{(\lambda_1,\ldots,\lambda_d)\in\R^d \mid \lambda_i\geq 0 \text{ and } \lambda_i+\lambda_{i+1}\leq \ell_i \text{ for all }i\in\{1,\ldots,d\}\big\}.
\end{equation}

\begin{Proposition}\label{prop:toric}
 Let $\De$ be any Delzant polygon. Then
 \[
 \Lambda(\PT(\De)) = \PTtilde(\De)
 \]
 and $\Lambda$ is a bijection from $\PT(\De)$ to the compact convex polytope $\PTtilde(\De)$.
 Furthermore,
 \begin{equation}\label{eqn:LambdaP-magnitude}
 \area(P) = \frac{1}{2}||\Lambda(P)||^2
 \end{equation}
 for any $P\in\PT(\De)$.
 \end{Proposition}

Proposition~\ref{prop:toric} is proved in Section~\ref{sec:toric-computing}.
This result already appeared in the work of Pelayo--Schmidt~\cite{PeSc2008}, but we include a complete proof nonetheless, since we use an argument which we will also extend to the semitoric case.
Note that the map $\Lambda$ depends on the labeling we chose of the vertices $p_1,\ldots, p_d$, but in the end we are most interested in the magnitude of the resulting vectors in $\R^d$, which is independent of reordering the coordinates, so any clockwise or counter-clockwise labeling will suffice.

Due to Proposition~\ref{prop:toric}, to study the packings of $\De$ we can instead study the polytope ${\PTtilde(\Delta)\subset\R^d}$.
In particular, finding a maximal packing of $\De$ is equivalent to finding a vector in $\PTtilde(\Delta)$ of maximal magnitude.
This already makes the problem much more tractable.
Furthermore, the magnitude of vectors in a compact polytope in $\R^d$ is maximized on one of its vertices (this is well known, and we prove it explicitly in Proposition~\ref{prop:vert}), and therefore the problem of finding a maximal packing can be reduced to checking the finitely many vertices of $\PTtilde(\De)$, which can be performed exactly by a computer.
We were thus able to design an algorithm to find the maximal packing of a given Delzant polygon $\Delta$.
An implementation of this algorithm in Python is available online, see~\cite{alg-url}.

Even without the use of a computer, Proposition~\ref{prop:toric} reduces the problem of finding a maximal packing to the inspection of a short list of inequalities.
Making use of Proposition~\ref{prop:toric}, in Section~\ref{sec:toric-computing}, we compute the toric packing density of several important examples of Delzant polygons.
Combining equations~\eqref{eqn:pack-density-triangle}, \eqref{eqn:pack-density-rect}, and~\eqref{eqn:pack-density-hirz}, derived in that section, we have the following.

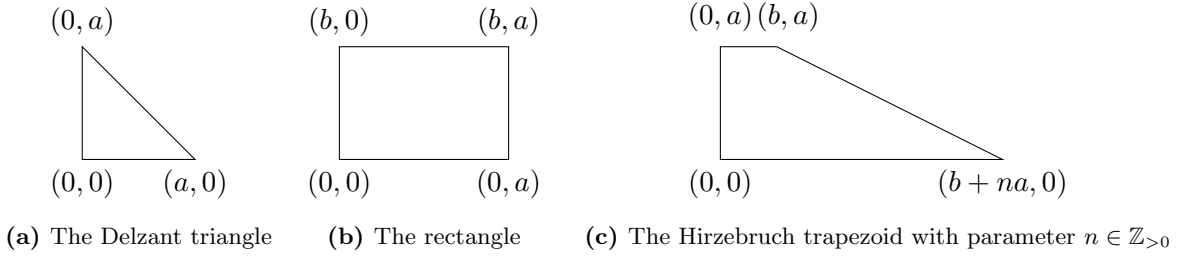
\begin{figure}[t]
 \centering
 \begin{subfigure}{.24\textwidth}
 \centering
 \begin{tikzpicture}[scale = 0.75]
 \draw (0,0) -- (0,2) -- (2,0) -- cycle;
 \node at (0,0) [anchor = north] {$(0,0)$};
 \node at (0,2) [anchor = south] {$(0,a)$};
 \node at (2,0) [anchor = north] {$(a,0)$};
 \end{tikzpicture}
 \caption{The Delzant triangle}
 \label{fig:del-ex-triangle}
 \end{subfigure}\,
 \begin{subfigure}{.21\textwidth}
 \centering
 \begin{tikzpicture}[scale = 0.75]
 \draw (0,0) -- (0,2) -- (3,2) -- (3,0) -- cycle;
 \node at (0,0) [anchor = north] {$(0,0)$};
 \node at (0,2) [anchor = south] {$(b,0)$};
 \node at (3,2) [anchor = south] {$(b,a)$};
 \node at (3,0) [anchor = north] {$(0,a)$};
 \end{tikzpicture}
 \caption{The rectangle}
 \label{fig:del-ex-rect}
 \end{subfigure}\,
 \begin{subfigure}{.52\textwidth}
 \centering
 \begin{tikzpicture}[scale = 0.75]
 \draw (0,0) -- (0,2) -- (1,2) -- (5,0) --cycle;
 \node at (0,0) [anchor = north] {$(0,0)$};
 \node at (0,2.1) [anchor = south] {$(0,a)$};
 \node at (1.2,2.1) [anchor = south] {$(b,a)$};
 \node at (5,0) [anchor = north] {$(b+na,0)$};
 \end{tikzpicture}
 \caption{The Hirzebruch trapezoid with parameter $n\in\Z_{>0}$}
 \label{fig:del-ex-hirz}
 \end{subfigure}
 \caption{Three examples of Delzant polygons, with parameters $a,b>0$.}
 \label{fig:del-ex}
\end{figure}

\begin{Theorem}\label{thm:toric-examples}
 Let $\De$ be a Delzant polygon and let $\rho_{{\rm T}}(\De)$ denote its toric packing density, as in Definition~{\rm \ref{def:toric-density}}. Then the toric packing densities of the examples from Definition~{\rm \ref{def:delz-examples}} $($and shown in Figure~{\rm \ref{fig:del-ex}}$)$ are as follows:
 \begin{itemize}\itemsep=0pt
 \item If $\De$ is the Delzant triangle with parameter $a>0$, then $\rho_{{\rm T}}(\De) = 1$.
 \item If $\De$ is the rectangle with parameters $a,b>0$, then \[\rho_{{\rm T}}(\De) = \frac{\min \{a,b\}}{\max\{a,b\}}.\]
 \item If $\De$ is the Hirzebruch trapezoid with parameters $a,b>0$ and $n\in\Z_{>0}$, then
 \begin{equation*}
 \rho_{{\rm T}}(\De) = \begin{cases}
 \dfrac{2a}{na+2b} & \text{if }n>1,\vspace{1mm}\\
 \dfrac{2a}{a+2b} & \text{if }n=1 \text{ and }a\leq b,\vspace{1mm}\\
 \dfrac{a^2+b^2+ (a-b)^2}{(a+2b)a} & \text{if }n=1 \text{ and } b < a < 2b,\vspace{1mm}\\
 \dfrac{a^2+2b^2}{(a+2b)a} & \text{if }n=1 \text{ and } a \geq 2b. \end{cases}
\end{equation*}
 \end{itemize}
\end{Theorem}

The cases of the Delzant triangle with any parameter $a>0$ and the rectangle with parameters~$a=b$ were already known to Pelayo~\cite{Pe2006}, who listed all Delzant polytopes (in all dimensions) which admit a packing of density 1.

\subsubsection{Semitoric packing}

We can carefully extend the techniques from the case of packing Delzant polygons to the more complicated case of packing semitoric polygons.
Let $\stpoly$ denote a semitoric polygon, this notation is explained in Section~\ref{sec:semitoric}.
Similar to the Delzant case, let $\PST(\stpoly)$ denote the set of admissible semitoric packings of $\stpoly$, and define $\Lambda \colon \PST(\stpoly)\to\R^d$ by taking the $\SL$-size of the triangles included in the packing. That is,
\begin{align}\label{eqn:LambdaST} \Lambda \colon\ \PST(\stpoly)   \to \R^d, \qquad \bigcup_{i=1}^d [B_i] & \mapsto (\lambda_1,\ldots, \lambda_d),
\end{align}
where $d$ is the number of Delzant plus hidden vertices of the semitoric polygon in question, and~$\lambda_i$ is the $\SL$-size of the semitoricly packed triangle $[B_i]$, as described in Definition~\ref{def:STadm-packed-triangle}. We denote a semitoric packing by $[P]$, as described in Definition~\ref{def:st-packing}.
Again, the map $\Lambda$ depends on a choice of labeling $p_1,\ldots,p_d$ of the vertices, which we assume is clockwise starting with the lexicographically minimal (i.e.,~bottom left) vertex, but the magnitude of the resulting vector is independent of this choice.

In Section~\ref{sec:semitoric-edge},
 we describe how we can also obtain real numbers $\ell_1,\ldots, \ell_d>0$ from a~semitoric polygon $\stpoly$, which is not quite the same as the $\SL$-lengths of the edges of $\De$.
Some of the vertices in semitoric polygons are so-called ``fake'' vertices, and $\ell_i$ is the length of the $i^\mathrm{th}$ \emph{semitoric edge}, which is a chain of edges connected by fake corners, as described in Definition~\ref{def:st-edge}.
We also obtain an additional constraint when compared to the Delzant case, which is encoded in one extra parameter $\alpha_i\in [0,\infty]$, given explicitly in equation~\eqref{eqn:alphai} in the statement of Lemma~\ref{lem:avoid-markedpt}. This extra constraint comes from the requirement that the packed triangles in a semitoric polygon avoid the marked points, and we take $\alpha_i=\infty$ in the case in which it does not add an additional constraint. Now, define
\begin{equation}\label{eqn:PSTtilde}
 \PSTtilde(\stpoly) = \left\{(\lambda_1,\ldots,\lambda_d)\in\R^d \,\bigg|\, \begin{matrix}\lambda_i\geq 0,\  \lambda_i+\lambda_{i+1}\leq \ell_i \text{, and}\\  \lambda_i\leq \alpha_i\text{ for all }i\in\{1,\ldots,d\}\end{matrix}\right\},
\end{equation}
where $\alpha_i$ is as in equation~\eqref{eqn:alphai}.

\begin{Theorem}\label{thm:semitoric}
 Let $\stpoly$ be any semitoric polygon. Then
 \[
 \Lambda(\PST(\stpoly)) = \PSTtilde(\stpoly),
 \]
 and $\Lambda$ is a bijection from $\PST(\stpoly)$ to the compact convex polytope $\PSTtilde(\stpoly)$.
 Furthermore,
 \[
 \area([P]) = \frac{1}{2}||\Lambda([P])||^2
 \]
 for any packing $[P]\in\PST(\stpoly)$.
\end{Theorem}

Theorem~\ref{thm:semitoric} is proven in Section~\ref{sec:semitoric-compute}. Theorem~\ref{thm:semitoric} is the direct analogue of Proposition~\ref{prop:toric} in the semitoric case, allowing the problem of finding the maximal packing density of a semitoric polygon to be reduced to finding a vector of maximal magnitude in a polygonal set, which must occur at one of the vertices. Thus, one can also obtain an algorithm in this scenario. The only additional difficulty is finding the constraints $\alpha_i$.

In Section~\ref{sec:pack-st-examples}, we use Theorem~\ref{thm:semitoric} to compute the packing density of several important semitoric polygons,
labeled as types (1)--(3) and depending on parameters $a,b, h, h_1, h_2\in\R$ and~$n\in\Z$ as shown in Figure~\ref{fig:semitoric-polygon-examples}.
Combining equations~\eqref{eqn:stpack-1}, \eqref{eqn:stpack-2}, \eqref{eqn:stpack-3a}, \eqref{eqn:stpack-3b}, and~\eqref{eqn:stpack-3c}, derived in that section, we have the following.

\begin{Theorem}\label{thm:semitoric-examples}
 Let $\stpoly$ be a semitoric polygon and let $\rho_{{\rm ST}}:=\rho_{{\rm ST}}(\stpoly)$ denote its semitoric packing density, as in Definition~{\rm \ref{def:semitoric-density}}.
 Then the semitoric packing densities of the examples from Definition~{\rm \ref{def:st-examples}} $($and shown in Figure~{\rm \ref{fig:semitoric-polygon-examples}}$)$ are as follows:
 \begin{itemize}\itemsep=0pt
 \item If $\stpoly$ is minimal of type $(1)$ with parameters $a>0$ and $h\in(0, a/2)$, then
 \[
 \rho_{{\rm ST}} = \frac{a^2+2h^2-2ah}{a^2}.
 \]
 \item If $\stpoly$ is minimal of type $(2)$ with parameters $a>0$, $b\geq 0$ and $h_1,h_2\in(0,a)$, then
 \[
 \rho_{{\rm ST}} = \frac{a}{a+b}.
 \]
 \item If $\stpoly$ is minimal of type $(3a)$ with parameters $a>0$, $b>0$, $h\in(0, a)$, and $n\in\Z_{\geq 1}$, then
 \begin{equation*}
 \rho_{{\rm ST}} = \begin{cases}
 \dfrac{2a}{na+2b} & \text{if }n>1,\vspace{1mm}\\
 \dfrac{2a}{a+2b} & \text{if }n=1 \text{ and }a\leq b,\vspace{1mm}\\
 \dfrac{a^2+b^2+ (a-b)^2}{(a+2b)a} & \text{if }n=1 \text{ and } b < a < 2b,\vspace{1mm}\\
 \dfrac{a^2+2b^2}{(a+2b)a} & \text{if }n=1 \text{ and } a \geq 2b. \end{cases}
 \end{equation*}
 \item If $\stpoly$ is minimal of type $(3b)$ with parameters $a>0$, $h \in (0,a)$, and $n\in\Z_{\geq 2}$, then
 \[
 \rho_{{\rm ST}} = \frac{2}{n}.
 \]
 \item If $\stpoly$ is minimal of type $(3c)$ with parameters $a>0$, $b\in (-a,0)$, $n\in\Z_{\geq 2}$, and $h\in (0,a+\frac{b}{n-1})$, then
 \[
 \rho_{{\rm ST}} = \frac{(a+b)^2+\min\big\{ ((n-1)a+b - (n-2)h )^2,a^2\big\}}{(na+2b)a}.
 \]
 \end{itemize}
\end{Theorem}

Note that the density of type (2) is independent of the parameters $h_1,h_2$ and the density of types (3a) and (3b) is independent of the parameter $h$.
The polygons from Definition~\ref{def:st-examples} and Figure~\ref{fig:semitoric-polygon-examples} are important because all semitoric polygons can be generated from these by applying a sequence of operations called corner chops and wall chops, see Remark~\ref{rmk:minimal} and~\cite{LFP-fam2}.

\subsubsection{Perfect packings of semitoric polygons}

In~\cite{Pe2006}, Pelayo lists all symplectic toric manifolds (of all dimensions) which admit a toric packing of density 1, which is called a \emph{perfect toric packing}.
Pelayo's proof comes from exploring admissible toric packings of Delzant polytopes. We do the same in the semitoric case, generalizing his result in dimension 4. We say that a semitoric polygon admits a perfect packing if its packing density is 1. In Section~\ref{sec:perfect}, we prove:

\begin{Theorem}\label{thm:perfect}
A semitoric polygon with at least one marked point admits a perfect packing if and only if it is:
\begin{itemize}\itemsep=0pt
 \item minimal of type $(2)$ with parameter $b=0$ $($and any $a>0)$,
 \item minimal of type $(3b)$ with parameter $n=2$ $($and any $a>0)$.
 \item or minimal of type inverted $(3b)$ with $n=2$ $($and any $a>0)$.
\end{itemize}
\end{Theorem}

The terminology in Theorem~\ref{thm:perfect} refers to the minimal semitoric polygons explicitly described in Definitions~\ref{def:st-examples} and \ref{def:inverted-3b}. See Figures~\ref{fig:semitoric-polygon-examples} and \ref{fig:inverted3b}.

\subsubsection{Equivariant capacities}
As described above, the packing densities that we compute for Delzant and semitoric polygons are closely related to certain equivariant symplectic capacities of the associated four dimensional integrable systems.
In Section~\ref{sec:capacities}, we explain how to translate our results about the packing densities of polygons into results about equivariant capacities of integrable systems
 by simply combining Theorems~\ref{thm:toric-examples} and~\ref{thm:semitoric-examples} with the results of Pelayo~\cite{Pe2006, Pe2007} and Figalli--Palmer--Pelayo~\cite{FiPaPe2016} which relate the packing density of polygons with the packing capacities of the associated manifolds.

\subsection{Structure of the paper}
In Section~\ref{sec:symplectic-background}, we describe the background and motivation for problems that we attack in this paper.
In particular, in Section~\ref{sec:capacities}, we explain how to translate our results about packing densities of polygons into results about packing densities and packing capacities of certain four-dimensional integrable systems.
In Section~\ref{sec:polygons-intro}, we give a concise introduction to Delzant and semitoric polygons independent of symplectic geometry.
In Section~\ref{sec:toric}, we describe toric packing and prove Proposition~\ref{prop:toric} and Theorem~\ref{thm:toric-examples}.
In Section~\ref{sec:semitoric}, we describe semitoric packing and prove Theorems~\ref{thm:semitoric} and~\ref{thm:semitoric-examples}.
In Section~\ref{sec:perfect}, we discuss perfect packings of semitoric polygons and prove Theorem~\ref{thm:perfect}.

\section{Background and motivation: equivariant ball packing}
\label{sec:symplectic-background}

In this section we review the background and describe the equivariant ball packing problem for toric and semitoric integrable systems.
We will make use of the results of several authors (primarily those of Atiyah~\cite{Atiyah}, Guillemin--Sternberg~\cite{GuSt1982}, Delzant~\cite{De1988}, Pelayo--V\~{u} Ng\d{o}c~\cite{PeVN2009,PeVN2011}, Pelayo~\cite{Pe2006,Pe2007}, and Figalli--Palmer--Pelayo~\cite{FiPaPe2016}) to translate these packing problems completely into questions purely concerned with polygons, and thus the present section is the only section in this paper which assumes any background in symplectic geometry.
Furthermore, this section is strictly motivational and thus can be skipped if desired.

Moreover, in Section~\ref{sec:capacities} we describe how to translate Theorems~\ref{thm:toric-examples} and~\ref{thm:semitoric-examples} into the language of equivariant symplectic capacities.

\subsection{Integrable systems}

Recall that any function $f$ on a symplectic manifold $(M,\om)$ induces a vector field $\mathcal{X}^f$ on $M$ by the equation $-\mathrm{d}f = \om\big(\mathcal{X}^f, \cdot\big)$.
An \emph{integrable system} is a~triple $(M,\om,F)$ where $(M,\om)$ is a~$2n$-dimensional symplectic manifold and $F=(f_1,\ldots, f_n)\colon M\to \R^n$ satisfies:
\begin{enumerate}\itemsep=0pt
 \item[(1)] $\om\big(\mathcal{X}^{f_i},\mathcal{X}^{f_j}\big)=0$ for all $i,j\in\{1,\ldots,n\}$ (i.e.,~the components of $F$ \emph{Poisson commute}), and
 \item[(2)] $\mathrm{d}f_1, \ldots, \mathrm{d}f_n$ are linearly independent almost everywhere.
\end{enumerate}
The map $F$ is called the \emph{momentum map} of the system.

\subsection{Toric integrable systems}
\label{sec:sympl-toric}
An integrable system $(M,\om,F)$ is called \emph{toric} if $(M,\om)$ is a compact\footnote{Authors sometimes omit the compactness requirement in the definition of toric integrable systems.
Non-compact toric integrable systems were classified by Karshon--Lerman~\cite{karshon-lerman}.} and connected $2n$-dimensional symplectic manifold, each component of the momentum map generates a periodic flow, and the composition of these flows produces an effective $n$-torus action.
That is, the flow of each of $\mathcal{X}^{f_1},\ldots,\mathcal{X}^{f_n}$ is periodic of the same period
and, since the components of $F$ Poisson commute, the flows of $\mathcal{X}^{f_1},\ldots,\mathcal{X}^{f_n}$ commute and thus generate an $n$-torus action, which we require to be effective.

In the 1980s, it was shown by Atiyah~\cite{Atiyah} and Guillemin--Sternberg~\cite{GuSt1982} that given a toric integrable system $(M,\om,F)$, the image $F(M)$ is a convex polytope, and soon after it was shown by Delzant~\cite{De1988} that $\Delta:= F(M)\subset \R^n$ furthermore satisfies three conditions:
\begin{itemize}\itemsep=0pt
 \item $\De$ is \emph{rational} (each edge is directed along an integer vector);
 \item $\De$ is \emph{simple} (exactly $n$ edges meet at each vertex);
 \item $\De$ is \emph{smooth} (for each vertex, the span of the primitive integer vectors directing the edges emanating from that vertex is equal to the lattice $\Z^n$).
\end{itemize}
The polytopes satisfying these conditions are now called \emph{Delzant polytopes}.
Delzant also showed that given any Delzant polytope $\Delta$ there exists exactly one toric integrable system (up to symplectomorphisms which intertwine the momentum maps) with $\Delta$ as its moment image (up to translations and the action of $\SLn$).
That is, up to the appropriate sense of equivalence on each side, \emph{there is a natural bijection from toric integrable systems to the set of Delzant polytopes given by $(M,\om,F)\mapsto F(M)$}. Thus, to study toric integrable systems we can instead study the associated Delzant polytopes.

We will now concentrate on toric integrable systems in dimension four (i.e.,~the case that $n=2$), and denote the components of the momentum map by $F=(J,H)$. In this case the system is classified by its two-dimensional image $F(M)$, which is thus called a \emph{Delzant polygon}. Delzant polygons will be one of the primary objects of interest for this paper, and there is a~detailed description of them in Section~\ref{sec:delzant-poly}.

\subsection{Semitoric integrable systems}

Toric integrable systems are very well understood, but also they are very rigid and thus there are not many systems of interest in physics which are toric.
About ten years ago, Delzant's classification of toric integrable systems was extended in dimension four by Pelayo--V\~{u} Ng\d{o}c~\mbox{\cite{PeVN2009, PeVN2011}} to a much broader and richer class of systems known as \emph{semitoric}.
If $(M,\om,(J,H))$ is a toric integrable system, then $J$ and $H$ both generate a periodic Hamiltonian flow, but in semitoric systems only $J$ is required to generate a periodic Hamiltonian flow.

\begin{Definition} \label{def:semitoricsystem}
A \emph{semitoric integrable system} is a four-dimensional integrable system \[(M,\om,F=(J,H))\] such that
\begin{itemize}\itemsep=0pt
 \item $J$ is proper\footnote{A function is \emph{proper} if the preimage of each compact set is compact; this is automatic if $M$ is compact.} and the Hamiltonian flow of $J$ is periodic of minimal period $2\pi$;
 \item the singular points of $F=(J,H)$ are non-degenerate and have no components of hyperbolic type.
\end{itemize}
\end{Definition}

The second condition refers to the notion of non-degenerate singularities of integrable systems, due to~\cite{Bolsinov-Fomenko}.
We do not explain this condition here, except to say that semitoric systems admit a type of singularities called \emph{focus-focus singularities}, which cannot exist in toric integrable systems. This condition on semitoric systems is discussed in detail in several papers, for instance~\cite[Sections 2.1--2.4]{LFPal}.
Semitoric systems are much more common than toric systems, for instance the coupled angular momentum system~\cite{ADH-momenta,LFP,SaZh1999} and the coupled spin-oscillator~\cite{ADH-spin-osc,PVN-spin-osc} are both semitoric systems, as are various generalizations of the coupled angular momenta~\cite{AHP-twist,HoPa2017}.
The spherical pendulum satisfies all conditions to be semitoric except that $J$ is not proper; such systems are called \emph{generalized semitoric systems} and are studied in~\cite{PeRaVN2015}.

Semitoric systems were classified in terms of a set of invariants by Pelayo and V\~{u} Ng\d{o}c \cite{PeVN2009,PeVN2011}.\footnote{The original classification of Pelayo--V\~{u} Ng\d{o}c imposes the additional generic assumption of \emph{simplicity}, but the classification was recently generalized to include all semitoric systems, simple or not, by Palmer--Pelayo--Tang~\cite{PPT-nonsimple}.}
One of these invariants directly generalizes the Delzant polygon, but there are choices when creating a polygon from a semitoric system. Thus, the \emph{semitoric polygon} is instead an equivalence class of polygons related by a certain group action.
Semitoric polygons are one of the main objects of interest for this paper, and a detailed description of them is given in Section~\ref{sec:semitoric-poly}.

\subsection{Toric and semitoric ball packing}

In the 1980s Gromov proved the non-squeezing theorem~\cite{Gromov1985}, which led to the introduction of an important class of invariants of symplectic manifolds called \emph{symplectic capacities} by Ekeland and Hofer~\cite{EkHo1989, Hocap1990}.
One of the first symplectic capacities discovered was the \emph{Gromov width}, which is the radius of the largest ball which can be symplectically embedded into the given symplectic manifold.

Given a symplectic manifold equipped with a group action, such as an integrable system, one can study the equivariant version of these invariants.
Given a Lie group $G$, a \emph{symplectic $G$-capacity} is an invariant of symplectic manifolds equipped with a symplectic action of the group~$G$, satisfying certain properties which are analogues of the properties of a symplectic capacity.
Symplectic $G$-capacities were introduced in~\cite{FiPaPe2016}.
One type of symplectic $G$-capacities introduced in~\cite{FiPaPe2016} are equivariant packing capacities.

Let $\mathbb{B}^{2n}(r) :=\{z \in \C^n\colon |z| < r\}$ be the $2n$-dimensional open ball of radius $r$. We equip it with the standard symplectic form $\om_0 = \frac{\mathrm{i}}{2} \sum_{j=1}^n {\rm d}z_j \wedge {\rm d}\bar{z_j}$ and an $n$-torus action by component-wise multiplication.

\begin{Definition} \label{def:equivariantembed}
Let $(M,\om)$ be a $2n$-dimensional toric manifold. An embedding $f\colon \mathbb{B}^{2n}(r) \to M $ is \emph{equivariant} if there exists an automorphism $\mathcal{A}\colon \mathbb{T}^n \to \mathbb{T}^n$ such that the following diagram commutes:
\[
\begin{tikzcd}
\mathbb{T}^n \times \mathbb{B}^n(r) \arrow[r, "\mathcal{A} \times f"] \arrow[d]
& \mathbb{T}^n \times M \arrow[d ] \\
\mathbb{B}^n(r) \arrow[r, "f" ]
& M,
\end{tikzcd}
\]
where the two vertical arrows are the torus actions.
\end{Definition}

An equivariant packing of a symplectic toric manifold is a set
\[
 \mathcal{B} = \bigcup_{\alpha\in I} B_\alpha,
\]
where $I$ is an index set, $B_\alpha\subset M$ is the image of an equivariantly embedded ball for all $\alpha\in I$, and the $B_\alpha$ are disjoint.
The \emph{toric packing density of $(M,\om,F)$} is defined to be
\begin{equation}\label{eqn:maximal-toric-packing}
 \rho_{{\rm T}}(M,\om,F) := \mathrm{sup} \left\{\frac{\text{vol}(\mathcal{B})}{\text{vol}(M)} \,  \bigg| \,  \mathcal{B} \text{ is an equivariant packing of $M$}\right\},
\end{equation}
where $\operatorname{vol}$ refers to the typical symplectic volume.

The works of Pelayo~\cite{Pe2006,Pe2007} and Pelayo--Schmidt~\cite{PeSc2008} use the fact that if $B\subset M$ is the image of an equivariantly embedded ball, then $F(B) \subset \De$ is a certain type of simplex, and conversely that any such simplex is the image of an equivariantly embedded ball (for a proof see for instance~\cite[Lemma 2.13]{PeSc2008}).
Thus, an equivariant packing of the toric manifold $(M,\om,F)$ is equivalent to a certain type of packing of the associated Delzant polytope $\De = F(M)$.
Specifically in dimension four, the equivariant ball packing problem is equivalent to packing a Delzant polygon by triangles obtained from affine transformations of the standard isosceles right triangle, i.e., the convex hull of $(0,0)$, $(\lambda,0)$, and $(0,\lambda)$ for some $\lambda>0$.
Furthermore, as a special case of the Duistermaat--Heckman theorem~\cite{DH-measure}, the symplectic volume of $B\subset M$ is proportional to the area of $F(B)\subset F(M)$, so $\rho_{{\rm T}}(M,\om,F)$ is equal to the fraction of $F(M)$ which can be filled by admissibly packed simplices (which are triangles if $n=2$).

Similar to the toric case, we can study equivariant ball packing problems in semitoric systems. V\~{u} Ng\d{o}c~\cite{VN2007} defined moment polygons of semitoric systems and proved that they are an invariant of semitoric systems, and this invariant was combined with several others to form a complete invariant of semitoric systems in~\cite{PeVN2009, PeVN2011} (extended to include so-called ``non-simple'' semitoric systems in~\cite{PPT-nonsimple}).
In semitoric systems there is no global $2$-torus action, but there are certain preferred local actions of the $2$-torus. In~\cite[Section 5]{FiPaPe2016}, the authors describe these local actions and given a semitoric system $(M,\om,F)$ define a \emph{semitoric embedding} $\mathbb{B}^4(r)\hookrightarrow M$ to be an embedding which respects one of these local actions.

We can now ask packing questions about semitoric manifolds.
Similar to the toric case, a~\emph{semitoric packing of $M$} is a union $\mathcal{B} = \bigcup_{\alpha\in I} B_\alpha$ where the $B_\alpha$ are disjoint images of semitoricly embedded balls.
The \emph{semitoric packing density} of a semitoric system $(M,\om,F)$ is defined by
\begin{equation}\label{eqn:maximal-semitoric-packing}
 \rho_{{\rm ST}}(M,\om,F) := \mathrm{sup} \left\{\frac{\text{vol}(\mathcal{B})}{\text{vol}(M)} \,  \bigg| \,  \mathcal{B} \text{ is a semitoric packing of $M$}\right\}.
\end{equation}
Following similar techniques to~\cite{Pe2006, Pe2007}, Figalli--Palmer--Pelayo~\cite{FiPaPe2016} showed that packings of semitoric systems are equivalent to certain types of packings of semitoric polygons.
We describe the rules for packing semitoric polygons in Section~\ref{sec:semitoric}.

We finish this section with an example of how an equivariantly embedded ball in a symplectic toric 4-manifold corresponds to a packed triangle in the associated Delzant polygon.

\begin{figure}
\centering
\begin{subfigure}{.95\textwidth}
\centering
\begin{tikzpicture}
 \draw (0,0) circle (1cm);
 \begin{scope}[even odd rule]
 \clip (-1,0) arc (180:360:1)
 -- (1,0) arc (0:180:1 and 0.4)
 -- cycle;
 \end{scope}
 \draw[dashed] (1,0) arc (0:180:1 and 0.4);
 \draw (-1,0) arc (180:360:1 and 0.4);
 \draw (1.5,-0.4) node[anchor = south] {$\times$};
 \draw (3,0) circle (1cm);
 \begin{scope}[even odd rule]
 \clip (2,0) arc (180:360:1)
 -- (4,0) arc (0:180:1 and 0.4)
 -- cycle;
 \end{scope}
 \draw (2,0) arc (180:360:1 and 0.4);
 \draw[dashed] (4,0) arc (0:180:1 and 0.4);
 \draw[->] (4.5,0) -- (7,0);
 \draw (5.7,0) node[anchor = south] {$F = ( h_1, h_2)$};
 \draw (8,1) node[anchor = south] {$(0,1)$} -- (8,-1) node[anchor = north] {$(0,0)$}-- (10,-1) node[anchor = north] {$(1,0)$} -- (10,1) node[anchor = south] {$(1,1)$} -- (8,1) ;
\end{tikzpicture}
\caption{The image of the moment map is a Delzant polygon, which is a square in this case.}
\label{fig:S2xS2-example-a}
\end{subfigure}\\[.7em]

\begin{subfigure}{.95\textwidth}
 \centering
 \begin{tikzpicture}[scale = .9]
 \draw (0,0)--(0,2)--(2,2)--(2,0)--cycle;
 \draw [line1,line width = 1,fill = fill1] (0,1.5) -- (0,0) -- (1.5,0);
 \draw [dashed, line1,line width = 1] (0,1.5) -- (1.5,0);
 \draw [decorate, line width=1pt, decoration={brace}] (-0.2,0)--(-0.2,1.5);
 \node at (-0.6, 0.75) {$\lambda$};
 \end{tikzpicture}
 \caption{The image under $F$ of an equivariantly embedded ball of radius $R=\sqrt{2\lambda}$.}
 \label{fig:S2xS2-example-b}
\end{subfigure}
\caption{The image of an equivariantly embedded ball in $S^2\times S^2$ is a triangle in the moment image, and the center of the ball is sent to a corner of the moment polygon. Notice that only two out of the three edges of the triangle are included.}
\label{fig:S2xS2-example}
\end{figure}
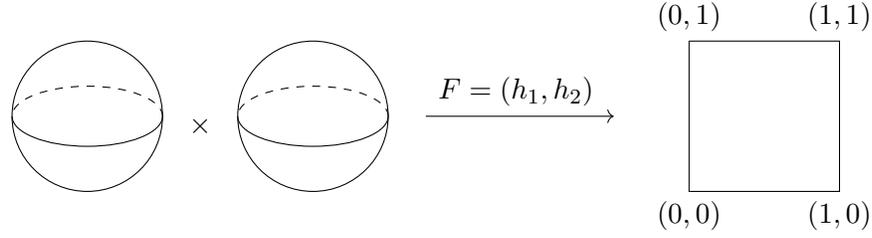
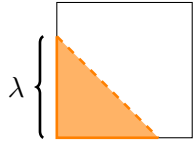

\begin{example}
Endow the usual $2$-sphere $S^2$ with $(h,\theta)$ coordinates where $h\in [0,1]$ and $\theta\in S^1 := \R/ (2\pi \Z)$,
and consider the product $S^2\times S^2$ with coordinates $(h_1, \theta_1, h_2, \theta_2)$.
Consider the integrable system $\big(S^2 \times S^2, \mathrm{d} h_1 \wedge \mathrm{d} \theta_1 + \mathrm{d} h_2 \wedge \mathrm{d}\theta_2, F\big)$ where
\[
 F(h_1, \theta_1, h_2, \theta_2) = (h_1,h_2).
\]
This $F$ generates the 2-torus action which rotates each sphere independently, i.e.,
\[(t_1,t_2) \cdot (h_1, \theta_1, h_2, \theta_2) = (h_1,t_1+ \theta_1, h_2, t_2 +\theta_2) \]
for $(t_1,t_2)\in\mathbb{T}^2:= \R^2/(2\pi \mathbb{Z})^2$.

The corresponding Delzant polygon is a square
\[F\big(S^2\times S^2\big) = [0, 1] \times [0,1],\]
as shown in Figure~\ref{fig:S2xS2-example-a}.
For any $\lambda\in(0,1]$, a $\mathbb{T}^2$-equivariant embedding $f\colon \mathbb{B}^4\big(\sqrt{2\lambda}\big)\to S^2\times S^2$
is given by $f\big(r_1e^{i\theta_1},r_2e^{i\theta_2}\big) = \big(\frac{1}{2}r_1^2, \theta_1, \frac{1}{2}r_2^2, \theta_2\big)$.
The center of the ball is sent to the corner $(0,0)$ of the moment image, and the image of the ball is the set
\[F\circ f\big(\mathbb{B}^4\big(\sqrt{2\lambda}\big)\big) = \{(x,y)\colon x \geq 0, y \geq 0, x+y < \lambda\}.\]
The image of $F\circ f$ in the moment polygon is shown in Figure~\ref{fig:S2xS2-example-b}.
\end{example}

\subsection{Densities and capacities of integrable systems}
\label{sec:capacities}

The results we have obtained for packing polygons (see Theorems~\ref{thm:toric-examples} and~\ref{thm:semitoric-examples}) automatically imply results about certain invariants called packing densities and packing capacities of toric and semitoric integrable systems.
This is because of results by Pelayo~\cite{Pe2006} and Figalli--Palmer--Pelayo~\cite{FiPaPe2016} which relate the packing densities and capacities of integrable systems to the densities and capacities of the associated polygons.
In this section we explain how our results can be translated into this context, directing the interested reader to~\cite{FiPaPe2016} for details.

Due to~\cite[Lemma 2.12]{Pe2006} for the toric case and~\cite[Proposition 7.4]{FiPaPe2016} for the semitoric case, we have that
\begin{equation}\label{eqn:toric-density-relation}
 \rho_{{\rm T}}(M,\om,F) = \rho_{\rm T}(\De)
\end{equation}
if $(M,\om,F)$ is a 4-dimensional toric integrable system with Delzant polygon $\De$ and
\begin{equation}\label{eqn:semitoric-density-relation}
 \rho_{{\rm ST}}(M,\om,F) = \rho_{\rm ST}(\stpoly)
\end{equation}
if $(M,\om,F)$ is a semitoric integrable system with semitoric polygon $\stpoly$.
For both of these results, the proofs make use of the concept of \emph{Duistermaat--Heckman measure}~\cite{DH-measure}, which relates the symplectic volume of a set with the measure of its image under the momentum map.
Using equations~\eqref{eqn:toric-density-relation} and~\eqref{eqn:semitoric-density-relation}, the results stated in Theorems~\ref{thm:toric-examples} and~\ref{thm:semitoric-examples} can thus be immediately translated into results about the packing densities of the integrable systems associated to the given toric and semitoric polygons.

In~\cite{FiPaPe2016}, the toric packing capacity $c_{{\rm T}}$ and semitoric packing capacity $c_{{\rm ST}}$ are defined.
Using~\cite[equation (2)]{FiPaPe2016}, we obtain the following formulas relating the capacities to the packing densities given in equations~\eqref{eqn:maximal-toric-packing} and~\eqref{eqn:maximal-semitoric-packing}:
 \begin{align}
& c_{{\rm T}}(M,\om,F) = \big( 2 \text{vol}(M) \rho_{{\rm T}}(M,\om,F)\big)^{\frac{1}{4}},\label{eqn:toric-cap-to-density} \\
& c_{{\rm ST}}(M,\om,F) = \big( 2 \text{vol}(M) \rho_{{\rm ST}}(M,\om,F) \big)^{\frac{1}{4}}.\label{eqn:semitoric-cap-to-density}
 \end{align}
 The exponent and scaling are needed to make the invariant satisfy the capacity axioms of \emph{monotonicity} and \emph{conformality}, discussed in detail in~\cite{FiPaPe2016}.

Thus, applying equations~\eqref{eqn:toric-density-relation},
\eqref{eqn:semitoric-density-relation}, \eqref{eqn:toric-cap-to-density}
and~\eqref{eqn:semitoric-cap-to-density} to Theorems~\ref{thm:toric-examples}
and~\ref{thm:semitoric-examples} allows one to obtain the toric and semitoric packing capacities (and packing densities) of any example of a~toric integrable system or semitoric integrable system with these polygons.

\begin{Remark}[obtaining a system from the polygon]
 Given a Delzant polygon there is an algorithm, due to Delzant~\cite{De1988}, which produces the associated toric integrable system.
 Unlike the toric case, it is non-trivial to find an explicit\footnote{Here by ``explicit'' we mean that both the symplectic manifold and the momentum map are given globally.} semitoric system with a prescribed semitoric polygons, but explicit systems have been found for many of the semitoric polygons discussed in the present paper.
 For the following we use terminology which will be introduced in Definition~\ref{def:st-examples}.
 For example, the coupled angular momenta system~\cite{ADH-momenta,LFP,SaZh1999} is minimal of type (3a) with parameter $n=2$, and for certain choices of parameters the system introduced in~\cite{HoPa2017} is minimal of type (3a) with $n=2$, of type (2) with $a>0$, or of type (2) with $a=0$.
 The systems on the~$n^{\text{th}}$ Hirzebruch surface introduced in~\cite{LFPal} (which are well understood but not as explicit as the other examples) are minimal of type (3a) with parameter $n$. In~\cite{LFPal}, the authors also introduce a~strategy to construct an explicit semitoric system with a given semitoric polygon, which is used in~\cite{HohMeu} to construct a system whose semitoric polygon is an octagon with four marked points.
 \end{Remark}

\subsubsection{Relation to almost-toric fibrations and infinite staircases}
\label{sec:ATFs}

Semitoric integrable systems are a special case of so-called \emph{almost-toric fibrations $($ATFs$)$}, as in Symmington~\cite{Sy2003}, and there is also a type of marked polygon associated to each ATF of a~symplectic 4-manifold.
Roughly, the main difference between the polygons associated to ATFs and those associated to semitoric systems is that for an ATF the dotted lines (the ``cuts'') do not have to be parallel. Any ATF with parallel cuts (and in particular, any ATF with exactly one marked point) can be viewed as a semitoric system.

Recently, ATFs have emerged as an important tool for computing capacities of certain 4-manifolds. The general idea is that the existence of certain triangles in the polygon associated to the ATF corresponds to balls symplectically embedded in the symplectic manifold which are in some sense compatible with the ATF (this is essentially the same as we have done in the present paper for the special case of semitoric systems).
When trying to compute the maximum radius of a ball which can be embedded into $(M,\om)$, \emph{a priori} there is no reason to assume that it is compatible with a given ATF, but the balls obtained from ATFs at least give a lower bound on the maximum radius.
This argument also works for embedding ellipsoids or polydisks instead of balls: ATFs can be used to find examples of symplectic embeddings which gives a~lower bound on the maximum size which can be embedded.
Equivalently, some authors use these techniques instead to obtain upper bounds on how much a given symplectic manifold has to be scaled by to admit a symplectic embedding of a certain fixed ball, ellipsoid, or polydisk.
Typically, the strategy is to use a~sequence of different ATFs on the same symplectic manifold to obtain a~sequence of improving bounds, hopefully approaching the true value.

The idea of using ATFs to obtain such bounds on these types of capacities was introduced in~\cite{CV,CGHMP} and used recently in~\cite{four-periodic,macgill,MMW-staircasepatterns} to obtain results on so-called `infinite staircase' behavior in capacities, which was first discovered and studied by McDuff--Schlenk~\cite{MS12}. Of course, while ATFs can be used to obtain lower (respectively upper) bounds, actually computing the capacities (as is done in those papers) additionally requires some technique for obtaining an upper (respectively lower) bound as well, which typically involves \emph{$J$-holomorphic curves}.
The equivariant capacities computed in the present paper, and those that appeared in~\cite{FiPaPe2016}, similarly imply bounds for traditional (i.e.,~non-equivariant) capacities.

\section{Delzant and semitoric polygons}
\label{sec:polygons-intro}

In this section we focus on the combinatorial aspects of Delzant and semitoric polygons, bypassing the symplectic geometry.

\subsection{Delzant polygons}
\label{sec:delzant-poly}

 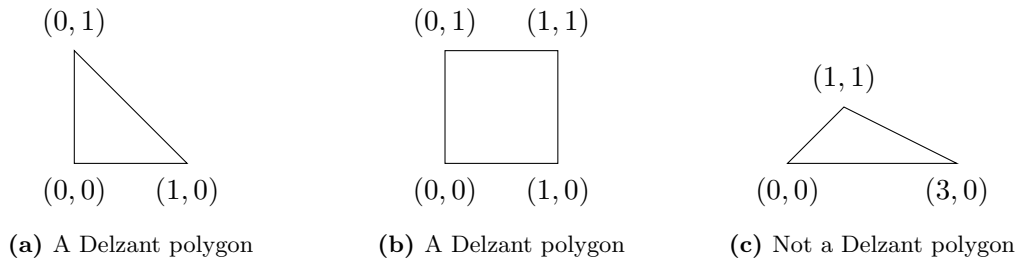
\begin{figure}
 \centering
 \begin{subfigure}{.3\textwidth}
 \centering
 \begin{tikzpicture}[scale = 0.75]
 \draw (0,0) -- (0,2) -- (2,0) -- cycle;
 \node (n1) at (0,-0.5) {$(0,0)$};
 \node (n2) at (2,-0.5) {$(1,0)$};
 \node (n3) at (0,2.5) {$(0,1)$};
 \end{tikzpicture}
 \caption{A Delzant polygon}
 \label{fig:delz1}
 \end{subfigure}
 \begin{subfigure}{.3\textwidth}
 \centering
 \begin{tikzpicture}[scale = 0.75]
 \draw (0,0) -- (0,2) -- (2,2) -- (2,0) -- cycle;
 \node (n1) at (0,-0.5) {$(0,0)$};
 \node (n2) at (2,-0.5) {$(1,0)$};
 \node (n2) at (2, 2.5) {$(1,1)$};
 \node (n3) at (0,2.5) {$(0,1)$};
 \end{tikzpicture}
 \caption{A Delzant polygon}
 \label{fig:delz2}
 \end{subfigure}
 \begin{subfigure}{.3\textwidth}
 \centering
 \begin{tikzpicture}[scale = 0.75]
 \draw (0,0) -- (1,1) -- (3,0) -- cycle;
 \node (n1) at (0,-0.5) {$(0,0)$};
 \node (n2) at (1,1.5) {$(1,1)$};
 \node (n3) at (3, -0.5) {$(3,0)$};
 \end{tikzpicture}
 \caption{Not a Delzant polygon}
 \label{fig:not-delz}
 \end{subfigure}
 \caption{Two Delzant polygons (see Figures~\ref{fig:delz1} and~\ref{fig:delz2}) and a polygon that fails the Delzant condition at the vertex $(1,1)$ Figure~\ref{fig:not-delz}.}
 \label{fig:delz-poly-examples}
 \end{figure}

A vector $u\in\mathbb{Z}^2$ is called \emph{primitive} if $u = kw$ for $w\in\mathbb{Z}^2$ and $k\in\Z$ implies that $k=\pm 1$ (i.e.,~$u$ is the shortest non-zero integral vector in the line that it describes). Given two vectors $u_1,u_2\in\mathbb{R}^2$, we will use $\mathrm{det}(u_1,u_2)$ to denote the determinant of the matrix whose first column is $u_1$ and whose second column is $u_2$.
A polygon is called \emph{rational} if each edge is directed along an integer vector.

\begin{definition}\label{def:delzant_cond}
Let $u_1,u_2\in\mathbb{Z}^2$ be the primitive integer vectors directing the edges emanating from a vertex of a rational polygon. Then we say that the vertex satisfies the \emph{Delzant condition} if the integer span of $u_1$ and $u_2$ is all of $\Z^2$. This is equivalent to the condition that ${\det (u_1,u_2) = \pm 1}$.
\end{definition}

\begin{Definition}\label{def:delzant} A \emph{Delzant polygon} is a rational convex polygon $\De\subset\R^2$ such that each vertex satisfies the Delzant condition.
\end{Definition}

Two specific examples of Delzant polygons, and one non-example, are shown in Figure~\ref{fig:delz-poly-examples}.

\begin{Remark}
The definition of Delzant polygon that we use is equivalent to the more general notion of a Delzant polytope (described in Section~\ref{sec:sympl-toric}) in the case that the polytope is two dimensional, which is the setting of this paper.
Of the three conditions to be a Delzant polytope the second condition, simplicity, is automatic when the dimension is $n=2$, and the third condition, smoothness, is equivalent to requiring all vertices to satisfy the Delzant condition.
\end{Remark}


We now explain three important examples of Delzant polygons. Any Delzant polygon can clearly be translated in $\R^2$ without affecting the Delzant condition on the corners, so for simplicity in each example we assume that one of the corners is at the origin.

\begin{Definition}\label{def:delz-examples}
 Let $a,b>0$. Then we have the following examples of Delzant polygons:
\begin{itemize}\itemsep=0pt
\item the \emph{Delzant triangle of side length $a$} is the right triangle which has vertices $(0,0)$, $(a,0)$, and $(0,a)$. See Figure~\ref{fig:del-ex-triangle}.

\item the \emph{Delzant rectangle with parameters $a,b$} is the typical rectangle, which has vertices at $(0,0)$, $(0,a)$, $(b,a)$, and $(b,0)$. See Figure~\ref{fig:del-ex-rect}.

\item given any $n\in \Z_{\geq 1}$ the \emph{Hirzebruch trapezoid} is the polygon which has vertices $(0,0)$, $(0,a)$, $(b,a)$,
and $(b+na,0)$. See Figure~\ref{fig:del-ex-hirz}.
\end{itemize}
\end{Definition}

Taking the determinant of the primitive vectors directing each pair of adjacent edges, it is straightforward to check that each of these examples is indeed a Delzant polygon.

\subsection{Semitoric polygons}
\label{sec:semitoric-poly}

Delzant polygons (and the associated toric integrable systems) are well understood, but represent a relatively restrictive class of integrable systems.
On the other hand, semitoric polygons, though they are more complicated, represent a much broader class of integrable systems.
Thus, though they can be more difficult to work with, it is worthwhile to extend techniques from the Delzant setting to the more general semitoric setting.

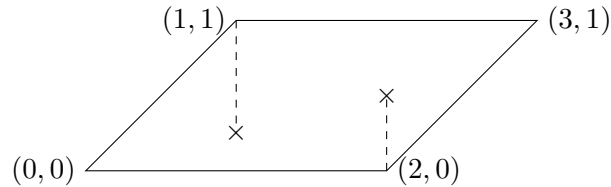
\begin{figure}
 \centering
 \begin{tikzpicture}[scale=1]
 \draw (0,0)--(2,2)--(6,2)--(4,0)--cycle;
 \node (n1) at (2,0.5) {$\times$};
 \draw [dashed](2,2)--(2,0.5);
 \node (n2) at (4, 1) {$\times$};
 \draw [dashed](4,0)--(4,1);
 \node at (0,0) [anchor = east] {$(0,0)$};
 \node at (2,2) [anchor = east] {$(1,1)$};
 \node at (6,2) [anchor = west] {$(3,1)$};
 \node at (4,0) [anchor = west] {$(2,0)$};
 \end{tikzpicture}
 \caption{A representative of a semitoric polygon which has two marked points and two cuts (one up and one down). The marked points are denoted by the $\times$, the cuts are denoted by the dashed lines.}
 \label{fig:representative}
\end{figure}

Semitoric polygons are similar to Delzant polygons, since they are also convex rational polygons with certain conditions on the corners, but they are more complicated for two reasons:
\begin{enumerate}\itemsep=0pt
 \item[(1)] semitoric polygons are decorated with marked points in their interior and ``cuts'' which go up or down from each marked point;
 \item[(2)] A ``semitoric polygon'' is actually an infinite family of polygons related by certain transformations which change the cut direction or skew the entire polygon by a certain linear transformation.
\end{enumerate}
Any single representative of the semitoric polygon determines the entire infinite family, but nevertheless it is natural to consider all representatives together.
Roughly, this is because the cuts in a semitoric polygon distort the boundary, so the most natural way to study local properties of the polygon is to consider an equivalent polygon where the cuts are moved away from the area of interest. In fact, this is exactly what we will do when considering packings of semitoric polygons.

Recall the Delzant corner condition of Definition~\ref{def:delzant}. Here we describe some other conditions which will be relevant for the corners on the cuts of semitoric polygons.
Let
\begin{equation}\label{eqn:T}
T := \begin{pmatrix} 1 & 0 \\ 1 & 1 \end{pmatrix}.
\end{equation}

\begin{Definition}\label{def:fake_hidden}
Let $v$ be a vertex of a rational polygon and $u_1$, $u_2$ be the primitive integral vectors directing the edges emanating from $v$ ordered so that $\det (u_1, u_2) > 0$. Let $k\in\Z_{> 0}$. Then \begin{itemize}\itemsep=0pt
 \item $v$ satisfies the \emph{fake corner condition for $k$ cuts} if $\det \big(u_1,T^k u_2\big) = 0$, and
 \item $v$ satisfies the \emph{hidden corner condition for $k$ cuts} if $\det \big(u_1,T^k u_2\big) = 1$.
\end{itemize}
\end{Definition}

We will start by describing a single representative of a semitoric polygon.
Given a point~$c\in\R^2$ and $\epsilon\in\{1,-1\}$, let $L_c^\epsilon$ denote the ray which starts at $c$ and is directed along the vector $(0,\epsilon)$, so it goes up if $\epsilon=1$ and down if $\epsilon=-1$.

\begin{Definition}\label{def:stpoly-representative} A \emph{semitoric polygon representative} is a triple $\strep$, where
\begin{enumerate}\itemsep=0pt
 \item[(1)] $\Delta \subset \mathbb{R}^{2}$ is a convex rational polygon;
 \item[(2)] $\vec c = (c_1,\dots,c_m)$ is the set of \emph{marked points}, and $c_i \in \textrm{int}(\Delta)$ for $i = 1,\dots,m$;
 \item[(3)] $\vec{\epsilon} = (\epsilon_1,\dots,\epsilon_m) \in \{ \pm1 \}^m$ is the collection of \emph{cut directions}.
\end{enumerate}
For each $i\in\{1,\ldots,m\}$, we call the ray $L_{c_i}^{\epsilon_i}\cap\Delta$ a \emph{cut}.
Furthermore, we require that:
 \begin{itemize}\itemsep=0pt
 \item each point where exactly $k$ cuts intersect the boundary $\partial \Delta$ is a vertex of $\Delta$ which satisfies either the fake or hidden corner condition for $k$ cuts (and are thus known as \emph{fake corners} or \emph{hidden corners}, respectively);
 \item all other vertices of $\Delta$ satisfy the Delzant corner condition (and are known as \emph{Delzant corners}).
 \end{itemize}
\end{Definition}

\begin{figure}
 \begin{center}
 \begin{tikzpicture}[scale = 1.5]
 \draw (0,0) -- (4,0) -- (2,1) -- cycle;
 \draw [line width = 2pt,-stealth,blue] (2,1) -- (0,0);
 \draw [line width = 2pt,-stealth,blue] (2,1) -- (4,0);
 \draw [dashed] (2,1) -- (2,0.5);
 \node (X) at (2,0.5) {$\times$};
 \node (c1) at (-0.5,0) {$(0,0)$};
 \node (c2) at (4.5,0) {$(4,0)$};
 \node (c3) at (2,1.3) {$(2,1)$};
 \node (w) at (0.4,0.9) {$u_1 = \begin{pmatrix} -2\\-1 \end{pmatrix}$};
 \node (v) at (3.3,0.9) {$u_2 = \begin{pmatrix} 2\\-1 \end{pmatrix}$};
 \end{tikzpicture}
 \caption{A semitoric polygon with one marked point, two Delzant corners, and one fake corner (on the cut). The vertex at $(2,1)$ satisfies the fake corner condition since $\mathrm{det}(u_1,Tu_2) = 0$.}
 \label{fig:sem1}
 \end{center}
 \end{figure}
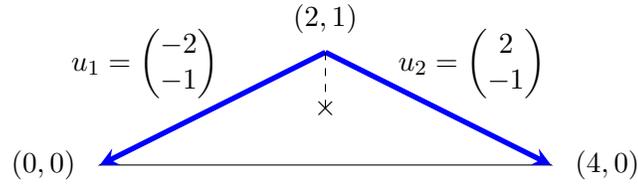

We also assume that the $c_i$ are ordered lexicographically. That is, if $c_i = (x,y)$ and $c_j = (x',y')$ then $i<j$ if and only if either $x<x'$ or both $x=x'$ and $y < y'$.
See Figures~\ref{fig:representative} and~\ref{fig:sem1} for examples of semitoric polygon representatives.

 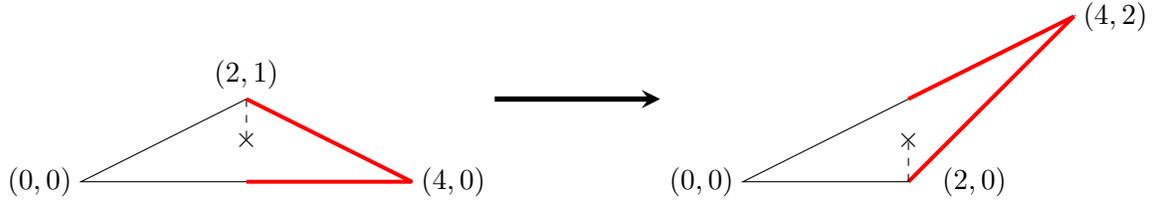
\begin{figure}
 \centering
 \begin{tikzpicture}[scale = 1.1]
 \draw (0,0) -- (4,0) -- (2,1) -- cycle;
 \draw [dashed] (2,1) -- (2,0.5);
 \draw [red,line width = 1.5] (2,1) -- (4,0);
 \draw [red,line width = 1.5] (2,0) -- (4,0);
 \node (X) at (2,0.5) {$\times$};
 \node (c1) at (-0.5,0) {$(0,0)$};
 \node (c2) at (4.5,0) {$(4,0)$};
 \node (c3) at (2,1.3) {$(2,1)$};

 \draw [-stealth, line width = 2] (5,1) -- (7,1);

 \draw (8,0) -- (10,0) -- (12,2) -- cycle;
 \draw [red,line width = 1.5] (10,0) -- (12,2);
 \draw [red,line width = 1.5] (10,1) -- (12,2);
 \draw [dashed] (10,0) -- (10,0.5);
 \node (X) at (10,0.5) {$\times$};
 \node (c1) at (7.5,0) {$(0,0)$};
 \node (c2) at (10.8,0) {$(2,0)$};
 \node (c3) at (12.5,2) {$(4,2)$};

 \end{tikzpicture}
 \caption{Changing an up cut to a down cut by applying $T$ to the right side (in red) of the polygon. Notice how the vertex originally at $(2,1)$ disappears after the transformation. This is precisely what is meant by it being a ``fake corner''.} \label{fig:change-cut}
 \end{figure}

There is a natural notion of equivalence for semitoric polygon representatives, and the actual semitoric polygon will be a set of equivalent representatives.
Given any function $f\colon \R^2\to \R^2$, we define an action of $f$ on $\strep$ by $f\cdot \strep = (f(\Delta),f(\vec{c}), \vec{\epsilon}\,)$,
where $f(\vec{c}) = (f(c_1), \ldots, f(c_m))$.
We consider three types of transformations on a semitoric polygon representative:
\begin{itemize}\itemsep=0pt
 \item \emph{Changing a cut direction:} replacing $\epsilon_i$ by $-\epsilon_i$ and, after shifting in the $x$-direction to place the origin on the same vertical line as the marked point $c_i$, acting on $\strep$ by the transformation which acts as $T^{\epsilon_i}$ on the half-plane to the right of $c_i$
 and as the identity on the half-plane to the left of $c_i$, where $x_i$ is the $x$-coordinate of $c_i$ (see Figure~\ref{fig:change-cut});
 \item \emph{$T$ acting globally:} $T$ acts on $(\Delta, \vec c, \vec\epsilon)$;
 \item \emph{Shifting:} A vertical translation acts on $\strep$.
\end{itemize}

To formally capture this behavior, we will now define a group action on the set of semitoric polygon representatives, and the semitoric polygon will be the orbit under this action.
Let $\mathcal{T}$ be the group generated by powers of $T$ and vertical translations of $\R^2$, and for $j\in\R$ let $\mathfrak{t}_j \colon \R^2\to\R^2$ be given by
\begin{equation}\label{eqn:frakt}
\mathfrak{t}_j(x,y)= \begin{cases}(x, y + x - j) & \textrm{if } x\geq j, \\ (x, y) & \textrm{otherwise.} \end{cases}
\end{equation}
Applying $\mathfrak{t}_j$ is equivalent to applying $T$ in coordinates whose origin is on the line $\{x=j\}$.
Let $G_m = \{1, -1\}^m$.
Then $\mathcal{T}\times G_m$ acts on a semitoric polygon representative by
\begin{equation}\label{eqn:group-action}
 \big(\tau,\vec{\epsilon}\,{}'\big)\cdot \strep = \big(\sigma(\De), \sigma(\vec{c}), (\epsilon_1' \epsilon_1,\ldots, \epsilon_m' \epsilon_m)\big)
\end{equation}
where \smash{$\sigma = \tau \circ \mathfrak{t}_{\pi_1(c_1)}^{u_1}\circ \cdots \circ \mathfrak{t}_{\pi_1(c_m)}^{u_m}$} and
$u_k = \epsilon_k (1-\epsilon_k')/2$.
Thus, acting by $(\tau,\vec{\epsilon}\,{}')\in \mathcal{T}\times G_m$ applies~$\tau$ globally on the polygon, shifting and/or skewing the entire polygon (along with the marked points), and changes the direction of the $i^{\mathrm{th}}$ cut between being up and down for all $i$ such that $\epsilon'_i = -1$.
When flipping the $i^{\mathrm{th}}$ cut, it also skews the right-hand side of the polygon, leaving the left-hand side alone (see Figure~\ref{fig:change-cut}).
Also note that the effect of the matrix $T$, and thus also of the changing of a cut direction, is to increase all slopes by $1$.

It is not obvious that acting on a semitoric polygon representative by any element of $\mathcal{T}\times G_m$ will yield another polygon satisfying the requirements to be a semitoric polygon representative, but this was proven in~\cite[Lemma 4.2]{PeVN2011}.

\begin{Definition}\label{def:semitoric-polygon}
Let $\strep$ be a semitoric polygon representative, as in Definition~\ref{def:stpoly-representative}. Then the orbit of $\strep$ under the action of $\mathcal{T}\times G_m$ described in equation~\eqref{eqn:group-action} is denoted by~$\stpoly$ and is called a \emph{semitoric polygon}.
\end{Definition}

A semitoric polygon is an infinite family of polygons, but notice that the entire family is determined by \emph{any single representative}. Thus, it is typical to just provide a single polygon to represent the entire family, which is what we will do in the present paper.

\begin{figure}
 \centering
 \begin{tikzpicture}[scale=0.8]
 \draw (0,0)--(2,2)--(4,2)--(6,0)--cycle;
 \node (n1) at (2,0.5) {$\times$};
 \draw [dashed](2,2)--(2,0.5);
 \node (n2) at (4, 1) {$\times$};
 \draw [dashed](4,2)--(4,1);
 \begin{scope}[shift={(7,0)}]
 \draw (0,2)--(4,2)--(6,0)--(2,0)--cycle;
 \node (n1) at (2,0.5) {$\times$};
 \draw [dashed](2,0)--(2,0.5);
 \node (n2) at (4, 1) {$\times$};
 \draw [dashed](4,2)--(4,1);
 \end{scope}
 \begin{scope}[shift={(7,-3)}]
 \draw (0,2)--(6,2)--(4,0)--(2,0)--cycle;
 \node (n1) at (2,0.5) {$\times$};
 \draw [dashed](2,0)--(2,0.5);
 \node (n2) at (4, 1) {$\times$};
 \draw [dashed](4,0)--(4,1);
 \end{scope}
 \begin{scope}[shift={(0,-3)}]
 \draw (0,0)--(2,2)--(6,2)--(4,0)--cycle;
 \node (n1) at (2,0.5) {$\times$};
 \draw [dashed](2,2)--(2,0.5);
 \node (n2) at (4, 1) {$\times$};
 \draw [dashed](4,0)--(4,1);
 \end{scope}
 \end{tikzpicture}
 \caption{Four representatives of the same semitoric polygon, which has two marked points and is minimal of type (2), see Definition~\ref{def:st-examples}. In each example, the vertices at the far left and right satisfy the Delzant corner conditions, and the two vertices on the cuts satisfy the fake corner condition. See Figure~\ref{fig:representative} for a figure of just one representative with its vertices labeled.}
 \label{fig:semitoric-4reps-example}
\end{figure}
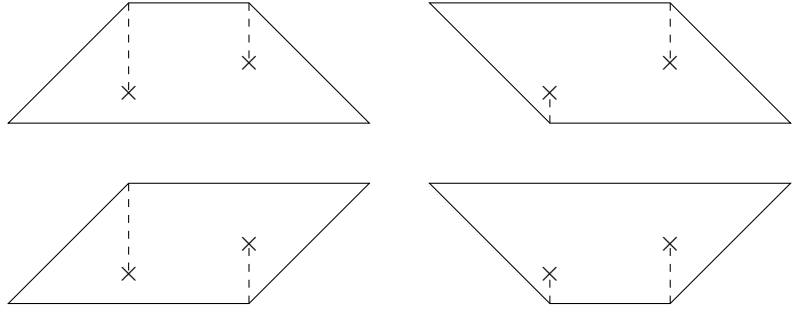

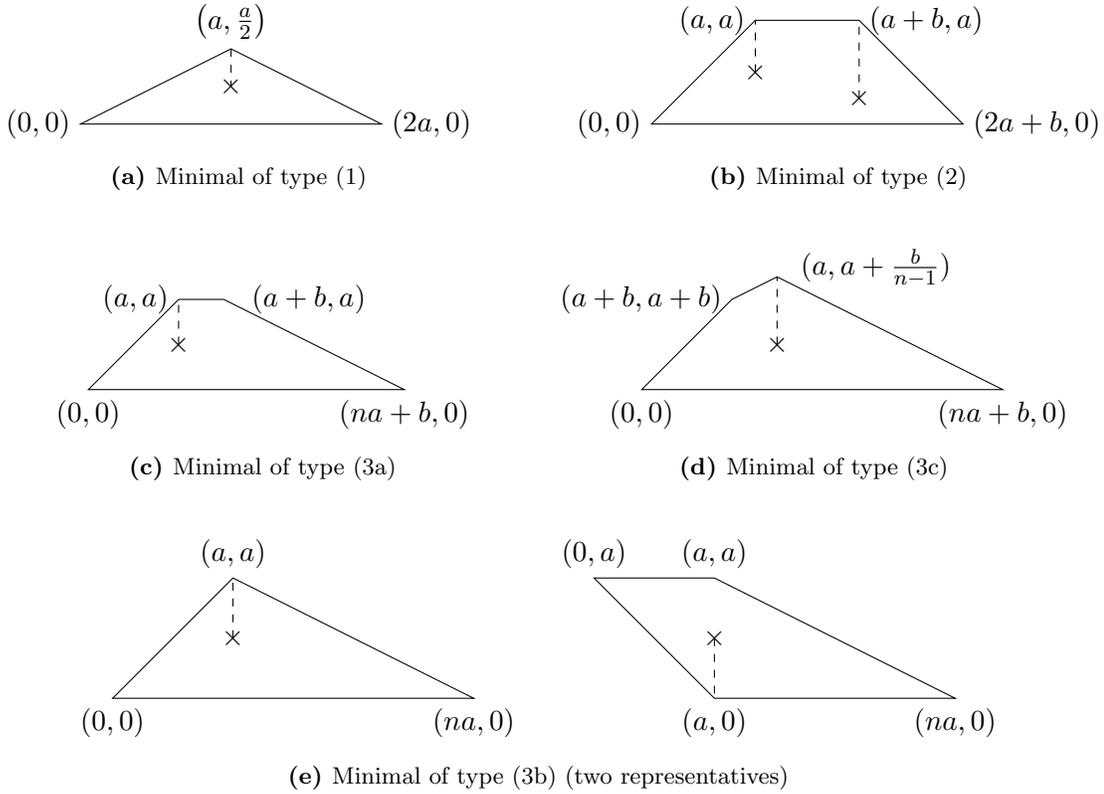
\begin{figure}[t]
 \centering
 \begin{subfigure}{.45\textwidth}
 \centering
 \begin{tikzpicture}[scale = 1]
 \draw (0,0) -- (2,1) -- (4,0) -- cycle;
 \node at (2,0.5) {$\times$};
 \draw [dashed](2,1)--(2,0.5);
 \node at (0,0) [anchor = east] {$(0,0)$};
 \node at (2,1) [anchor = south] {$\left(a,\frac{a}{2}\right)$};
 \node at (4,0) [anchor = west] {$(2a,0)$};
 \end{tikzpicture}
 \caption{Minimal of type (1)}
 \label{fig:stmin-1}
 \end{subfigure}\qquad
 \begin{subfigure}{.45\textwidth}
 \centering
 \begin{tikzpicture}[scale = .69]
 \draw (0,0) -- (2,2) -- (4,2) -- (6,0) -- cycle;
 \node at (2,1) {$\times$};
 \draw [dashed](2,1)--(2,2);
 \node at (4,.5) {$\times$};
 \draw [dashed](4,.5)--(4,2);
 \node at (0,0) [anchor = east] {$(0,0)$};
 \node at (2,2) [anchor = east] {$(a,a)$};
 \node at (4,2) [anchor = west] {$(a+b,a)$};
 \node at (6,0) [anchor = west] {$(2a+b,0)$};
 \end{tikzpicture}
 \caption{Minimal of type (2)}
 \label{fig:stmin-2}
 \end{subfigure}\\[1.5em]
 \begin{subfigure}{.45\textwidth}
 \centering
 \begin{tikzpicture}[scale = .6]
 \draw (0,0) -- (2,2) -- (3,2) -- (7,0) -- cycle;
 \node at (2,1) {$\times$};
 \draw [dashed](2,1)--(2,2);
 \node at (0,0) [anchor = north] {$(0,0)$};
 \node at (2,2) [anchor = east] {$(a,a)$};
 \node at (3.4,2) [anchor = west] {$(a+b,a)$};
 \node at (7,0) [anchor = north] {$(na+b,0)$};
 \end{tikzpicture}
 \caption{Minimal of type (3a)}
 \label{fig:stmin-3a}
 \end{subfigure}
 \begin{subfigure}{.45\textwidth}
 \centering
 \begin{tikzpicture}[scale = .6]
 \draw (0,0) -- (2,2) -- (3,2.5) -- (8,0) -- cycle;
 \node at (3,1) {$\times$};
 \draw [dashed](3,1)--(3,2.5);
 \node at (0,0) [anchor = north] {$(0,0)$};
 \node at (2,2) [anchor = east] {$(a+b,a+b)$};
 \node at (3.3,2.7) [anchor = west] {$(a,a+\frac{b}{n-1})$};
 \node at (8,0) [anchor = north] {$(na+b,0)$};
 \end{tikzpicture}
 \caption{Minimal of type (3c)}
 \label{fig:stmin-3c}
 \end{subfigure}\\[1.5em]
 \begin{subfigure}{.9\textwidth}
 \centering
 \begin{tikzpicture}[scale = .8]
 \draw (0,0) -- (2,2) -- (6,0) -- cycle;
 \node at (2,1) {$\times$};
 \draw [dashed](2,1)--(2,2);
 \node at (0,0) [anchor = north] {$(0,0)$};
 \node at (2,2) [anchor = south] {$(a,a)$};
 \node at (6,0) [anchor = north] {$(na,0)$};
 \begin{scope}[shift = {(8,0)}]
 \draw (0,2) -- (2,2) -- (6,0) -- (2,0) -- cycle;
 \node at (2,1) {$\times$};
 \draw [dashed](2,1)--(2,0);
 \node at (0,2) [anchor = south] {$(0,a)$};
 \node at (2,2) [anchor = south] {$(a,a)$};
 \node at (6,0) [anchor = north] {$(na,0)$};
 \node at (2,0) [anchor = north] {$(a,0)$};
 \end{scope}
 \end{tikzpicture}
 \caption{Minimal of type (3b) (two representatives)}
 \label{fig:stmin-3b}
 \end{subfigure}
 \caption{One representative of each of the minimal semitoric polygons of type (1), (2), (3a), and (3c) with the coordinates of the vertices labeled.
 In Figure~\ref{fig:stmin-3b}, there are two representatives of the minimal semitoric polygon of type (3b) which are related by a change of cut direction operation composed with a~global skewing. The representative shown on the left has a hidden corner at $(a,a)$. Notice that when the cut direction is changed away from the point $(a,a)$, in the representative on the right, a Delzant corner is revealed.}
 \label{fig:semitoric-polygon-examples}
\end{figure}

Four representatives of the same semitoric polygon are shown in Figure~\ref{fig:semitoric-4reps-example}.
All vertices on cuts in Figure~\ref{fig:semitoric-4reps-example} are fake corners, and notice that when the cut direction is changed there is no corner remaining there. Hidden corners have the property that when the cut direction is changed what remains is a Delzant corner, as in Figure~\ref{fig:stmin-3b}.

\begin{Remark}
If there is at most one marked point in each vertical line (i.e.,~if each of the $c_i$ have distinct $c$-coordinates), then the associated semitoric system is called a \emph{simple semitoric system}, and the definition of a semitoric polygon simplifies somewhat. In this case, at most one cut will meet any corner and thus $k=1$ in the definition for the hidden and fake corner conditions.
Simplicity was assumed in the original classification of semitoric systems by Pelayo and V\~{u} Ng\d{o}c~\cite{PeVN2009, PeVN2011}, but this assumption was removed when the classification was extended to include all semitoric systems, simple or not, by Palmer--Pelayo--Tang~\cite{PPT-nonsimple}.
When semitoric polygons were first introduced by V\~{u} Ng\d{o}c~\cite{VN2007} simplicity was not assumed.
\end{Remark}

\begin{Remark}
The semitoric polygons of Definition~\ref{def:semitoric-polygon} are sometimes called \emph{marked semitoric polygons}, as in~\cite{LFPal}, since they include the information of the marked points. The original invariant described and used in~\cite{VN2007} does not include the height of each of the marked points, and in the classification of semitoric systems~\cite{PeVN2009, PeVN2011} the height of each marked point is a separate invariant. Marked semitoric polygons are introduced in~\cite{LFPal}
as a way to organize the invariants together, and in~\cite{HP-extend} the definition is extended to allow for non-simple systems (with multiple marked points in the same vertical line), following~\cite{PPT-nonsimple}.
\end{Remark}

Now we give several important examples of semitoric polygons, all shown in Figure~\ref{fig:semitoric-polygon-examples}.

\begin{Definition}
\label{def:st-examples}
A semitoric polygon $[\De,((a,h)),(+1)]$ is called:
\begin{itemize}\itemsep=0pt
 \item \emph{minimal of type $(1)$} if $\De$ is the polygon with vertices $(0,0)$, $(a,a/2)$, $(2a,0)$, where $a>0$ and $0<h<a/2$ (see Figure~\ref{fig:stmin-1});
 \item \emph{minimal of type $(3a)$} if $\De$ is the polygon with vertices $(0,0)$, $(a,a)$, $(a+b,a)$, $(na+b, 0)$,
 where $a>0$, $b>0$, $n\in\Z$, $n\geq 1$, and $0<h<a$ (see Figure~\ref{fig:stmin-3a});
 \item \emph{minimal of type $(3b)$} if $\De$ is the polygon with vertices $(0,0)$, $(a,a)$, $(na, 0)$,
 where $a>0$, $n\in\Z$, $n\geq 2$, and $0<h<a$ (see Figure~\ref{fig:stmin-3b});
 \item \emph{minimal of type $(3c)$} if $\De$ is the polygon with vertices $(0,0)$, $(a+b,a+b)$, $\big(a,a+\frac{b}{n-1}\big)$, $(na+b,0)$,
 where $a>0$, $-a<b<0$, $n\in\Z$, $n\geq 2$, and $0<h<a+\frac{b}{n-1}$ (see Figure~\ref{fig:stmin-3c}).
\end{itemize}
Furthermore, a semitoric polygon $[\De,((a,h_1), (a+b,h_2)),(+1,+1)]$ is
\begin{itemize}\itemsep=0pt
 \item \emph{minimal of type $(2)$} if $\De$ is the polygon with vertices $(0,0)$, $(a,a)$, $(a+b,a)$, $(2a+b,0)$,
 where $a>0$, $b\geq 0$, and $0<h_i<a$ for $i=1,2$ (see Figure~\ref{fig:stmin-2}).
\end{itemize}
\end{Definition}

In Definition~\ref{def:st-examples} notice that we only have to describe one representative of the semitoric polygon to determine the entire family, so for simplicity we choose a representative with all cuts pointing upwards, which is why $\epsilon_1=\epsilon_2=+1$ in all cases.
Also notice that the parameter $b\in\R$ from type (3c) is strictly negative, to keep the notation in alignment with (3a).

\begin{Remark}\label{rmk:minimal}
 A semitoric polygon is called \emph{minimal} if it cannot be obtained from another semitoric polygon by an operation known as a \emph{corner chop}, see~\cite{KPP2018} for a complete description of this operation (a corner chop on a semitoric polygon corresponds to a \emph{toric type blowup} on the associated semitoric system~\cite[Section 4]{LFPal}).
 The classification of minimal semitoric polygons into types comes from~\cite{KPP2018} (expanding on techniques from~\cite{KPP2015}), in which the authors prove that all minimal semitoric polygons come in one of seven types, numbered (1)--(7).
 In~\cite{LFP-fam2} this theory is further refined to take into account an operation on the polygons called a \emph{wall chop} (corresponding to a \emph{semitoric type blowup} of the system), and any semitoric polygon which cannot be obtained from another one by either a wall chop or a corner chop is called \emph{strictly minimal}. In~\cite{LFP-fam2}, it is shown that each of the polygons of types (4)--(7) can be obtained from another system by a wall chop, and thus all strictly minimal polygons are of type (1), (2), or~(3) (and type (3) splits into three subcases, (3a), (3b), and (3c)).
 Strictly minimal polygons are particularly fundamental (all semitoric polygons can be obtained from these via a sequence of corner chops and wall chops), which is why we choose to focus on them for this paper.
 Note that types (3a), (3b), and (3c) are only actually minimal if $n\neq 3$, since if $n=3$ they can be obtained from type (1) via a corner chop. All other polygons listed in Definition~\ref{def:st-examples} are actually (strictly) minimal.
 Moreover, there is a change of convention between the present paper and the work of Kane--Palmer--Pelayo~\cite{KPP2018}: in the present paper the minimal system of type (3a), (3b), or (3c) with parameter $n$ would instead have parameter $k=n-1$ in~\cite{KPP2018}.
\end{Remark}

\section{Toric packing}\label{sec:toric}

In this section, we describe the packing problem for Delzant polygons, prove Proposition~\ref{prop:toric} about the space of possible packings for a given Delzant polygon, and compute the packing density in some examples, obtaining Theorem~\ref{thm:toric-examples}.

\subsection{Set up for toric packing}
We start with a definition.

\begin{Definition}\label{def:packedtriangle}
For $\lambda>0$, let
\[
B(\lambda) := \big\{(x,y)\in\R^2\mid x\geq 0, \,y\geq 0,\, x+y<\lambda\big\}.
\]
We call $B(\lambda)$ the \emph{model triangle of size $\lambda$}.
For a Delzant polygon $\De$ and $p$ a vertex of $\De$, we say that a set $B\subset \De$ is an \emph{$($admissibly$)$ packed triangle at $p$} if there exists $M \in \SL$ and some~$\lambda>0$ such that
\[B = M(B(\lambda)) + p\]
and furthermore that $B$ and $\De$ are equal in a small neighborhood of $p$.
\end{Definition}

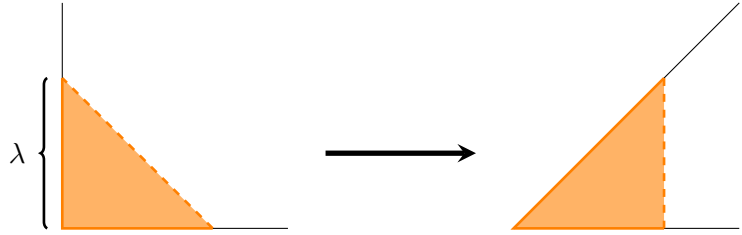
\begin{figure}
 \centering
 \begin{tikzpicture}
 \draw (0,0) -- (0,3);
 \draw (0,0) -- (3,0);
 \draw [line1,line width = 1,fill = fill1] (2,0) -- (0,0) -- (0,2);
 \draw [dashed, line1,line width = 1] (2,0) -- (0,2);
 \draw (6,0) -- (9,3);
 \draw (6,0) -- (9,0);
 \draw [line1,line width = 1,fill = fill1] (8,0) -- (6,0) -- (8,2);
 \draw [dashed, line1,line width = 1] (8,2) -- (8,0);
 \draw [-stealth, line width = 2] (3.5,1) -- (5.5,1);
 \draw [decorate, line width=1pt, decoration={brace}] (-0.2,0)--(-0.2,2);
 \node (n1) at (-0.6,1) {$\lambda$};
 \end{tikzpicture}
 \caption{Packing the model triangle $B(\lambda)$ (left) into a corner of a Delzant polygon (right) by an element of $\SL$.}
 \label{fig:packing_ex}
\end{figure}

Notice that $B(\lambda)$ includes exactly two out of the three segments of its boundary, and also notice that $B$ has to fit ``perfectly into the corner'', see Figure~\ref{fig:packing_ex}.

\begin{Definition}\label{def:DPacking}
A \emph{packing} $P$ of a Delzant polygon $\De$ is the union of a set of pairwise disjoint admissibly packed triangles, i.e.,
\[P = \bigcup_{i=1}^d B_i,\]
where each $B_i$ is an admissibly packed triangle in $\De$ and $B_i \cap B_j = \varnothing$
if $i \neq j$.
Given a Delzant polygon $\De$, let $\PT(\De) = \{P \mid P \text{ is a packing of } \De \}$.
\end{Definition}

 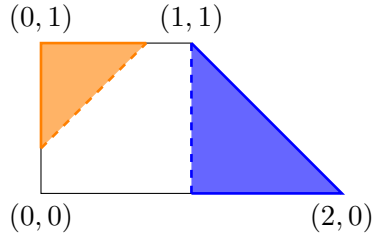
\begin{figure}
 \centering
 \begin{tikzpicture}[scale=2]
 \draw (0,0)--(2,0)--(1,1)--(0,1)--cycle;

 \draw [line1,line width = 1,fill = fill1] (0,.3) -- (0,1) -- (.7,1);
 \draw [\triangleline, line1,line width = 1] (0,.3) -- (.7,1);
 \draw [line2,line width = 1,fill = fill2] (1,1) -- (2,0) -- (1,0);
 \draw [\triangleline, line2,line width = 1] (1,1) -- (1,0);
 \node at (0,0) [anchor = north] {$(0,0)$};
 \node at (2,0) [anchor = north] {$(2,0)$};
 \node at (1,1) [anchor = south] {$(1,1)$};
 \node at (0,1) [anchor = south] {$(0,1)$};

 \end{tikzpicture}
 \caption{A packing of a Delzant polygon by two triangles.}
 \label{fig:delz_pack}
 \end{figure}

Since a Delzant polygon has only finitely many vertices, we can pack at most finitely many triangles inside a Delzant polygon, which is why we have assumed that the union is finite in Definition~\ref{def:DPacking}.
An example of a Delzant polygon packed with two triangles is shown in Figure~\ref{fig:delz_pack}.

For any measurable subset $S\subset \R^2$, let $\area(S)$ denote its usual Lesbegue measure.

\begin{Definition}\label{def:toric-density}
The \emph{packing density} of a Delzant polygon $\De$ is denoted by $\rho_{{\rm T}}(\De)$ and given by
\[
 \rho_{{\rm T}}(\De): = \sup_{P \in \PT(\De)}\left( \frac{\area(P)}{\area(\De)}\right).
\]
\end{Definition}

Given any Delzant polygon $\De$, notice that $P\subset \De$ for any packing $P$, and thus $0<\rho_{{\rm T}}(\De)\leq 1$.
We call any packing with density equal to 1 a ``perfect packing''. Pelayo~\cite{Pe2006} showed that squares and right isosceles triangles are the only Delzant polygons (in dimension $4$) which admit a perfect packing.
We discuss perfect semitoric packings in Section~\ref{sec:perfect}.

We endow the space of possible packings $\PT(\De)$ with the topology pulled back from $\PTtilde(\De)$ via the bijection $\Lambda$ described in Proposition~\ref{prop:toric}.
Since $\PTtilde(\De)$ is compact, $\PT(\De)$ is also compact, and therefore
a packing which achieves the upper bound always exists.

\begin{Definition}
A packing $P_0 \in \PT(\De)$ is a \emph{maximal packing} of $\De$ if \[\area(P_0) = \sup_{P \in \PT(\De)} \area(P).\]
\end{Definition}

The packing density is the proportion of the polygon covered by triangles in a maximal packing.
A Delzant polygon can have several distinct maximal packings, but only finitely many, as shown by~\cite[Proposition 1.4]{PeSc2008}.

\subsection[SL2Z-length]{\texorpdfstring{$\boldsymbol{\SL}$}{SL2Z}-length}
Here we introduce the notion of $\SL$-length, which is the natural way to measure the length of the edges of a Delzant polygon.
It will also be the natural language in which to discuss packings.
The results in this section are not new, but we could not find a good reference for their proofs in the literature, so we develop them here.

The following lemma is well known, but we provide a proof for completeness.

\begin{Lemma}\label{lem:SL2Z-vect}
 Let $v\in \R^2$ be a vector which is either vertical or has rational slope. Then there exists a matrix $M\in\SL$ such that $Mv$ is in the positive $x$-axis $($i.e.,~$Mv\in \R_{\geq 0}\times \{0\})$.
 Moreover, there exists an $\ell\geq 0$ such that for any such $M$, $Mv = \ell \left(\begin{smallmatrix}1 \\ 0\end{smallmatrix}\right)$.
\end{Lemma}

\begin{proof}
The case that $v$ is the zero vector is trivial, so we assume that $v$ is not the zero vector for the duration of the proof.
Let $w=\vect{w_1}{w_2}\in\Z^2$ be a primitive vector in the same direction as~$v$.
Thus, $\mathrm{gcd}(w_1,w_2)=1$ and $v = \ell w$ for some $\ell>0$.

By B\'{e}zout's theorem, there exist integers $x$, $y$ such that $xw_1 + yw_2 = \gcd (w_1,w_2) = 1$. Let~$M = \begin{psmallmatrix} x & y \\ -w_2 & w_1 \end{psmallmatrix}$. Then~$\det(M) = xw_1 +yw_2 =1$, so $M \in \SL$, and
\[
Mw =\begin{pmatrix}
xw_1+yw_2 \\ -w_2w_1+w_1w_2 \end{pmatrix} = \begin{pmatrix} 1\\ 0\end{pmatrix}.
\]
Thus, $Mv = \ell \left(\begin{smallmatrix}  1\\ 0\end{smallmatrix} \right)$, so $M$ is the required matrix.

Suppose that there is another matrix $M' \in \SL$ such that $M' v = \ell' \vect{1}{0}$ for some $\ell'>0$. We will show that $\ell = \ell' $. Let $B = M'M^{-1}$ and denote its entries by $B_{ij}$. Then
\begin{align*}\label{eq:sl_h2}
 B \begin{pmatrix} 1\\ 0\end{pmatrix} = M' M^{-1}\begin{pmatrix} 1\\ 0\end{pmatrix} = M' \left( \frac{1}{\ell} v \right) = \frac{\ell'}{\ell} \begin{pmatrix} 1\\ 0\end{pmatrix} ,
\end{align*}
which implies that $B_{11} = \frac{\ell'}{\ell}$ and $B_{21} = 0$. Since $B \in \SL$, we know that $\det(B) = 1$ so~$B_{11} B_{22} = 1$. Since $B_{11}, B_{22}\in \mathbb{Z}$, we obtain that $B_{11} = \pm 1$ so $\ell = \ell'$, since both are positive.
\end{proof}

An important part of Lemma~\ref{lem:SL2Z-vect} is that $\ell\geq 0$ is independent of the choice of $M$, and thus the following notion is well defined.

\begin{Definition}\label{def:sl2z-length} Let $v\in\R^2$ be a vector which is either vertical or has rational slope.
The \emph{$\SL$-length of $v$} is the unique $\ell\geq 0$
such that $Mv = \ell \vect{1}{0}$
for some $M\in\SL$.
If $L\subset \R^2$ is a line segment which is either vertical or has rational slope then the \emph{$\SL$-length of $L$} is the $\SL$-length of $v := q-p$, where $p$ and $q$ are the end points of $L$.
\end{Definition}

By definition, $\SL$-length is preserved under the action of $\SL$.
Notice that $\SL$-length is not the same as the Euclidean length in $\R^2$.
For instance, the vector $\begin{psmallmatrix} 2 \\ 1 \end{psmallmatrix}$ has $\mathrm{SL}_2(\mathbb{Z})$ length $1$.

\begin{Corollary} \label{cor:SL_gcd}
Let $ x = \vect{x_1}{x_2} \in \mathbb{Z}^2$. Then the $\SL$-length of $x$ is $\mathrm{gcd}(x_1, x_2)$.
\end{Corollary}

In practice, Corollary \ref{cor:SL_gcd} is sufficient for computing the $\SL$-lengths of the sides of a~rational polygon (such as Delzant or semitoric polygons) since it is possible to scale the entire polygon such that every vertex lies on an integer point and thus all the vectors directing the edges have integer entries. This does not affect the packing density, since it is invariant under such scaling.

\begin{Remark}
 Let $\De$ be a Delzant polygon and $P\in \PT(\De)$ be a toric packing.
 Then for any $A\in\SL$, it can be verified that $A(\De)$ is a Delzant polygon and $A(P)\in\PT(A(\De))$, and furthermore that $\area(\De) = \area(A(\De))$
 and $\area(P) = \area(A(P))$. Thus, we conclude that the packing density is invariant under the action of $\SL$.
\end{Remark}

\subsection{Computing the toric packing density}
\label{sec:toric-computing}

In this section, we explain how to transform the packing problem to a quadratic optimization problem over a set of linear constraints.
This result already appeared in~\cite[Proposition 1.4]{PeSc2008}, but we prove it here for completeness.

A \emph{compact convex polytope} in $\R^d$ is a closed bounded set in $\R^d$ which is the intersection of a~finite number of closed half spaces.

\begin{proof}[Proof of Proposition~\ref{prop:toric}]
Let $\De$ be a Delzant polygon. Label the vertices of $\De$ by $p_1,\ldots,p_d$, arranged in clockwise order, let $e_i$ denote the edge connecting $p_i$ to $p_{i+1}$,
and let $\ell_i\geq 0$ be the $\SL$-length of $e_i$, for $i=1,\ldots, d$.
We assume that $p_1$ is the lexicographically minimal (i.e.,~bottom left) vertex.
For simplicity, we take these indices cyclically, so $p_{d+1} = p_1$, for instance.
We will show that $\Lambda\colon \PT(\De)\to\R^d$ as defined in equation~\eqref{eqn:Lambda} is a bijection onto the set $\PTtilde(\De)$, as defined in equation~\eqref{eqn:tildePT}.

First, we will show that $\Lambda(\PT(\De))\subset \PTtilde(\De)$.
Let $P = \bigcup_{i=1}^d B_i$ be a packing where $B_i$ is a~triangle packed at the vertex $p_i$ (and some $B_i$ may be empty).
Let $(\lambda_1,\ldots, \lambda_d):=\Lambda(P)$ using the definition of the map $\Lambda$ from equation~\eqref{eqn:Lambda}, which means that $\lambda_i$ is the $\SL$-size of $B_i$, taking $\lambda_i=0$ if $B_i=\varnothing$.

By acting on $\De$ and $P$ by an element of $\SL$ and a translation, we may assume that $p_i = (0,0)$ and $p_{i+1} = (0,\ell_i)$. Now, it is clear to see that the two triangles packed at $v_i$ and $v_{i+1}$ do not intersect if and only if the sum of their side lengths do not exceed $\ell_i$. The $\SL$-length of their sides are $\lambda_i$ and $\lambda_{i+1}$, respectively, and $\SL$-length is invariant under the action of~$\SL$. Since the $\SL$-length of a vertical line segment is equal to its actual length, we conclude that $\lambda_i + \lambda_{i+1} \leq \ell_i$ for $i\in\{1,\ldots,d\}$,
see Figure~\ref{fig:toric-proof}.
The side length of a triangle cannot be negative, so $\lambda_i \geq 0$ for $i \in \{1,2,\ldots,d\}$. Thus, $\Lambda(P) \in \PTtilde(\De)$.

\begin{figure}[t] \centering
 \begin{tikzpicture}
 \begin{scope}[shift={(-2,0)}]
 \draw (3,0) -- (0,0) -- (3,3) -- (5,3);
 \draw [line1,line width = 1,fill = fill1] (1.5,0) -- (0,0) -- (1.5,1.5);
 \draw [dashed, line1,line width = 1] (1.5,1.5) -- (1.5,0);
 \draw [line2,line width = 1,fill = fill2] (2,2) -- (3,3) -- (4,3);
 \draw [dashed, line2,line width = 1] (4,3) -- (2,2);
 \node at (1.5,2) {$e_i$};
 \end{scope}
 \draw [-stealth, line width = 2] (3.5,1) -- (5.5,1);
 \begin{scope}[shift = {(8,0)}]
 \draw (3,0) -- (0,0) -- (0,3) -- (2,3);
 \draw (0,0) -- (3,0);
 \draw [line1,line width = 1,fill = fill1] (1.5,0) -- (0,0) -- (0,1.5);
 \draw [dashed, line1,line width = 1] (1.5,0) -- (0,1.5);
 \draw [line2,line width = 1,fill = fill2] (0,2) -- (0,3) -- (1,3);
 \draw [dashed, line2,line width = 1] (1,3) -- (0,2);
 \draw [decorate, line width=1pt, decoration={brace}] (-0.2,0)--(-0.2,1.5);
 \node at (-0.7,.75) {$\lambda_i$};
 \draw [decorate, line width=1pt, decoration={brace}] (-0.2,2)--(-0.2,3);
 \node at (-0.7,2.5) {$\lambda_{i+1}$};
 \draw [decorate, line width=1pt, decoration={brace}] (-1.2,0)--(-1.2,3);
 \node at (-1.6,1.5) {$\ell_i$};
 \end{scope}
 \end{tikzpicture}
 \caption{In the proof of Proposition~\ref{prop:toric}, we use an element of $\SL$ to send the edge $e_i$ onto the $y$-axis, and then it is clear that the two adjacent packed triangles are disjoint if and only if the sum of their $\SL$-sizes is less than the $\SL$-length of the edge between them. That is, we have $\lambda_{i}+\lambda_{i+1}\leq \ell_i$.} \label{fig:toric-proof}
\end{figure}
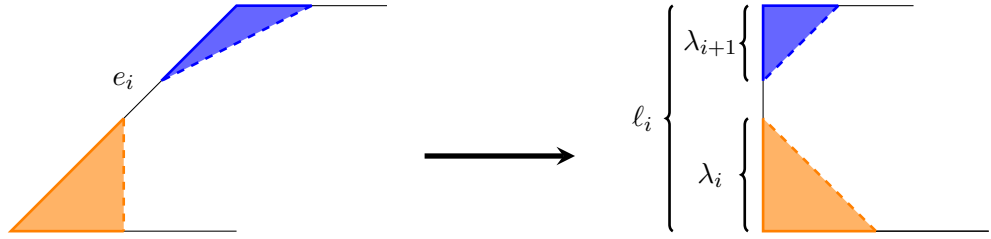

We will now show $\Lambda$ is bijective.
It is clear that $\Lambda$ is injective, as two different packings will differ by at least one admissibly packed triangle $B_i$ and the only difference is the $\SL$-size~$\lambda_i$. To show it is surjective, let $(\lambda_1,\ldots,\lambda_d) \in \PTtilde(\De)$. We claim that there exists a packing $P$ in which the triangle $B_i$ packed at the vertex~$p_i$ has $\SL$-size $\lambda_i$. Notice that the inequalities imply that $\lambda_i \leq \ell_i$ and $\lambda_i \leq \ell_{i-1}$, so a $B_i$ of such a size can be admissibly packed at~$p_i$. As the argument in the previous paragraph, we know that $\lambda_i+\lambda_{i+1}\leq \ell_i$ implies that adjacent triangles do not intersect, and convexity guarantees non-adjacent triangles will not intersect either.
Thus, such a $P\in\PT(\De)$ exists, so $\Lambda$ is surjective.
Hence, we conclude that $\Lambda$ is a bijection from $\PT(\De)$ to $\PTtilde(\De)$.

We will now prove equation~\eqref{eqn:LambdaP-magnitude}.
The area of the model triangle $B(\lambda)$ is $\frac{1}{2}\lambda^2$, and area is invariant under the action of $\SL$ and translation.
Thus, any triangle of $\SL$-size equal to~$\lambda$ has area $\frac{1}{2}\lambda^2$,
and given any packing $P = \bigcup_{i=1}^d B_i$ where the $\SL$-size of $B_i$ is $\lambda_i$, we obtain
\[
 \area(P) = \sum_{i=1}^d \area(B_i) = \sum_{i=1}^d \frac{1}{2} \lambda_i^2 = \frac{1}{2} || \Lambda(P) ||^2,
\]
using the fact that the $B_i$ are disjoint.

To conclude that $\PTtilde(\De)$ is a compact convex polytope, we simply notice that $\PTtilde(\De)$ is defined by a finite number of non-strict linear inequalities. It is bounded because it is a subset of the bounded set $[0,\ell]^d\subset \R^d$, where $\ell = \mathrm{max}\{ \ell_i \mid i = 1,\ldots d\}$.
\end{proof}

Let $\De$ be a Delzant polygon. Then applying Proposition~\ref{prop:toric} to the definition of the packing density (see Definition~\ref{def:toric-density}), we have
\[
 \rho_{{\rm T}}(\De) = \sup \left\{ \frac{\area(P)}{\area(\De)}\colon  P \in \PT(\De) \right\} = \sup \left\{ \frac{||q||^2}{2 \area(\De)}\colon  q \in \PTtilde(\De) \right\}.
\]

Since $\PTtilde(\De)$ is a compact subset of $\R^d$, and $q\mapsto ||q||^2$ is continuous, the supremum is actually a maximum. Our goal now is to reduce the number of candidates for the maximum down to a~finite set. This set will be the vertices of the polytope $\PTtilde(\De)$, as we now see.

Let $\mathcal{P}$ be a compact convex polytope in $\R^d$. A point $p\in\mathcal{P}$ is called \emph{extreme} for $\mathcal{P}$ if it cannot be written in the form $ p = t q + (1- t) r$ for any $q,r \in \mathcal{P}$ and $t \in (0,1)$. In this case, $p\in\mathcal{P}$ being extreme is equivalent to $p$ being a $0$-dimensional face of $\mathcal{P}$, otherwise known as a \emph{vertex} of $\mathcal{P}$.

\begin{Proposition}
\label{prop:vert}
 Let $\mathcal{P}\subset \R^d$ be a compact convex polytope of dimension $d$. Then the function $q\mapsto ||q||^2$, for $q\in\R^d$, achieves its maximum on a vertex of $\mathcal{P}$.
\end{Proposition}

\begin{proof}
 Define $f\colon \mathcal{P}\to \R^d$ by $f(q) = ||q||^2$.
 The idea of this proof is that if $p\in\mathcal{P}$ is not a~vertex, then there is a line segment through $p$ which is contained in $\mathcal{P}$, and we can compute that the maximum of $f$ on this line segment does not occur at $p$, and therefore the maximum of $f$ in $\mathcal{P}$ cannot occur at $p$.

 Suppose that $p$ is not a vertex, and thus is not extreme. Then
 there exist $q,r \in \mathcal{P}$ and $t_0 \in (0,1)$ such that $p = (1-t_0)q +t_0r$. Since $\mathcal{P}$ is convex and $q,r \in \mathcal{P}$, the line segment $L$ from~$q$ to $r$ is contained in $\mathcal{P}$.
 Let $v$ be the vector from $q$ to $r$, so that points in $L$ can be written as $q+tv$ for $t\in[0,1]$.
 Notice that
 \[
 \frac{\mathrm{d}^2}{\mathrm{d}t^2}\big( f(q+tv)\big)|_{t=t_0} = 2||v||^2>0,
\]
 so $p$ is not the maximum of $f$ restricted to $L\subset \mathcal{P}$, and thus $p$ is not the maximum of $f$ in $\mathcal{P}$.

 We have concluded that $f$ cannot take its maximum value on $\mathcal{P}$ at any non-vertex $p\in\mathcal{P}$, but since $\mathcal{P}$ is compact $f$ must take its maximum somewhere. We conclude that $f$ achieves its maximum at a vertex of $\mathcal{P}$.
 \end{proof}

Applying Proposition \ref{prop:vert} to the polytope
$\PTtilde(\Delta)$ from Proposition~\ref{prop:toric},
we obtain the following:
\begin{Corollary}
 \label{cor:toric-vertex}
 Let $\De$ be a Delzant polygon.
 Then
 \[
 \rho_{{\rm T}} (\De) = \max \left\{ \frac{||q||^2}{2 \area(\De)}\colon  q\in V\big(\PTtilde(\De)\big)\right\},
 \]
 where $V\big(\PTtilde(\De)\big)$ is the finite set of vertices of the compact convex polytope $\PTtilde(\De)$, given in equation~\eqref{eqn:tildePT}.
 \end{Corollary}

 \subsection{A program to calculate the toric packing density}
\label{sec:program}

We have produced a python program, available online at~\cite{alg-url}, which takes as an input the vertices of a Delzant polygon and gives as an output its maximal packing and maximal packing density.
The key idea that allows the development of such an algorithm is that Proposition~\ref{prop:toric} and Corollary~\ref{cor:toric-vertex} reduce the number of candidates for the maximum packing down to a finite list, which thus can be checked by a computer.
Roughly, the algorithm follows these steps:
\begin{enumerate}\itemsep=0pt
 \item[(1)] take as an input the vertices of the Delzant polygon $\De$;
 \item[(2)] calculate the $\SL$-lengths of the edges of $\De$ using Corollary~\ref{cor:SL_gcd};
 \item[(3)] find the vertices $v_1,\ldots,v_k\in\R^d$ of $\PTtilde(\De)$;
 \item[(4)] return the largest of $||v_1||^2,\ldots,||v_k||^2$ divided by $2\area(\De)$.
\end{enumerate}
By Corollary~\ref{cor:toric-vertex}, the value returned by this procedure is the packing density of $\De$.

In many cases, such as when the vertices of the polygon lie on rational points, the Delzant polygon can be scaled so that its vertices lie on the integer lattice, and this scaling does not affect the packing density.
Therefore, in these cases each number in our program~\cite{alg-url} is rational and stored exactly as a fraction or integer, so no information is lost due to floating-point representation.


\subsection{The alternating edge condition}
Corollary~\ref{cor:toric-vertex} is useful for having a computer calculate the packing density of a specific Delzant polygon, but in practice $\PTtilde(\De)$ can have many vertices and it is not always easy to find the vertices of this polytope in the case that $\De$ depends on parameters.

In some cases, though, there are useful shortcuts to computing the packing density. Let $\De$ be a Delzant polygon with $d=2k$ edges for some $k\in\Z_{>0}$. We say $\De$ satisfies the \emph{alternating edge condition} if the edges can be labeled in clockwise or counterclockwise order as $e_1,\ldots,e_d$ such that $\ell_{2j-1} \leq \ell_{2j}$ for $j\in\{1,\ldots,k\}$, where $\ell_i$ is the $\SL$-length of $e_i$.

 \begin{Proposition}\label{prop:alt}
 Let $\De$ be a Delzant polygon that satisfies the alternating edge condition and let $\ell_1, \ell_2,\dots,\ell_{2d}$ be the $\SL$-length of its edges. Then
 \[\rho_{{\rm T}}(\De) = \frac{1}{2\area(\De)}\sum_{j=1}^k \ell_{2j-1}^2.\]
 \end{Proposition}

\begin{proof}
Let $\PTtilde(\De)$ be as in equation~\eqref{eqn:tildePT} and let $q=(\lambda_1,\ldots,\lambda_d)\in\PTtilde(\De)$.
Then
\begin{equation}\label{eqn:alt-edge-proof}
 ||q||^2
 = \sum_{i=1}^d \lambda_i^2
 = \sum_{j=1}^k \big(\lambda_{2j-1}^2 + \lambda_{2j}^2\big)
 \leq \sum_{j=1}^k (\lambda_{2j-1} + \lambda_{2j})^2
 \leq \sum_{j=1}^k \ell_{2j-1}^2.
\end{equation}
Let $P\in\PT(\De)$ be a maximal packing of $\De$. By Proposition~\ref{prop:toric}, $\Lambda(P)\in\PTtilde(\De)$ and $\area(P) = \frac{1}{2}||\Lambda(P)||^2$, so equation~\eqref{eqn:alt-edge-proof} implies that
 \[\rho_{{\rm T}}(\De) = \frac{\area(P)}{\area(\De)} \leq \frac{1}{2\area(\De)} \sum_{j=1}^k \ell_{2j-1}^2.\]

Now, we show that we can achieve this bound. Let $q := (0,\ell_1,0,\ell_3,\dots,0,\ell_{d-1})$. We observe that the entries of $q$ are non-negative, and for any $j\in\{1,\ldots,k\}$ we have that \[\lambda_{2j-1} + \lambda_{2j} = 0 + \ell_{2j-1} \leq \ell_{2j-1}\qquad\textrm{ and }\qquad \lambda_{2j} + \lambda_{2j+1} = \ell_{2j-1} + 0 \leq \ell_{2j},\] with the last inequality by the assumption on $\De$. This implies that $q \in \PTtilde(\De)$.
Thus, by Proposition~\ref{prop:toric}, there exists some packing $P\in\PT(\De)$ such that $\Lambda(P) = q$,
and
\[
 \frac{\area(P)}{\area(\De)} = \frac{1}{2\area(\De)}||q||^2 = \frac{1}{2\area(\De)}\sum_{j=1}^d \ell_{2j-1}^2,
\]
as desired.
\end{proof}

Notice that we defined the alternating edge condition in the situation that the odd indexed edges are shorter than the even indexed edges, but the case in which the even indexed edges are the shorter ones can be dealt with by shifting the labels of the edges.


\subsection{Examples}

Now we compute the packing density of the examples from Definition~\ref{def:delz-examples}.

\subsubsection{The Delzant triangle}
Let $\De$ be the Delzant triangle with parameter $a>0$, as in Figure~\ref{fig:del-ex-triangle}.
Notice that packing a~single triangle of $\SL$-size $a$ at any of the three vertices of $\De$ will achieve the upper bound of 1 for the packing density.
Thus, if $\De$ is the Delzant triangle with parameter $a$, then
\begin{equation}\label{eqn:pack-density-triangle}
 \rho_{{\rm T}}(\De) = \frac{\frac{1}{2}a^2}{\area(\De)} = 1,
\end{equation}
which was already known to Pelayo~\cite{Pe2006}.

\subsubsection{The Rectangle} Next we consider the rectangle with parameters $a,b>0$, as in Figure~\ref{fig:del-ex-rect}. The $\SL$-lengths of the edges are $a$, $b$, $a$, $b$, so in both possible cases of $a\geq b$ or $b\geq a$, it satisfies the alternating edge condition. By Proposition~\ref{prop:alt}, we conclude that if $\De$ is the rectangle with parameters~$a,b>0$ then
\begin{equation}\label{eqn:pack-density-rect}
 \rho_{{\rm T}}(\De) = \frac{2(\min\{a,b\})^2}{2\area(\De)} = \frac{(\min\{a,b\})^2}{ab}=\frac{\min\{a,b\}}{\max\{a,b\}}.
\end{equation}
Notice that if $a=b$ in this case we obtain $\rho_{{\rm T}}(\De) =1$, which was already known to Pelayo~\cite{Pe2006}.

\subsubsection{The Hirzebruch trapezoid}
\label{sec:hirz-packing}
Let $\De$ be the Hirzebruch trapezoid of parameters $a,b>0$ and $n\in \Z_{>0}$, as in Definition~\ref{def:delz-examples} and Figure~\ref{fig:del-ex-hirz}.
The $\SL$-lengths of the edges are $a$, $b$, $a$, $an+b$, and thus $\lambda = (\lambda_1,\lambda_2,\lambda_3,\lambda_4)\in\PTtilde(\De)$ if and only if $\lambda_i\geq 0$ for all $i\in\{1,2,3,4\}$ and
\begin{equation}\label{eqn:hirz-inequalities}
 \lambda_1+\lambda_2 \leq a,\qquad
 \lambda_2+\lambda_3 \leq b,\qquad
 \lambda_3+\lambda_4 \leq a,\qquad
 \lambda_4+\lambda_1 \leq an+b.
\end{equation}
For such a $\lambda$, we have that
\begin{equation}\label{eqn:hirz-upper}
 ||\lambda||^2 = \big(\lambda_1^2 + \lambda_2^2\big) + \big(\lambda_3^2 + \lambda_4^2\big) \leq (\lambda_1 + \lambda_2)^2 + (\lambda_3 + \lambda_4)^2 \leq a^2 + a^2 = 2a^2.
\end{equation}
The upper-bound in equation~\eqref{eqn:hirz-upper} is achievable in the case that $n>1$ by taking $\lambda = (a, 0, 0, a)$, and it can easily be checked that $(a,0,0,a)\in\PTtilde(\De)$.

If $n=1$ there are several cases.

{\it Case $1$: $n=1$ and $b\geq a$.} In this case the $\SL$-lengths of the edges are $a$, $b$, $a$, $a+b$. Since $b\geq a$ the polygon satisfies the alternating edge condition, and by Proposition~\ref{prop:alt}, $\lambda = (0,a,0,a)\in\PTtilde(\De)$ has the maximal magnitude which is $2a^2$.

{\it Case $2$: $n=1$ and $b\leq a/2$.} Since $\lambda_1\leq a$, $\lambda_4\leq a$, $a>b$, and $\lambda_1+\lambda_4\leq a+b$ we conclude that $\lambda_1^2+\lambda_4^2 \leq a^2 + b^2$.
Thus we have that
\[
 \lambda_1^2 + \lambda_2^2 + \lambda_3^2 + \lambda_4^2 = \big(\lambda_2^2 + \lambda_3^2\big) + \big(\lambda_1^2 + \lambda_4^2\big) \leq a^2 + 2b^2.
\]
Since $b\leq a/2$, $(a,0,b,b)\in\PTtilde(\De)$ will give a valid packing that achieves this bound.

{\it Case $3$: $n=1$ and $a/2 < b < a$.} Since $b>a/2$ the packing $(a,0,b,b)$ from the previous case is not possible, so we compute the maximal magnitude of a vector in $\PTtilde(\De)$ using Corollary~\ref{cor:toric-vertex}. Suppose that $\lambda\in\PTtilde(\De)$ is a vector of maximal magnitude in $\PTtilde(\De)$, and we will show that $||\lambda||^2 = a^2+b^2+(a-b)^2$. Corollary~\ref{cor:toric-vertex} states that $\lambda$ is a vertex of $\PTtilde(\De)$. Thus, at least four of the eight inequalities defining $\PTtilde(\De)$ must actually evaluate to equalities at $\lambda$.

It is clear that if $\lambda$ is such that $||\lambda||^2$ is maximal then $\lambda_1, \lambda_4>0$, and furthermore that we do not have simultaneously that both $\lambda_2=0$ and $\lambda_3=0$.
Thus, there are two subcases:
\begin{itemize}\itemsep=0pt
 \item Suppose that $\lambda_2=0$ and $\lambda_3\neq 0$.
Then for $\lambda$ to be a vertex of $\PTtilde(\De)$ we must have at least three of the four equations in the system~\eqref{eqn:hirz-inequalities} as equalities.
The only cases in which this is possible are $\lambda = (a,0,a-b,b)$ or $\lambda = (a,0,b,a-b)$, both of which have magnitude~$a^2+b^2+(a-b)^2$.

\item The other case is that $\lambda_3=0$ and $\lambda_2\neq 0$. Similarly, the only vertices with these values are~$\lambda = (b, a-b, 0, a)$ and $\lambda = (a-b, b, 0, a)$, which again both have magnitude $a^2+b^2+(a-b)^2$.
\end{itemize}

Putting all of these cases together, using the fact that $\area(\De) = \frac{a}{2}(na+2b)$ in this case, and applying Proposition~\ref{prop:toric}, we obtain
\begin{equation}\label{eqn:pack-density-hirz}
 \rho_{{\rm T}}(\De) = \begin{cases}
 \dfrac{2a}{na+2b}, & \textrm{if }n>1,\vspace{1mm}\\
 \dfrac{2a}{a+2b}, & \textrm{if }n=1 \textrm{ and }a\leq b,\vspace{1mm}\\
 \dfrac{a^2+b^2+ (a-b)^2}{(a+2b)a}, & \textrm{if }n=1 \textrm{ and } b<a<2b,\vspace{1mm}\\
 \dfrac{a^2+2b^2}{(a+2b)a}, & \textrm{if }n=1 \textrm{ and } a \geq 2b. \end{cases}
\end{equation}

\begin{figure}[t]
 \centering
 \begin{subfigure}{.32\textwidth}
 \centering
 \begin{tikzpicture}[scale = .55]
 \draw (0,0)--(0,4)--(5,4)--(9,0)--cycle;
 \draw [line1,line width = 1,fill = fill1] (0,4) -- (0,0) -- (4,0);
 \draw [\triangleline, line1,line width = 1] (4,0) -- (0,4);
 \draw [line2,line width = 1,fill = fill2] (5,0) -- (9,0) -- (5,4);
 \draw [\triangleline, line2,line width = 1] (5,4) -- (5,0);
 \end{tikzpicture}
 \caption{$b>a$}
 \end{subfigure}\,
 \begin{subfigure}{.32\textwidth}
 \centering
 \begin{tikzpicture}[scale = .55]
 \draw (0,0)--(0,4)--(3,4)--(7,0)--cycle;
 \draw [line1,line width = 1,fill = fill1] (0,4) -- (0,0) -- (4,0);
 \draw [\triangleline, line1,line width = 1] (4,0) -- (0,4);
 \draw [line2,line width = 1,fill = fill2] (4,0) -- (7,0) -- (4,3);
 \draw [\triangleline, line2,line width = 1] (4,3) -- (4,0);
 \draw [line3,line width = 1,fill = fill3] (2,4) -- (3,4) -- (4,3);
 \draw [\triangleline, line3,line width = 1] (4,3) -- (2,4);
 \end{tikzpicture}
 \caption{$\frac{a}{2}<b<a$}
 \end{subfigure}
 \begin{subfigure}{.25\textwidth}
 \centering
 \begin{tikzpicture}[scale = .55]
 \draw (0,0)--(0,4)--(1,4)--(5,0)--cycle;
 \draw [line1,line width = 1,fill = fill1] (0,4) -- (0,0) -- (4,0);
 \draw [\triangleline, line1,line width = 1] (4,0) -- (0,4);
 \draw [line2,line width = 1,fill = fill2] (4,0) -- (5,0) -- (4,1);
 \draw [\triangleline, line2,line width = 1] (4,1) -- (4,0);
 \draw [line3,line width = 1,fill = fill3] (0,4) -- (1,4) -- (2,3);
 \draw [\triangleline, line3,line width = 1] (2,3) -- (0,4);
 \end{tikzpicture}
 \caption{$b<\frac{a}{2}$}
 \end{subfigure}
 \caption{Various cases of a maximal packing of the Hirzebruch surface with $n=1$, depending on the relative values of the parameters $a>0$ and $b>0$.}
 \label{fig:hirz-pack}
\end{figure}
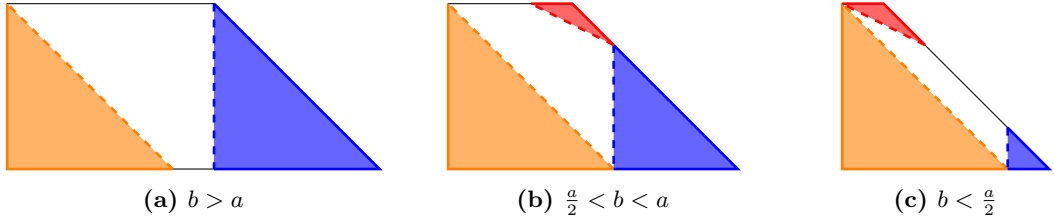

In Figure~\ref{fig:hirz-pack}, we show maximal packings of the Hirzebruch trapezoid for various choices of parameters.
Theorem~\ref{thm:toric-examples} now follows from equations~\eqref{eqn:pack-density-triangle}, \eqref{eqn:pack-density-rect}, and~\eqref{eqn:pack-density-hirz}.

\section{Semitoric packing}\label{sec:semitoric}

Packing semitoric polygons is a much more subtle problem than packing Delzant polygons, but we will see in this section that we can carefully adapt our methods from the previous section to the semitoric case.

\subsection{Set up for semitoric packing}
\label{sec:semitoric-setup}

Recall that a semitoric polygon is an infinite family of polygons, so the analogue of a ``subset'' in this context is a subset of each representative of the semitoric polygon.
The idea of an admissibly packed triangle in a semitoric polygon is that \textit{in at least one representative of the semitoric polygon} it should satisfy exactly the same conditions as the Delzant case from Definition~\ref{def:DPacking}.
The definition of an admissibly packed triangle in a semitoric polygon is obtained from the study of embedded 4-balls in semitoric systems which are equivariant with respect to a certain local 2-torus action, as described in~\cite{FiPaPe2016}.
As in the previous section, we will not focus on the motivation here but instead we will simply describe the conditions at the level of polygons.

Recall that a semitoric polygon $\stpoly$ is defined as an orbit of the action of the group~$\mathcal{T}\times G_m$ on a semitoric polygon representative $\strep$, see Definition~\ref{def:semitoric-polygon}.
The idea of the following definition is that a subset of a given semitoric polygon representative is an admissibly packed triangle if we can pass to an equivalent representative (by the action of $\mathcal{T}\times G_m$) and in this equivalent representative the corresponding subset avoids the marked points and also satisfies the same conditions as in the Delzant case from Definition~\ref{def:packedtriangle}.

\begin{Definition}\label{def:STadm-packed-triangle}
Let $(\De,\vec c, \vec \epsilon)$ be a semitoric polygon representative with $m$ marked points, and let $p\in\Delta$ be a vertex. We say $B \subset \De$ is an \emph{admissibly packed triangle at $p$} if there exists a~$\sigma \in \mathcal{T}\times G_m$ such that $\sigma(B) \subset \sigma(\De)$ satisfies:
\begin{enumerate}\itemsep=0pt
 \item[(1)] $\sigma(p)\in\sigma(\De)$ is a vertex which satisfies the Delzant corner condition;
 \item[(2)] $\sigma(c_i)\notin \sigma(B)$ for all $i=1,\ldots,m$;
 \item[(3)]\label{item:lambda} there exists $M\in\SL$ and some $\lambda>0$ such that
 \[
 \sigma(B) = M(B(\lambda)) + \sigma(p),
 \]
 where $B(\lambda)$ is the model triangle of size $\lambda$, as in Definition~\ref{def:DPacking};
 \item[(4)] $\sigma(B)$ and $\sigma(\De)$ are equal in a small neighborhood of $\sigma(p)$.
\end{enumerate}
The number $\lambda>0$ from item~\eqref{item:lambda} is called the \emph{$\SL$-size of $B$}.
\end{Definition}

Let $\strep$ be a semitoric polygon representative and let $B\subset \De$ be an admissibly packed triangle at a vertex $p\in\De$.
If $(\De',\vec{c}\,{}',\vec{\epsilon}\,{}')\in\stpoly$ is another representative of the same semitoric polygon, then there exists some $\sigma \in \mathcal{T}\times G_m$ such that
\[(\De',\vec{c}\,{}',\vec{\epsilon}\,{}' ) = \sigma \cdot \strep,\]
and then it is straightforward to check that $\sigma(B)\subset\sigma(\De)$
is an admissibly packed triangle at the vertex $\sigma(p)\in\sigma(\De)$.
Thus, by using the action of the group $\mathcal{T}\times G_m$,
an admissibly packed triangle $B$ in a single representative $\strep$ induces an admissibly packed triangle in all representatives of the semitoric polygon $\stpoly$.
So, an \emph{admissibly packed triangle in a semitoric polygon $\stpoly$} is a set
\[
 [B] = \{\sigma(B)\colon \sigma\in\mathcal{T}\times G_m\},
\]
where $B\subset \De$ is an admissibly packed triangle in $\strep$.
As with semitoric polygons, this entire infinite family is determined by any single representative, and thus we will work with a single representative of an admissibly packed triangle, passing to equivalent representatives when necessary.

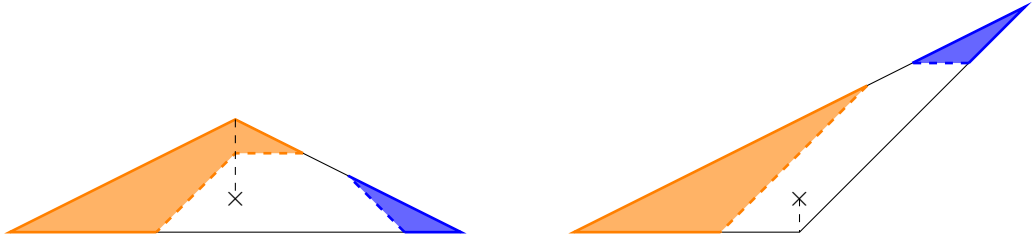
\begin{figure}[t]
 \centering
 \begin{tikzpicture}[scale = 1.5]
 \draw (0,0)--(2,1)--(4,0)--cycle;
 \fill[fill=fill1] (1.30,0)--(0,0)--(2,1)--(2.6,0.70)--(2,0.70);
 \draw[line1, line width = 1] (1.30,0)--(0,0)--(2,1)--(2.6,0.70);
 \draw [\triangleline, line1,line width = 1] (2.6,0.70) -- (2,0.70)--(1.30,0);
 \draw[line2, line width = 1, fill=fill2] (3,0.5)--(4,0)--(3.5,0);
 \draw [\triangleline, line2,line width = 1] (3.5,0) -- (3,0.5);
 \node at (2,0.3) {$\times$};
 \draw [dashed](2,1)--(2,0.3);
 \draw (5,0)--(9,2)--(7,0)--cycle;
 \fill[fill=fill1] (6.30,0)--(5,0)--(7,1)--(7.6,1.30)--(7,0.70);
 \draw[line1, line width = 1] (6.30,0)--(5,0)--(7,1)--(7.6,1.30);
 \draw [\triangleline, line1,line width = 1] (7.6,1.30) -- (7,0.70)--(6.30,0);
 \draw[line2, line width = 1, fill=fill2] (8,1.5)--(9,2)--(8.5,1.5);
 \draw [\triangleline, line2,line width = 1] (8.5,1.5) -- (8,1.5);
 \node at (7,0.3) {$\times$};
 \draw [dashed](7,0.3)--(7,0);
 \end{tikzpicture}
 \caption{Two representatives of the same packing of a semitoric polygon, related by a~change of cut direction operation. Note that the orange ``triangle" in the representative on the left becomes an actual triangle in the equivalent polygon shown on the right.}
 \label{fig:sem_weird_triangle}
 \end{figure}

An admissibly packed triangle in a semitoric polygon does not always look like a triangle in every representative, as shown in Figure~\ref{fig:sem_weird_triangle}.

Now that we have a notion of an admissibly packed triangle in a semitoric polygon, the definitions of packings and the packing density are essentially the same as in the toric case.
The main difference is that it is more precise to view a packing of a semitoric polygon as a packing of all representatives of the polygon.

\begin{Definition}\label{def:st-packing}
A \emph{packing} $P$ of a semitoric polygon representative $\strep$ is the union of a~set of pairwise disjoint admissibly packed triangles, i.e.,
\[P = \bigcup_{i=1}^d B_i,\]
where each $B_i$ is an admissibly packed triangle in $\strep$ and $B_i \cap B_j = \varnothing$
if $i \neq j$.
A packing of a semitoric polygon $\stpoly$ is
\[
 [P] = \{\sigma(P)\colon P \text{ is a packing of }\strep\text{ and }\sigma\in\mathcal{T}\times G_m\},
\]
in which case $\sigma(P)$ is a packing of $\sigma\cdot\strep$
for each $\sigma\in\mathcal{T}\times G_m$.
Let \[\PST(\stpoly) = \{[P] \mid [P] \text{ is a packing of } \stpoly \}.\]
\end{Definition}

The elements of $\mathcal{T}\times G_m$ are all piecewise affine transformations whose linear parts have determinant 1.
Thus, they preserve area.
We conclude that all representatives of a single semitoric polygon have the same area, and all representatives of a given admissibly packed triangle have the same area.
Thus, the following notion is independent of the choice of representative of the semitoric polygon and of the packing.

\begin{Definition}\label{def:semitoric-density}
The \emph{packing density} of a semitoric polygon $\stpoly$ is denoted by $\rho_{{\rm ST}}(\stpoly)$ and given by
\[
 \rho_{{\rm ST}}(\stpoly): = \sup_{[P] \in \PST(\stpoly)}\left( \frac{\area(P)}{\area(\De)}\right).
\]
\end{Definition}

Again, it is sufficient to consider packings of a single representative of $\stpoly$ to compute $\rho_{{\rm ST}}(\stpoly)$, and this is how we will proceed for the remainder of the paper.

\subsection{Edges of semitoric polygons}
\label{sec:semitoric-edge}
Recall that when changing the cut direction to move a cut away from a fake corner in a semitoric polygon, the fake corner ``disappears'' in the new representative and two edges of the original polygon merge into one edge, see for instance Figure~\ref{fig:sem_weird_triangle}.
The two representatives of the same semitoric polygon shown in Figure~\ref{fig:stmin-3b} even have a different number of edges.
This motivates the idea that the analogue of an ``edge'' in a semitoric polygon should be able to continue through a fake corner, but not through a hidden or Delzant corner.

This is particularly relevant for the present paper, since an admissibly packed triangle in a~semitoric polygon is allowed to pass through a fake corner, but not a hidden or Delzant corner, as described in Definition~\ref{def:st-packing}.
Thus, we will see that the bound on the $\SL$-sizes of adjacent triangles does not depend on the length of the actual edge between them, but rather on the length of any string of adjacent edges connected by fake corners.

\begin{figure}[t]
 \centering
 \begin{subfigure}{.45\linewidth}
 \centering
 \begin{tikzpicture}[scale = 1.0]
 \draw (0,0)--(2,2)--(5,2)--(3,0)--cycle;
 \node at (2,1) {$\times$};
 \draw [dashed](2,1)--(2,2);
 \draw[red, line width = 1.75] (0,0)--(2,2)--(5,2);
 \draw[darkgreen, line width = 1.75] (5,2)--(3,0);
 \draw[blue, line width = 1.75] (3,0)--(0,0);
 \end{tikzpicture}
 \caption{}
 \label{fig:st-edges-a}
 \end{subfigure}
 \begin{subfigure}{.45\linewidth}
 \centering
 \begin{tikzpicture}[scale = 1.0]
 \draw (0,2)--(5,2)--(3,0)--(2,0)--cycle;
 \node at (2,1) {$\times$};
 \draw [dashed](2,1)--(2,0);
 \draw[red, line width = 1.75] (0,2)--(5,2);
 \draw[darkgreen, line width = 1.75] (5,2)--(3,0);
 \draw[blue, line width = 1.75] (3,0)--(2,0)--(0,2);
 \end{tikzpicture}
 \caption{}
 \label{fig:st-edges-b}
 \end{subfigure}
 \caption{Two representatives of the same semitoric polygon with the semitoric edges shown in three different colors. Notice that the semitoric edges start and end only on Delzant or hidden corners (although there are no hidden corners in this example) but continue through fake corners. Also notice that the semitoric edges are preserved by the group action, while the usual edges of the polygon can be broken by it.}
 \label{fig:st-edges}
 \end{figure}
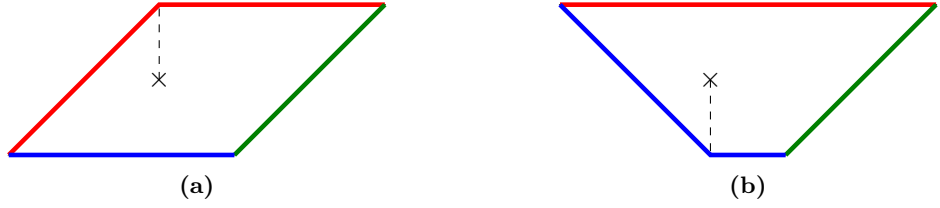

\begin{Definition}\label{def:st-edge}
Suppose that $\strep$ is a semitoric polygon representative.
A \emph{semitoric edge} is a union of adjacent edges $\tilde{e} = e_1\cup \dots \cup e_k$ such that
\begin{itemize}\itemsep=0pt
 \item one of the endpoints of $e_1$ is a hidden or Delzant corner;
 \item one of the endpoints of $e_k$ is a hidden or Delzant corner;
 \item $e_j$ is connected to $e_{j+1}$ by a fake corner for each $j=1,\ldots, k-1$.
\end{itemize}
\end{Definition}
If $k=1$, then this is simply an edge of the polygon connecting two adjacent hidden or Delzant corners.
Figure~\ref{fig:st-edges} shows two representatives of a semitoric polygon with each semitoric edge drawn in a different color.

Notice that elements of $\mathcal{T}\times G_m$ send hidden or Delzant corners to hidden or Delzant corners, and thus they send semitoric edges to semitoric edges.
We define the $\SL$-length of the semitoric edge to be the sum of the $\SL$-lengths of the edges it is the union of. That is, the $\SL$-length of the semitoric edge $\tilde{e} = e_1\cup \cdots \cup e_k$ is given by $\sum_{j=1}^k \ell_{e_j}$, where $\ell_{e_j}$ denotes the $\SL$-length of $e_j$.
Note that the action of $\mathcal{T}\times G_m$ also preserves $\SL$-length, so an element of $\mathcal{T}\times G_m$ always sends a semitoric edge to another semitoric edge of the same $\SL$-length.

Let $\stpoly$ be a semitoric polygon. We label the hidden and Delzant corners of a representative $\strep$ by $p_1,\ldots,p_d$ in clockwise order starting with the lexicographically minimal element (i.e.,~the lower left corner), and let $\tilde{e}_i$ denote the semitoric edge connecting $p_i$ to $p_{i+1}$, as usual taking $p_{d+1} = p_1$. For $i=1,\dots,d$, let $\ell_i$ denote the $\SL$-length of the semitoric edge $\tilde{e}_i$. Then, based on the discussion in the preceding paragraph, the tuple $(\ell_1,\ldots,\ell_d)$ is independent of the choice of representative.

\subsection{Computing the semitoric packing density}
\label{sec:semitoric-compute}

Let $\strep$ be a semitoric polygon representative. As in the previous subsection we assume that the hidden and Delzant corners are labeled in clockwise order by $p_1,\ldots,p_d$ and that the semitoric edges are labeled by $\tilde{e}_1,\ldots,\tilde{e}_d$ with $\SL$-lengths given by $\ell_1,\ldots,\ell_d$.
Note that the following proof is an instance in which it is helpful to change the directions of the cuts.

\begin{Lemma}\label{lem:st-firstbound}
 Let $P = \bigcup_{i=1}^d B_i$ be a packing of $\strep$ such that for each $i=1,\ldots, d$ either~${B_i=\varnothing}$ or $B_i$ is a triangle packed at $p_i$. Let $\lambda_i$ denote the $\SL$-size of $B_i$ {\rm (}taking $\lambda_i=0$ if~$B_i$ is empty$)$. Then for all $i\in\{1,\ldots,d\}$,
 \[
 \lambda_i + \lambda_{i+1}\leq \ell_i.
 \]
\end{Lemma}

\begin{proof}
Fix some $i\in\{1,\ldots,d\}$. Then there exists some $\sigma\in\mathcal{T}\times G_m$ such that the equivalent polygon representative $\sigma\cdot\strep$ has no cuts which intersect the semitoric edge $\sigma(\tilde{e}_i)$.
Since~$\sigma(B_i)\cap \sigma(B_{i+1}) = \varnothing$, the same argument as given in the proof of Proposition~\ref{prop:toric} in Section~\ref{sec:toric-computing} now shows that the sum of the $\SL$-size of $\sigma(B_i)$ and the $\SL$-size of $\sigma(B_{i+1})$ is bounded above by the $\SL$-length of $\sigma(\tilde{e}_i)$. Since the action of $\mathcal{T}\times G_m$ preserves $\SL$-length and $\SL$-size, this implies that $\lambda_i+\lambda_{i+1}\leq \ell_i$.
\end{proof}

Unlike the case for Delzant polygons, the bound given in Lemma~\ref{lem:st-firstbound} is not the only restriction on the $\SL$-sizes of triangles packed in a semitoric polygon. The requirement that the packed triangles also avoid all marked points introduces another bound, which we compute now.

\begin{figure}[t]
 \centering
 \begin{tikzpicture}[scale = 1.0]
 \draw (0,0)--(4,2)--(8,0)--cycle;
 \node at (4,1) {$\times$};
 \draw [dashed](4,1)--(4,2);
 \draw[red, line width = 1.75] (0,0)--(4,2)--(8,0);
 \draw[blue, line width = 1.75] (0,0)--(8,0);
 \node at (-0.6,0.4) {$p_1=(0,0)$};
 \node at (0,0) [circle,fill,inner sep=1.6pt]{};
 \node at (4.6,-0.4) {$\tilde{q}=(2,0)$};
 \node at (4,0) [circle,fill,inner sep=1.6pt]{};
 \node at (4.6,2.4) {$q_1^1=(2,1)$};
 \node at (4,2) [circle,fill,inner sep=1.6pt]{};
 \node at (8.3,0.4) {$(4,0)$};
 \end{tikzpicture}
 \caption{In this case the vertical line $L_1$ through the marked point intersects both the red and blue semitoric edges emanating from $p_1$, so there are two possible choices for $q_1^1$. Following the convention in the discussion before Lemma~\ref{lem:avoid-markedpt}, we choose the point on the top (red) edge as $q_1^1$, since the $\SL$-length from $p_1$ to the other candidate, labeled $\tilde{q}$, is $2$, while the $\SL$-length from $p_1$ to $q_1^1$ is $1$.}
 \label{fig:q-choice}
 \end{figure}
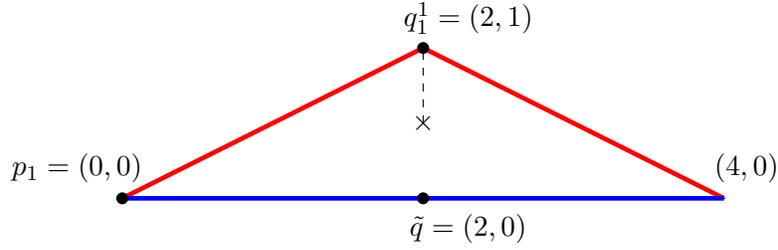

 Now we set up some notation for the following lemma.
 Let $\strep$ be a semitoric polygon representative, with hidden and Delzant corners $p_1,\ldots,p_d\in\De$.
 Let $v^i,w^i\in\Z^2$ denote the primitive vectors directing the edges emanating from the vertex $p_i$ with components
 \[v^i = \vectbig{v^i_1}{v^i_2}\qquad \text{ and }\qquad w^i = \vectbig{w^i_1}{w^i_2}.\]
 For $j\in\{1,\ldots, m\}$,
 let $L_j$ denote the vertical line through the marked point $c_j$.
 If $L_j$ intersects exactly one of the two semitoric edges emanating from $p_i$, we denote the intersection point by $q_i^j$.
 If $L_j$ intersects both of the semitoric edges emanating from $p_i$, then we let $q_i^j$ denote the intersection point for which the $\SL$-length of the piecewise linear path from $p_i$ to the intersection point along $\partial\De$ is minimal,
 if the paths have the same $\SL$-length then we may choose either point.
 See Remark~\ref{rmk:q} for further discussion on why this choice of $q_i^j$ is important, and see Figure~\ref{fig:q-choice} for an illustration of this situation.

\begin{Lemma}\label{lem:avoid-markedpt}
Let $P$ be a packing of $\strep$ and let $\lambda_1,\ldots, \lambda_d\geq 0$ denote the $\SL$-size of the triangles in the packing, where $\lambda_i$ is the $\SL$-size of the triangle packed at $p_i$.
Let~$i\in\{1,\ldots, d\}$, and $j\in\{1,\ldots, m\}$.
Suppose that $L_j$ intersects at least one of the two semitoric edges emanating from $p_i$, and
let
 \begin{itemize}\itemsep=0pt
 \item $a_i^j$ denote the $\SL$-length of the piecewise linear path in $\partial\De$ connecting $p_i$ to the intersection point $q_i^j$ which is contained in a single semitoric edge and
 \item $b_i^j$ denote the vertical distance from $c_j$ to the intersection point $q_i^j$.
 \end{itemize}
 Define
 \begin{equation}\label{eqn:alphai}
 \alpha_i = \min_{j} \left\{a_i^j + b_i^j \big|v^i_1-w^i_1\big| \,\bigg|\, \begin{matrix}L_j \text{ intersects at least one of the}\\ \text{semitoric edges emanating from $p_i$}\end{matrix}\right\},
 \end{equation}
 taking $\alpha_i = \infty$ if the set is empty and where the minimum is over all $j\in\{1,\ldots, m\}$.
 Then
 \[
 \lambda_i \leq \alpha_i
 \]
 for all $i\in \{1,\ldots,d\}$.
\end{Lemma}

\begin{proof}
It is sufficient to prove that
\begin{equation}\label{eqn:avoid-markedpt-proof}
 \lambda_i\leq a_i^j + b_i^j \big|v_1^i - w_1^i\big|
\end{equation}
for any $i\in\{1,\ldots, d\}$ and $j\in\{1,\ldots, m\}$ such that the vertical line $L_j$ through the marked point~$c_j$ intersects at least one of the semitoric edges emanating from $p_i$.
For the remainder of the proof we fix some such choice of $i$ and $j$ and we use simplified notation $p=p_i$, $c=c_j$, $L=L_j$, $a = a_i^j$, $b = b_i^j$, $v = v^i$, $w=w^i$, $q=q_i^j$, and $B=B_i$.

Notice that $a$ and $b$ measure the $\SL$-length of line segments, and thus they are invariant under the action of $\mathcal{T}\times G_m$.
Furthermore, notice that acting by an element of $\mathcal{T}\times G_m$ preserves the value $|v_1-w_1|$, which can be verified by examining equations~\eqref{eqn:T}, \eqref{eqn:frakt}, and~\eqref{eqn:group-action}.
Thus, we may freely change representatives of the polygon.
If $q$ is on the top boundary, by acting by an element of $\mathcal{T}\times G_m$ we may assume that all cuts go downwards in the representative~$\strep$.
Similarly, if $q$ is on the bottom boundary we assume that all cuts go upwards.
Thus, the portion of $\partial \De$ which connects $p$ to $q$ is a single line segment (i.e.,~it doesn't contain any other vertex of $\De$).
We assume that this segment is directed along the vector $v$, otherwise we switch $v$ and~$w$.
Furthermore, since all cuts point away from $q$ and since hidden corners intersect a cut, we conclude that the point $p$ must be a Delzant corner, so $\det(v,w) = \pm 1$.

Let $\lambda \geq 0$ denote the $\SL$-size of the triangle $B\subset \De$ which is packed at the point $p$.
If~$\lambda=0$ then equation~\eqref{eqn:avoid-markedpt-proof} is automatically satisfied, and we are done.
Since we have assumed that all cuts go away from the point $q$ on the boundary, there are no cuts which intersect $B$, and therefore
\[
 B = M(B(\lambda))+p
\]
for some $M\in\SL$ where $B(\lambda)$ is the model triangle of size $\lambda$.
By the definition of an admissibly packed triangle, $c\notin B$. Keeping the same $M$, we can define
\[
 \widetilde{B} = M\big(B\big(\widetilde{\lambda}\big)\big)+p,
\]
where
\[
 \widetilde{\lambda} = \sup\{\lambda_0>0 \mid c\notin M(B(\lambda_0))+p\}.
\]
Since $c\notin B$ we conclude that $\lambda\leq \widetilde{\lambda}$.
Note that $\widetilde{B}$ might not be an admissibly packed triangle, since for instance it may not be entirely contained within $\De$, but at the moment we are interested in the bound introduced by the presence of $c$, and in any case we still have that $\lambda\leq \widetilde{\lambda}$.
We will now complete the proof by showing that
\[
 \widetilde{\lambda} = a + b |v_1-w_1|,
\]
where $v_1$ and $w_1$ are the first components of the vectors $v$ and $w$, respectively.

Notice that $\widetilde{\lambda}$ is the value such that
$c\in\partial \widetilde{B}\setminus\widetilde{B}$.
The triangle $\widetilde{B}$ has vertices $p$, $p+\widetilde{\lambda} v$, and $p+\widetilde{\lambda}w$, and the point $c$ lies on the edge of $\widetilde{B}$ that connects $p+\widetilde{\lambda} v$ to $p+\widetilde{\lambda}w$. This is enough to determine the value of $\widetilde{\lambda}$, which we do now.
This entire scenario is shown in Figure~\ref{fig:one_cut}.

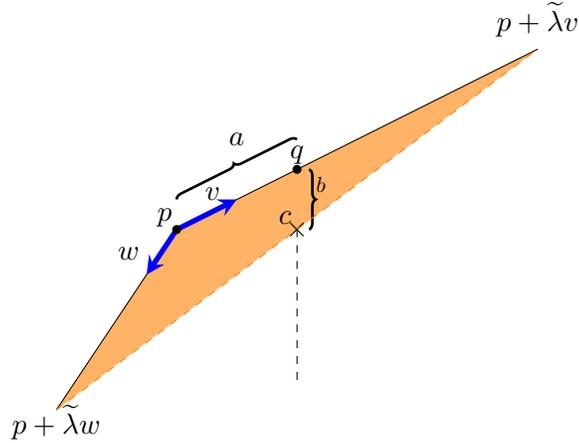
\begin{figure}[t]
\begin{center}
 \begin{tikzpicture}[scale = 0.8]
 \draw [fill=fill1](-2,-3) -- (0,0)--(6,3);
 \node at (2,0) {$\times$};
 \node at (1.8,.2) {$c$};
 \draw [dashed](2,-2.5)--(2,0);
 \node at (2,1.25) {$q$};
 \node at (2,1) [circle,fill,inner sep=1.2pt]{};
 \draw [decorate, line width=1pt, decoration={brace, mirror}] (2.2,0)--(2.2,1);
 \node[scale=0.75] at (2.4, .75) {$b$};
 \draw[-stealth, blue, line width = 2] (0,0) -- (-0.5,-0.75);
 \draw[-stealth, blue, line width = 2] (0,0) -- (1,0.5);
 \node at (0.6,0.6) {$v$};
 \node at (-0.8,-0.4) {$w$};
 \draw[\triangleline, line1] (-2,-3) -- (6,3);
 \node at (-0.2,0.2) {$p$};
 \node at (0,0) [circle,fill,inner sep=1.2pt]{};
 \draw [decorate, line width=1pt, decoration={brace}] (0,0.5) -- (2,1.5);
 \node at (1,1.5) {$a$};
 \node at (6,3.5) {$p+\widetilde{\lambda}v$};
 \node at (-2,-3.2) {$p+\widetilde{\lambda}w$};
 \end{tikzpicture}
 \caption{The situation described in Lemma~\ref{lem:avoid-markedpt}. The triangle $\widetilde{B}$ is taken to be as large as possible until it hits the marked point $c$, therefore determining a bound on the size of triangle that can be packed at $p$. Note that the braces indicating the lengths of intervals are showing the $\SL$-length of those intervals.} \label{fig:one_cut}
 \end{center}
\end{figure}

The line segment connecting $p$ with $q$ is directed along the primitive vector $v$ and has $\SL$-length $a$, so $q = p+av$.
Denote $u = \vect{0}{1}$ and let $\delta = 1$ if $q$ is on the top boundary of $\De$ and $\delta = -1$ otherwise. Then
\[
 c = q-b\delta u = p+av - b\delta u.
\]
Since $c$ lies on the line segment connecting $p+\widetilde{\lambda}w$ to $p+\widetilde{\lambda}v$, we conclude that $\widetilde{\lambda}w - \widetilde{\lambda}v$ and $p+\widetilde{\lambda}v-c$ are colinear, so $\det\big(w-v, p+\widetilde{\lambda}v-c\big)=0$.

Putting all of this together, we conclude that
\begin{align*}
 0 & = \det \big(w - v, p+\widetilde{\lambda}v-c\big) = \det \big(w-v, \big(\widetilde{\lambda} - a\big)v + \delta bu\big)\\
 & = \big(\widetilde{\lambda} - a\big)\det (w,v) + \delta b \det (w, u) - \big(\widetilde{\lambda}-a\big) \det (v, v) - \delta b \det (v, u)\\
 & = \big(\widetilde{\lambda} - a\big) \det (w,v) + \delta b(w_1 - v_1).
\end{align*}
Thus, combining this computation with the fact that $a,b>0$ and the fact that $p$ is a Delzant corner so $\det(w,v)=\pm 1$, we have that
\[
 \widetilde{\lambda} = a + \delta \det(w,v) (v_1-w_1)b.
\]
We already know that $\delta \det(w,v) = \pm 1$, and to prove the claim it is now sufficient to verify that~$\delta \det(w,v) (v_1-w_1) = |v_1-w_1|$.
We omit the details, but this follows from checking many similar cases ($q$ is on the top or bottom boundary of $\De$, $q$ is to the left or right of $p$, and whether or not $p$ is the extreme left or right vertex of the polygon). Note that in the case that the vertical line through $c$ intersects both semitoric edges emanating from $p$, which can only occur if $p$ is the extreme right or left vertex of $\De$, we must make use of the fact that $q$ was chosen to be the intersection points which is along the path of minimal $\SL$-length (see Remark~\ref{rmk:q}).

We conclude that $\widetilde{\lambda} = a + b|v_1-w_1|$, which completes the proof.
\end{proof}

\begin{Remark}\label{rmk:q}
In the case that $L_j$ intersects both of the semitoric edges emanating from $p_i$, it is important to choose the correct intersection point to be labeled by $q_i^j$, otherwise equation~\eqref{eqn:alphai} will be incorrect.
If $q_i^j$ is instead chosen to be the intersection which is along a path of longer $\SL$-length from $p_i$, then
in the last part of the proof we will obtain instead that $\delta \det(w,v) (v_1-w_1) = -|v_1-w_1|$, a change in sign.
The values of $a_i^j$ and $b_i^j$ will also be different, in a way which compensates for the change in sign, so the result would be the same but the formula would look slightly different. To avoid this confusion we simply specified which intersection point to choose.
For an example of this situation, consider
the minimal system of type~(1) whose maximal packing is computed in Section~\ref{sec:type1-packing}.
\end{Remark}

With the above results, we are now prepared to prove Theorem~\ref{thm:semitoric}.

\begin{proof}[Proof of Theorem~\ref{thm:semitoric}]
Let $\stpoly$ be a semitoric polygon, and define $\PST := \PST (\stpoly)$ and $\PSTtilde := \PSTtilde(\stpoly)$ as in equation~\eqref{eqn:PSTtilde}, using simplified notation which suppresses the dependence on $\stpoly$ for the duration of this proof.

Fix a representative $\strep$ of $\stpoly$.
Denote the hidden and Delzant corners by $p_1,\ldots, p_d$, labeled in clockwise order starting with the lexicographically minimal vertex, and as usual let $\ell_i$ denote the $\SL$-length of the semitoric edge connecting $p_i$ to $p_{i+1}$.
Then, as in equation~\eqref{eqn:LambdaST}, $\Lambda\colon \PST \to \R^d$ is given by
\[
\Lambda\left(\left[\bigcup_{i=1}^d B_i\right]\right) = (\lambda_1,\ldots,\lambda_d),
\]
where $\lambda_i$ is the $\SL$-size of $B_i$ and $B_i$ denotes the admissibly packed triangle at the point~$p_i$ in the representative $\strep$, taking $\lambda_i=0$ if $B_i=\varnothing$.
Then Lemmas~\ref{lem:st-firstbound} and~\ref{lem:avoid-markedpt} immediately imply that
\[
 \Lambda(\PST)\subset \PSTtilde.
\]
The map $\Lambda$ is clearly injective.
To show that $\Lambda$ is surjective onto $\PSTtilde$, let $(\lambda_1,\ldots,\lambda_d)\in\PSTtilde$.
Fix any $i\in\{1,\ldots,d\}$.
We can choose a representative of $\stpoly$ in which $p_i$ corresponds to a~vertex~$\widetilde{p}_i$ which is a Delzant corner with edges of $\SL$-length $\ell_{i-1}$ and $\ell_i$ emanating from it.
Since $\lambda_i\leq \ell_i$ and $\lambda_i\leq \ell_{i-1}$, a triangle of size $\lambda_i$ can be admissibly packed at the vertex~$\widetilde{p}_i$, and by the proof of Lemma~\ref{lem:avoid-markedpt} we see that since $\lambda_i\leq \alpha_i$ this triangle does not contain any marked points.
Returning to the original representative $\strep$, we conclude that for each~$i\in\{1,\ldots,d\}$ there exists an admissibly packed triangle at $p_i$ of size $\lambda_i$, call it $B_i$.
Since~$\lambda_i+\lambda_{i+1}\leq \ell_i$, we conclude that the $B_i$ are disjoint, and so $P = \bigcup_{i=1}^d B_i$ is an admissible packing of the representative~$\strep$.
Thus, $[P]\in\PST$ with the property that~$\Lambda([P]) = (\lambda_1,\ldots,\lambda_d)$, as desired. Thus $\Lambda$ is surjective onto $\PSTtilde$.

The set $\PSTtilde$ is a compact convex polygonal set because it is described by a finite number of linear inequalities and contained within the compact set $\{(\lambda_1,\ldots,\lambda_d)\in\R^d \mid 0\leq \lambda_i\leq \ell_i\}$.

For $[P]\in\PST$, the equality $\area ([P]) = \frac{1}{2}||\Lambda([P])||^2$ is clear, since the area of an admissibly packed triangle with $\SL$-size $\lambda$ is $\frac{1}{2}\lambda^2$, as in the proof of Proposition~\ref{prop:toric}.
\end{proof}

Just as with Proposition~\ref{prop:toric} in the case of Delzant polygons, Theorem~\ref{thm:semitoric} makes computing the packing density of a semitoric polygon much easier.
The following result is the semitoric analogue of Corollary~\ref{cor:toric-vertex} from the Delzant case.

\begin{Corollary} \label{cor:semitoric-vertex}
 Let $\stpoly$ be a semitoric polygon.
 Then
 \[
 \rho_{{\rm ST}} (\stpoly) = \sup \left\{ \frac{||q||^2}{2 \area(\De)}\colon q\in V\big(\PSTtilde(\stpoly)\big)\right\},
 \]
 where $V\big(\PSTtilde(\stpoly)\big)$ is the finite set of vertices of the compact convex polytope $\PTtilde(\stpoly)$, given in equation~\eqref{eqn:PSTtilde}.
 \end{Corollary}

\begin{proof}This follows immediately from Theorem~\ref{thm:semitoric} and Proposition~\ref{prop:vert} applied to the compact convex polytope $\PSTtilde(\stpoly)$.
\end{proof}

\begin{Remark}
 In Section~\ref{sec:program} we discussed an algorithm for finding the toric packing density of a given Delzant polygon.
 Making use of Theorem~\ref{thm:semitoric} and Corollary~\ref{cor:semitoric-vertex}, a similar algorithm can also be developed for computing the semitoric packing density of a given representative of a semitoric polygon.
 The only extra difficulty is computing the bound $\alpha_i$ from equation~\eqref{eqn:alphai}.
\end{Remark}

\subsection{Examples}\label{sec:pack-st-examples}
\begin{figure}[t]
 \centering
 \begin{subfigure}{.49\textwidth}
 \centering
 \begin{tikzpicture}[scale = 1]
 \draw (0,0) -- (2,1) -- (4,0) -- cycle;
 \draw[line1, fill=fill1] (.5,0)--(0,0)--(1,.5);
 \draw[\triangleline, line1] (1,.5)--(.5,0);
 \fill[fill=fill2] (4,0)--(2,1)--(1,0.5)--(2,0.5)--(2.5,0)--cycle;
 \draw[line2] (2.5,0)--(4,0)--(2,1)--(1,0.5);
 \draw[\triangleline,line2] (1,0.5)--(2,0.5)--(2.5,0);
 \node at (2,0.5) {$\times$};
 \draw [dashed](2,1)--(2,0.5);
 \node at (0,0) [anchor = east] {$(0,0)$};
 \node at (2,1) [anchor = south] {$\left(a,\frac{a}{2}\right)$};
 \node at (4,0) [anchor = west] {$(2a,0)$};
 \end{tikzpicture}
 \caption{Maximal packing of type (1)}
 \label{fig:stmaxpacking-1}
 \end{subfigure}\qquad
 \begin{subfigure}{.46\textwidth}
 \centering
 \begin{tikzpicture}[scale = .7]
 \draw (0,0) -- (2,2) -- (4,2) -- (6,0) -- cycle;
 \draw[line1, fill=fill1] (2,2)--(0,0) -- (2,0);
 \draw[\triangleline, line1] (2,2)--(2,0);
 \draw[line2, fill=fill2] (4,2) -- (6,0) -- (4,0);
 \draw[\triangleline,line2] (4,2)--(4,0);
 \node at (2,1) {$\times$};
 \draw [dashed](2,1)--(2,2);
 \node at (4,.5) {$\times$};
 \draw [dashed](4,.5)--(4,2);
 \node at (0,0) [anchor = east] {$(0,0)$};
 \node at (2,2) [anchor = east] {$(a,a)$};
 \node at (4,2) [anchor = west] {$(a+b,a)$};
 \node at (6,0) [anchor = west] {$(2a+b,0)$};
 \end{tikzpicture}
 \caption{Maximal packing of type (2)}
 \label{fig:stmaxpacking-2}
 \end{subfigure}\\[1.5em]
 \begin{subfigure}{.49\textwidth}
 \centering
 \begin{tikzpicture}[scale = .9]
 \draw (0,0) -- (2,2) -- (3,2) -- (7,0) -- cycle;
 \draw[line1, fill=fill1] (2,2)--(0,0) -- (2,0);
 \draw[\triangleline,line1] (2,0)--(2,2);
 \draw[line2, fill=fill2] (3,2) -- (7,0) -- (5,0);
 \draw[\triangleline,line2] (5,0)--(3,2);
 \node at (2,1) {$\times$};
 \draw [dashed](2,1)--(2,2);
 \node at (0,0) [anchor = north] {$(0,0)$};
 \node at (2,2) [anchor = east] {$(a,a)$};
 \node at (3.4,2) [anchor = west] {$(a+b,a)$};
 \node at (7,0) [anchor = north] {$(na+b,0)$};
 \end{tikzpicture}
 \caption{Maximal packing of type (3a)}
 \label{fig:stmaxpacking-3a}
 \end{subfigure}\qquad
 \begin{subfigure}{.46\textwidth}
 \centering
 \begin{tikzpicture}[scale = .7]
 \draw (0,0) -- (2,2) -- (6,0) -- cycle;
 \draw[line1, fill=fill1] (2,2) -- (0,0) -- (2,0);
 \draw[\triangleline,line1] (2,0)--(2,2);
 \draw[line2, fill=fill2] (2,2) -- (6,0) -- (4,0);
 \draw[\triangleline,line2] (4,0)--(2,2);
 \node at (2,1) {$\times$};
 \draw [dashed](2,1)--(2,2);
 \node at (0,0) [anchor = north] {$(0,0)$};
 \node at (2,2) [anchor = south] {$(a,a)$};
 \node at (6,0) [anchor = north] {$(na,0)$};
 \end{tikzpicture}
 \caption{Maximal packing of type (3b)}
 \label{fig:stmaxpacking-3b}
 \end{subfigure}\\[1.5em]
 \begin{subfigure}{.9\textwidth}
 \centering
 \begin{tikzpicture}[scale = .5]
 \draw (0,0) -- (2,2) -- (4,3) -- (10,0) -- cycle;
 \draw[line1, fill=fill1] (2,2) -- (0,0) -- (2,0);
 \draw[\triangleline, line1] (2,0) -- (2,2);
 \fill[fill=fill2] (2.4,2.2) -- (4,3) -- (10,0) -- (6.2,0) -- (4,2.2) -- cycle;
 \draw[line2] (2.4,2.2) -- (4,3) -- (10,0) -- (6.2,0);
 \draw[\triangleline, line2] (2.4,2.2) -- (4,2.2) -- (6.2,0);
 \node at (4,2.2) {$\times$};
 \draw [dashed](4,3)--(4,2.2);
 \node at (0,0) [anchor = north] {$(0,0)$};
 \node at (2,2) [anchor = east] {$(a+b,a+b)$};
 \node at (4.3,3.2) [anchor = west] {$(a,a+\frac{b}{n-1})$};
 \node at (10,0) [anchor = north] {$(na+b,0)$};
 \begin{scope}[shift = {(14,0)}]
 \draw (0,0) -- (2,2) -- (4,3) -- (10,0) -- cycle;
 \draw[line1, fill=fill1] (2,2) -- (0,0) -- (2,0);
 \draw[\triangleline, line1] (2,0) -- (2,2);
 \fill[fill=fill2] (2,2) -- (4,3) -- (10,0) -- (6,0) -- (4,2) -- cycle;
 \draw[line2] (2,2) -- (4,3) -- (10,0) -- (6,0);
 \draw[\triangleline, line2] (6,0) -- (4,2) -- (2,2);
 \node at (4,0.5) {$\times$};
 \draw [dashed](4,3)--(4,0.5);
 \node at (0,0) [anchor = north] {$(0,0)$};
 \node at (2,2) [anchor = east] {$(a+b,a+b)$};
 \node at (4.3,3.2) [anchor = west] {$(a,a+\frac{b}{n-1})$};
 \node at (10,0) [anchor = north] {$(na+b,0)$};
 \end{scope}
 \end{tikzpicture}
 \caption{Maximal packings of type (3c) for the cases of $\alpha_3<a$ and $\alpha_3\geq a$}
 \label{fig:stmaxpacking-3c}
 \end{subfigure}
 \caption{One representative of the maximal packing of each of the minimal semitoric polygons of types (1), (2), (3a), and (3b).
 In Figure~\ref{fig:stmaxpacking-3c} we show the two distinct cases of maximal packings of the minimal semitoric polygon of type (3c) dependent on the parameters. The representative shown on the left is when $\alpha_3 < a$ and the one on the right is when $\alpha_3 \geq a$.}
 \label{fig:semitoric-packings}
\end{figure}
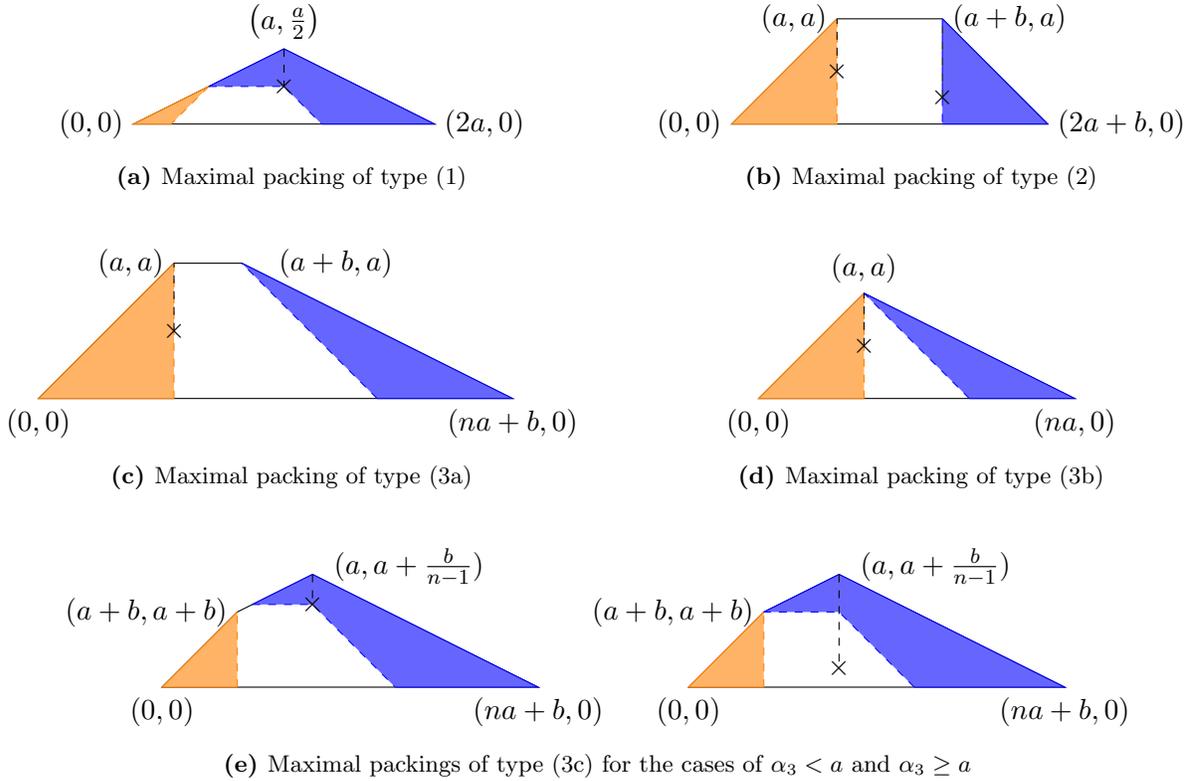

In this section, we make use of Theorem~\ref{thm:semitoric} to compute the packing density of the examples described in Definition~\ref{def:st-examples}.
Figure~\ref{fig:semitoric-packings} shows a maximal packing of one representative of each of these semitoric polygons.

\subsubsection{Minimal of type (1)}
\label{sec:type1-packing}
Let $[\De_{(1)}] = [\De,((a,h)),(+1)]$ denote the minimal semitoric polygon of type (1) with parameter~$a>0$ as defined in Definition~\ref{def:st-examples}.
To calculate the inequalities describing \smash{$\PSTtilde([\De_{(1)}])$} we must compute $\ell_1$, $\ell_2$, $\alpha_1$, and $\alpha_2$.
The main difficulty is finding the values for $\alpha_1$ and $\alpha_2$. For instance, let us compute $\alpha_1$.
This is exactly the difficult situation described in Lemma~\ref{lem:avoid-markedpt} and Remark~\ref{rmk:q} in which the vertical line through the marked point passes through both semitoric edges emanating from $p_1$, and therefore we must compare the $\SL$-lengths of the segments from $p_1$ to the vertical line to choose the point $q_1^1$.
The line segment from $p_1=(0,0)$ to $(a,0)$ has $\SL$-length $a$ and the line segment from $p_1$ to $\big(a,\frac{a}{2}\big)$ has $\SL$-length $\frac{a}{2}$, and thus following the discussion before Lemma~\ref{lem:avoid-markedpt} we define $q_1^1 = \big(a,\frac{a}{2}\big)$ and using equation~\eqref{eqn:alphai} obtain~$\alpha_1 = a-h$.
The computation for $\alpha_2$ is similar.

We have $\textrm{area}(\De) = \frac{a^2}{2}$, $\ell_1 = a$, $\ell_2 = 2a$, $\alpha_1 = \alpha_2 = a-h$. Thus, the inequalities describing~$\PSTtilde([\De_{(1)}])$ are $\lambda_1,\lambda_2\geq 0$ and
\begin{equation*}
\lambda_1+\lambda_2 \leq a,\qquad
 \lambda_2+\lambda_1 \leq 2a,\qquad
 \lambda_1 \leq a-h,\qquad
 \lambda_2 \leq a-h.
\end{equation*}
We find the maximal magnitude of $(\lambda_1,\lambda_2)\in\PSTtilde([\De_{(1)}])$ is achieved at $(\lambda_1,\lambda_2) = (a-h,h)$ and $(\lambda_1,\lambda_2) = (h,a-h)$.
Applying Theorem~\ref{thm:semitoric}, we obtain that the packing density is
\begin{equation}\label{eqn:stpack-1}
\rho_{{\rm ST}} \left([\De_{(1)}]\right) = \frac{a^2+2h^2-2ah}{a^2}.
\end{equation}
We observe that this packing density is always strictly greater than $\frac{1}{2}$ and strictly less than $1$, as $0< h < \frac{a}{2}$.

\subsubsection{Minimal of type (2)} Let
$[\De_{(2)}] := [\De,((a,h_1),(a+b,h_2)),(+1,+1)]$
denote the minimal semitoric polygon of type~(2) with parameters $a>0$ and $b\geq 0$ as defined in Definition~\ref{def:st-examples}. Using similar arguments as in Section~\ref{sec:type1-packing}, we have $\textrm{area}(\De) = a(a+b)$, $\ell_1 = 2a+b$, $\ell_2 = 2a+b$, $\alpha_1 = \alpha_2 = a$. Thus, the inequalities describing $\PSTtilde([\De_{(2)}])$ are $\lambda_1,\lambda_2\geq 0$ and
\begin{equation*}
\lambda_1+\lambda_2 \leq 2a+b,\qquad
 \lambda_1 \leq a,\qquad
 \lambda_2 \leq a.
\end{equation*}
We find the maximal magnitude of $(\lambda_1,\lambda_2)\in\PSTtilde([\De_{(2)}])$ is achieved at $(\lambda_1,\lambda_2) = (a,a)$.
Applying Theorem~\ref{thm:semitoric}, we obtain that the packing density is
\begin{equation}\label{eqn:stpack-2}
\rho_{{\rm ST}} ([\De_{(2)}]) = \frac{a}{a+b}.
\end{equation}

Notice that the packing density is independent of $h_1$ and $h_2$ and we obtain a perfect packing when $b=0$ (see Section~\ref{sec:perfect}).

\subsubsection{Minimal of type (3a)}
For the case of the minimal polygons of types (3a), (3b), and (3c), it is helpful to apply Corollary~\ref{cor:semitoric-vertex} so that only the vertices of $\PSTtilde(\stpoly)$ need to be checked to find the semitoric packing.
We suppress the details here, but the computation is very similar to the one for computing the toric packing density of the Hirzebruch trapezoid in Section~\ref{sec:hirz-packing}.

Let $[\De_{(3a)}]:=[\De,((a,h),(+1)]$ denote the minimal semitoric polygon of type (3a) with parameters $a>0$, $b>0$, and $n\in\Z_{\geq 1}$ as defined in Definition~\ref{def:st-examples}. We have \smash{$\textrm{area}(\De) = \frac{(na+2b)a}{2}$}, $\ell_1 = a+b$, $\ell_2 = a$, $\ell_3 = na+b$, $\alpha_1 = a$, $\alpha_2 = b+n(a-h)$. Thus, the inequalities describing $\PSTtilde([\De_{(3a)}])$ are $\lambda_1,\lambda_2,\lambda_3\geq 0$ and
\begin{gather*}
\lambda_1+\lambda_2 \leq a+b,\qquad
 \lambda_2+\lambda_3 \leq a,\qquad
 \lambda_3+\lambda_1 \leq na+b,\\
 \lambda_1 \leq a,\qquad
 \lambda_2 \leq b+n(a-h).
\end{gather*}
First suppose that $n>1$. Then we find the maximal magnitude of $(\lambda_1,\lambda_2,\lambda_3)\in\PSTtilde([\De_{(3a)}])$ is achieved at $(\lambda_1,\lambda_2,\lambda_3) = (a,0,a)$. Applying Theorem~\ref{thm:semitoric}, we obtain that the packing density is\looseness=-1
\begin{equation}\label{eqn:stpack-3a}
\rho_{{\rm ST}} ([\De_{(3a)}]) = \frac{2a}{na+2b}.
\end{equation}
As in the case for the toric packing density of the Hirzebruch trapezoid, if $n=1$ there are various cases of the semitoric packing density of type (3a) depending on the relative values of $a$ and $b$.
Making use of Corollary~\ref{cor:semitoric-vertex}, one can readily check that if $a\leq b$ the maximum packing is given by $(a,a,0)$, if $b<a<2b$ the maximal packing is given by $(a, b, a-b)$, and if $a\geq 2b$ the maximal packing is given by $(a,b,b)$. Computing the densities of these packings produces the results listed in the statement of Theorem~\ref{thm:semitoric-examples}.

\subsubsection{Minimal of type (3b)}
Let $[\De_{(3b)}]:=[\De,((a,h),(+1)]$ denote the minimal semitoric polygon of type (3b) with parameters $a>0$ and $n\in\Z_{\geq 2}$ as defined in Definition~\ref{def:st-examples}. We have \smash{$\textrm{area}(\De) = \frac{na^2}{2}$}, $\ell_1 = a$, $\ell_2 = a$, $\ell_3 = na$, $\alpha_1 = a$, $\alpha_3 = a$. Thus, the inequalities describing \smash{$\PSTtilde([\De_{(3b)}])$} are $\lambda_1,\lambda_2,\lambda_3\geq 0$ and
\begin{equation*}
\lambda_1+\lambda_2 \leq a,\qquad
 \lambda_2+\lambda_3 \leq a,\qquad
 \lambda_3+\lambda_1 \leq na,\qquad
 \lambda_1 \leq a,\qquad
 \lambda_3 \leq a.
\end{equation*}
We find that maximal magnitude of $(\lambda_1,\lambda_2,\lambda_3)\in\PSTtilde([\De_{(3b)}])$ is achieved at $(\lambda_1,\lambda_2,\lambda_3) = (a,0,a)$. Applying Theorem~\ref{thm:semitoric} we obtain that the packing density is
\begin{equation}\label{eqn:stpack-3b}
\rho_{{\rm ST}} ([\De_{(3b)}]) = \frac{2}{n}.
\end{equation}
Notice that the packing density is also independent of $h$ and we obtain a perfect packing if $n=2$ (see Section~\ref{sec:perfect}).

\subsubsection{Minimal of type (3c)}
Let $[\De_{(3c)}]:=[\De,((a,h),(+1)]$ denote the minimal semitoric polygon of type (3c) with parameters $a>0$, $b<0$, and $n\in\Z_{\geq 2}$ as defined in Definition~\ref{def:st-examples}. We have \smash{$\textrm{area}(\De) = \frac{(na+2b)a}{2}$}, $\ell_1 = a+b$, $\ell_2 = a$, $\ell_3 = na+b$, $\alpha_2 = n(a-h)+b$, $\alpha_3 = (n-1)a+b - (n-2)h$. Thus, the inequalities describing \smash{$\PSTtilde([\De_{(3c)}])$} are $\lambda_1,\lambda_2,\lambda_3\geq 0$ and
\begin{gather*}
 \lambda_1+\lambda_2 \leq a+b,\qquad
 \lambda_2+\lambda_3 \leq a,\qquad
 \lambda_3+\lambda_1 \leq na+b,\qquad
 \lambda_2 \leq n(a-h)+b,\\
 \lambda_3 \leq (n-1)a+b - (n-2)h.
\end{gather*}
We find the maximal magnitude of $(\lambda_1,\lambda_2,\lambda_3)\in\PSTtilde([\De_{(3c)}])$ is achieved at $(\lambda_1,\lambda_2,\lambda_3) = (a+b,0,\min\{\alpha_3,a\})$. Applying Theorem~\ref{thm:semitoric}, we obtain that the packing density is
\begin{equation}\label{eqn:stpack-3c}
\rho_{{\rm ST}} ([\De_{(3c)}]) = \frac{(a+b)^2+\min\bigl\{ ((n-1)a+b - (n-2)h )^2,a^2\bigr\}}{(na+2b)a}.
\end{equation}
Notice that this expression is independent of $h$ in the case that $n=2$.

Theorem~\ref{thm:semitoric-examples} now
follows from equations~\eqref{eqn:stpack-1}, \eqref{eqn:stpack-2}, \eqref{eqn:stpack-3a}, \eqref{eqn:stpack-3b}, and~\eqref{eqn:stpack-3c}.
In all of these cases we made use of Theorem~\ref{thm:semitoric} to compute the packing density.

\section{Perfect packings of semitoric polygons}
\label{sec:perfect}

We call a toric or semitoric packing \emph{perfect} if it has density 1, which means that the packed triangles cover all of the polygon except for a set of measure zero.
In~\cite{Pe2006}, Pelayo showed that the only Delzant polygons which admit perfect toric packings are the square and the Delzant triangle.

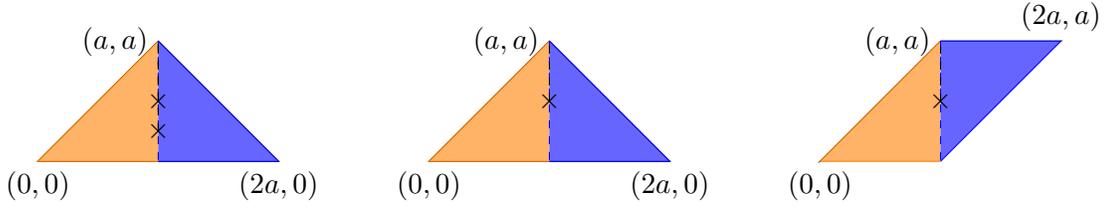
\begin{figure}[t]
 \centering
 \begin{tikzpicture}[scale = 0.8]
 \draw (0,0) -- (2,2) -- (4,0) -- cycle;
 \draw[line1, fill=fill1] (2,2)--(0,0) -- (2,0);
 \draw[\triangleline, line1] (2,2)--(2,0);
 \draw[line2, fill=fill2] (2,2) -- (4,0) -- (2,0);
 \draw[\triangleline,line2] (2,2)--(2,0);
 \node at (2,1) {$\times$};
 \draw [dashed](2,1)--(2,2);
 \node at (2,.5) {$\times$};
 \draw [dashed](2,.5)--(2,2);
 \node at (0,0) [anchor = north] {$(0,0)$};
 \node at (2,2) [anchor = east] {$(a,a)$};
 \node at (4,0) [anchor = north] {$(2a,0)$};
 \begin{scope}[shift={(6.5,0)}]
 \draw (0,0) -- (2,2) -- (4,0) -- cycle;
 \draw[line1, fill=fill1] (2,2)--(0,0) -- (2,0);
 \draw[\triangleline, line1] (2,2)--(2,0);
 \draw[line2, fill=fill2] (2,2) -- (4,0) -- (2,0);
 \draw[\triangleline,line2] (2,2)--(2,0);
 \node at (2,1) {$\times$};
 \draw [dashed](2,1)--(2,2);
 \node at (0,0) [anchor = north] {$(0,0)$};
 \node at (2,2) [anchor = east] {$(a,a)$};
 \node at (4,0) [anchor = north] {$(2a,0)$};
 \end{scope}
 \begin{scope}[shift={(13,0)}]
 \draw (0,0) -- (2,2) -- (4,2) -- cycle;
 \draw[line1, fill=fill1] (2,2)--(0,0) -- (2,0);
 \draw[\triangleline, line1] (2,2)--(2,0);
 \draw[line2, fill=fill2] (2,2) -- (4,2) -- (2,0);
 \draw[\triangleline,line2] (2,2)--(2,0);
 \node at (2,1) {$\times$};
 \draw [dashed](2,1)--(2,2);
 \node at (0,0) [anchor = north] {$(0,0)$};
 \node at (2,2) [anchor = east] {$(a,a)$};
 \node at (4,2) [anchor = south] {$(2a,a)$};
 \end{scope}
 \end{tikzpicture}
 \caption{Each of the three semitoric polygons which admits a perfect packing. Note that all cuts are going up; there are two cuts in the first figure and one in each of the others.
 From left to right, the polygons are: type (2), type (3b), and the inverted type (3b).}
 \label{fig:perfect-packings}
\end{figure}

From Theorem~\ref{thm:semitoric-examples}, we can already see two examples of semitoric polygons that admit a~perfect packing: the minimal polygon of type (2) with any parameter $a>0$ and taking $b=0$, and the minimal polygon of type (3b) with parameters $n=2$ and any $a>0$.
In this section, we will prove Theorem~\ref{thm:perfect}, which states that these two examples and one other (defined below) are the only semitoric polygons which admit a perfect packing, all of which are shown in Figure~\ref{fig:perfect-packings}.

\begin{figure}[t]
 \centering
 \begin{tikzpicture}[scale = .7]
 \draw (0,0) -- (2,2) -- (6,2) -- (2, 0) -- cycle;
 \node at (2,1) {$\times$};
 \draw [dashed](2,1)--(2,2);
 \node at (0,0) [anchor = north] {$(0,0)$};
 \node at (2,2) [anchor = south] {$(a,a)$};
 \node at (6,2) [anchor = south] {$(na,a)$};
 \node at (2,0) [anchor = north] {$(a,0)$};
 \end{tikzpicture}
 \caption{A representative of the inverted (3b) minimal semitoric polygon.}
 \label{fig:inverted3b}
\end{figure}
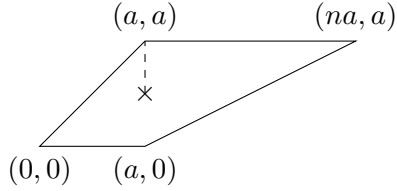

\begin{Definition}\label{def:inverted-3b}
 A semitoric polygon $[\De,((a,h)),(+1)]$ is called \emph{minimal of type inverted $(3b)$} if $\De$ is the polygon with vertices $(0,0)$, $(a,a)$, $(an, a)$ and $(a,0)$,
 where $a>0$, $n\in\Z$, $n\geq 2$, and $0<h<a$ (see Figure~\ref{fig:inverted3b}).
\end{Definition}

We start with a lemma about what a perfect packing has to look like nearby a marked point.
Recall that the model triangle $B(\lambda)$ contains two of its three edges and, in an admissibly packed triangle, it is those two edges that are required to be set to the boundary of the semitoric polygon. Thus, given some admissibly packed $B$ in a semitoric polygon, the set $e = \overline{B}\setminus B$ is the portion of the boundary of $B$ which is contained in the interior of the semitoric polygon.

\begin{Lemma}\label{lem:perfect-local}
 Let $\strep$ be a semitoric polygon representative.
 Suppose that $P\subset \De$ is a~perfect semitoric packing of $\De$, and let $c\in\De$ denote one of the marked points.
 Let $P = \bigcup_{i=1}^d B_i$, where each $B_i$ is an admissibly packed triangle in $\De$.
 Let $e_i = \overline{B}_i\setminus B_i$ for each $i=1,\ldots, d$.
 Then~$c\in e_i$ for exactly two values of $i\in\{1,\ldots,d\}$, and
for each such $i$ the set $e_i$ is a vertical line segment.
\end{Lemma}

\begin{proof}
 Recall that there are no marked points contained in an admissibly packed triangle, so since $\De = \overline{P}$ we conclude that $c\in e_i$ for at least one value of $i$, without loss of generality, we assume that $c\in e_1$.
 By Definition~\ref{def:STadm-packed-triangle}, we conclude that by changing to an equivalent representative of the semitoric polygon, we may arrange that $B_1$ does not intersect any of the cuts, and thus $e_1$ is a straight line segment.

 If $c\notin e_i$ for all $i\geq 2$, then the packing will not be perfect, since any open neighborhood of $c$ will contain a set of positive measure not contained in $P$. Furthermore, since the $B_i's$ are disjoint, it is impossible for $c$ to be contained in $e_i$ for three or more distinct $i$. We conclude that~$c\in e_i$ for exactly two values of $i$, and without loss of generality, we denote these by $e_1$ and~$e_2$. Since the packing is perfect, we must have that $e_1=e_2$ as sets in $\R^2$.

 Now, suppose that $e_1$ is not vertical, and thus since the cut emanating from $c$ does not intersect $B_1$ it must intersect $B_2$. Then $e_2$ is not a straight line segment, since a cut passes through it, but this contracts the fact that $e_1=e_2$ and $e_1$ is a straight line segment. We conclude that $e_1$ and $e_2$ are both vertical line segments.
\end{proof}

Now we will prove a useful lemma about packed triangles which contain a vertical edge.

\begin{Lemma}\label{lemma:slope}
Let $\strep$ be a representative of a semitoric polygon, and let $B\subset\De$ be an admissibly packed triangle.
If $e:= \overline{B}\setminus B$ is a vertical line segment, then
\begin{itemize}\itemsep=0pt
 \item $B$ does not intersect any of the cuts in $\De$ {\rm (}so $B$ is actually a triangle$)$, and
 \item the two edges of $B$ which are not $e$ each have integer slopes which differ by one.
\end{itemize}
\end{Lemma}

\begin{proof}
Since $e$ is the only edge of $B$ not contained in $\partial \De$, one end point of $e$ will be in the top boundary of $\De$ and one will be in the bottom boundary. Thus, $\partial B \setminus e$ is a path in $\partial \De$ which travels along the top boundary around to the bottom boundary. Since $B$ is convex, this means that $B$ includes all points from the bottom to top of $\De$ in this region, that is:
\[B = \{(x,y)\in\De\colon a \leq x < b\} \qquad \text{or} \qquad B = \{(x,y)\in\De\colon a < x \leq b\}\]
for some $a,b\in\R$ (the expression for $B$ depends on if $B$ is on the right or left side of $\De$).
Thus, if~$B$ intersects any of the cuts in $\De$, then $B$ would also have to contain the corresponding marked point, since it lies on the same vertical line segment in $\De$ as the cut does. Since an admissibly packed triangle in a semitoric polygon does not include any marked points, we conclude that $B$ does not intersect any cuts and we have established the first part of the claim.

Thus, by Definition~\ref{def:STadm-packed-triangle},
$B = M(B(\lambda))$ for some $\lambda>0$ and $M\in\SL$, where $B(\lambda)$ is the convex hull of the origin and the points $(0,\lambda)$, $(\lambda, 0)$, as above.
Note that $M$ sends the hypotenuse of $B(\lambda)$ to the vertical edge of $B$, so either
\[
 M \vectbig{-1}{1} = \vectbig{1}{0} \qquad \text{or}\qquad M \vectbig{1}{-1} = \vectbig{1}{0}.
\]
Without loss of generality, we assume the former.
Along with the fact that $\det(M) = 1$, this implies that $M$ is of the form
\[
M = \begin{pmatrix} 1 & 1\\ c & c+1\end{pmatrix}
\]
for some $c\in \Z$.
The other two edges of $B$ are thus directed by the vectors
\[
M \vectbig{1}{0} = \vectbig{1}{c} \qquad\text{ and }\qquad M \vectbig{0}{1} = \vectbig{1}{c+1},
\]
which have integer slopes $c$ and $c+1$, which differ by one.
\end{proof}

Now we are ready to prove Theorem~\ref{thm:perfect}.

\begin{proof}[Proof of Theorem~\ref{thm:perfect}]
Let $\stpoly$ be a semitoric polygon with at least one marked point which admits a perfect packing $[P]$.
We will show that it must be one of the polygons listed in the statement of the theorem.
For the duration of this proof, we will fix a representative~$\strep$ of $\stpoly$ with a corresponding perfect packing $P$.

{\it Only two triangles.}
First, we will show that $P$ consists of exactly two triangles.
Since $\strep$ contains at least one marked point $c_1 = (x_1,y_2)$, by Lemma~\ref{lem:perfect-local} $c_1$ is contained in the closure of exactly two distinct admissibly packed triangles in $P$, which we call $B_1$ and $B_2$.
Lemma~\ref{lem:perfect-local} also states that $e_1:=\overline{B}_1\setminus B_1$ and $e_2 = \overline{B}_2\setminus B_2$ are vertical line segments, and moreover that $e_1=e_2$.

For the remainder of the proof, let $\ell :=e_1$ denote the vertical line segment through $c_1$ contained in $\De$.
Let $p_{\text{top}}$ and $p_{\text{bottom}}$ be the endpoints of $\ell$, where $p_{\text{top}}$ is on the top boundary of $\De$ and $p_{\text{bottom}}$ is on the bottom boundary.
Furthermore, let
\[
\De_{\text{left}} = \{(x,y)\in\De\colon x<x_1\} \qquad \text{ and } \qquad
\De_{\text{right}} = \{(x,y)\in\De\colon x>x_1\}.
\]

Note that $\De\setminus\ell$ has exactly two components, so since and $B_1$ and $B_2$ are disjoint convex sets which share $\ell$ as a boundary component, this means that each of them lies on one side of $\ell$. That is, switching labels if necessary, $B_1\subseteq \De_{\text{left}}$ and $B_2\subseteq \De_{\text{right}}$.
By the definition of semitoric packing, $\partial B_1 \setminus \ell \subset \partial \De$, which implies that the boundary of $B_1$ encircles $\De_{\text{left}}$, and since $B_1$ is convex, we conclude that $B_1 = \De_{\text{left}}$. Similarly, $B_2 = \De_{\text{right}}$.
Thus, $B_1$ and $B_2$ cover all but a measure zero set of $\De$, so there can be no other sets included in the packing. That is, $P = B_1 \cup B_2$.

Furthermore, since $\De = B_1 \cup \ell \cup B_2$ and marked points do not lie in packed triangles, we conclude that all marked points must lie on $\ell$.

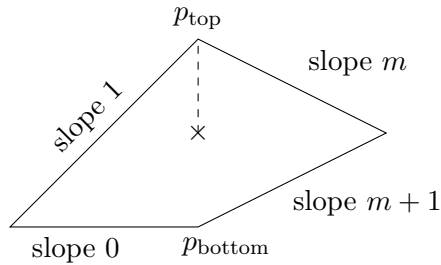
\begin{figure}[t]
 \centering
 \begin{tikzpicture}[scale = 1.25]
 \draw (0,0) -- (2,2) -- (4,1) -- (2,0) -- cycle;
 \node at (2,1) {$\times$};
 \draw [dashed](2,1)--(2,2);
 \node at (1,1) [rotate = 45, anchor = south] {slope $1$};
 \node at (3.7,1.5) [anchor = south] {slope $m$};
 \node at (3.8,.5) [anchor = north] {slope $m+1$};
 \node at (.7,0) [anchor = north] {slope $0$};
 \node at (2,2) [anchor =south] {$p_{\text{top}}$};
 \node at (2.3,0) [anchor =north] {$p_{\text{bottom}}$};
 \end{tikzpicture}
 \caption{In the proof of Theorem~\ref{thm:perfect}, we argue that any semitoric polygon which admits a~perfect packing must be of the above form, where the single ``$\times$'' represents any number $k\geq 1$ of marked points with all cuts upwards. The slopes of the edges are labeled, and the vertex at~$p_{\text{bottom}}$ will instead be a~straight edge if $m=-1$.
 In order for this polygon to be convex with a~vertex at~$p_{\text{top}}$, we see that $m<1$ and $m+1\geq 0$.}
 \label{fig:slopes}
\end{figure}

{\it Considerations of the slopes.}
By Lemma~\ref{lemma:slope}, the slopes of the edges of $B_1$ are integers which differ by one, so by using the actions of $T$ and $T^{-1}$ (which shift all slopes by $1$ or $-1$), we may assume that the edges of $B_1$ which are not $\ell$ have slopes 1 and 0. That is, after translating if necessary, we may assume that $B$ is the convex hull of $(0,0)$, $(a,a)$, and $(a,0)$ with the boundary segment $\ell$ removed, where $a>0$ is some parameter and in these coordinates $\ell$ is the line segment from $(a,0)$ to $(a,a)$.
We may also assume that all cuts in this representative are directed upwards.

Applying Lemma~\ref{lemma:slope} to $B_2$, we conclude that the slopes of the edges of $B_2$ are given by integers~$m$ and $m+1$.
The form of $\De$ is shown in Figure~\ref{fig:slopes}.
Thus, the top boundary of $\De$ is comprised of two segments separated by the point $p_{\text{top}}$ where the cut intersects $\partial \De$.
The segment on the left has slope 1 and the segment on the right has slope $m$.
Since the polygon is convex at $p_{\text{top}}$ we see that $m\leq 1$, and moreover since there must be a vertex where the cut meets the boundary, we actually have that $m<1$.
Similarly, on the bottom boundary the slope on the left is 0 and the slope on the right is $m+1$, and since the polygon must be convex at $p_{\text{bottom}}$ we see that $m+1\geq 0$.
We conclude that $m\in\Z$ satisfies $-1\leq m < 1$, and thus $m\in\{-1,0\}$.

{\it Completing the proof by cases.}
Suppose that the polygon has $k\geq 1$ marked points on the line $\ell$.
Now we complete the proof by considering cases depending on the values of $k$ and $m$.

{\it Case $1$: $m=0$ and $k=1$.}
In this case, $B_2$ is the convex hull of $(a,0)$, $(a,a)$, and $(2a,a)$ with the boundary segment $\ell$ removed. Thus, the polygon is minimal of type inverted (3b) with parameter $n=2$. In this case $p_{\text{top}}$ is a fake corner and $p_{\text{bottom}}$ is a Delzant corner.

{\it Case $2$: $m=0$ and $k\geq 2$.}
This case is impossible, since by acting by the transformation $\mathfrak{t}_\ell^{2}$ to transform two of the upwards cuts into downwards cuts the polygon will become non-convex.

{\it Case $3$: $m=-1$ and $k=1$.}
In this case, $B_2$ is the convex hull of $(a,0)$, $(a,a)$, and~$(0,a)$ with the boundary segment $\ell$ removed. Thus, the polygon is minimal of type (3b) with parameter~${n=2}$.
In this case, $p_{\text{top}}$ is a hidden corner and $p_{\text{bottom}}$ is not a corner.

{\it Case $4$: $m=-1$ and $k=2$.}
In this case, $B_2$ is the convex hull of $(a,0)$, $(a,a)$, and~$(0,a)$ with the boundary segment $\ell$ removed. Thus, the polygon is minimal of type (2) with parameter~$b=0$.
In this case, $p_{\text{top}}$ is a 2-fake corner and $p_{\text{bottom}}$ is not a corner.

{\it Case $5$: $m=-1$ and $k\geq 3$.}
Similar to Case~2 above, this case is impossible since applying~$\mathfrak{t}_\ell^{3}$ to transform three of the upward cuts into downward cuts yields a non-convex polygon.

Having checked all possible cases, we conclude that if $\stpoly$ admits a perfect packing then it must be one of the three types listed in the statement of the theorem. Verifying that these three types admit perfect packings is straightforward and follows from Theorem~\ref{thm:semitoric-examples}.
\end{proof}

\subsection*{Acknowledgements}

This paper is the result of three semesters of work supported by the Illinois Geometry Lab (IGL) program at the University of Illinois at Urbana Champaign.
We are very thankful to the IGL for their support.
We are also thankful to Parth Deshmukh, Felipe Pallo Rivadeneira, Haoxiang Sun, Deming Tian, and Tongshu Liu
for their hard work in
earlier IGL projects tangentially related to this one.
We also thank Marino Romero for encouraging us to compile our results into this paper, and Yohann Le Floch for many helpful comments and suggestions on an early version.
Finally, we also thank the referees who gave very helpful comments which improved the paper.\looseness=-1

\pdfbookmark[1]{References}{ref}
\LastPageEnding

\end{document}